\documentclass[11pt,reqno]{article}

\usepackage{newtxtext} 

\usepackage[T1]{fontenc}
\usepackage[english]{babel}
\usepackage[left=1in,right=1in,top=1in,bottom=1in]{geometry}
\usepackage{mathtools}
\usepackage{amsthm}
\usepackage[mathcal]{euscript} 
\usepackage{amssymb}
\usepackage{amscd}
\usepackage{amsmath}
\usepackage{mathrsfs}
\usepackage{amsbsy}
\usepackage{array}
\usepackage{cancel}
\usepackage{enumitem}
\usepackage{tikz-cd}
\usepackage{tikz}
\usetikzlibrary{matrix,positioning,arrows,fit,backgrounds}
\usepackage{url}
\usepackage{hyperref}
\usepackage{cite}
\usepackage{setspace}
\usepackage{extpfeil} 
\usepackage{fancyhdr}    

\newcommand{\headerpalatino}{\fontfamily{ppl}\selectfont}
\let\origmathcal\mathcal
\usepackage[charter]{mathdesign} 
\let\mathcal\origmathcal

\usetikzlibrary{matrix, positioning, fit}
\usetikzlibrary{arrows.meta, positioning}
\usetikzlibrary{arrows.meta, positioning, decorations.pathreplacing}
\binoppenalty=\maxdimen
\relpenalty=\maxdimen
\allowdisplaybreaks[1]

\title{Positivity of automorphic vector bundles on unitary Shimura varieties}
\author{Deding Yang} 
\date{}
\hypersetup{
	colorlinks,
	allcolors=[rgb]{0,0.45,0}
}

\renewcommand{\H}{H_1^{\mathrm{dR}}}

\DeclareMathOperator{\SHom}{\mathcal{H}\text{\kern -3pt {\calligra\large om}}\,}

\theoremstyle{plain}
\newtheorem{thm}{Theorem}[section]
\newtheorem{coro}[thm]{Corollary}
\newtheorem{prop}[thm]{Proposition}
\newtheorem{lem}[thm]{Lemma}
\newtheorem{conjecture}[thm]{Conjecture}

\theoremstyle{definition}
\newtheorem{definition}[thm]{Definition}
\newtheorem{rmk}[thm]{Remark}

\newtheorem{example}[thm]{Example}
\newtheorem{convention}[thm]{Convention}
\newtheorem{construction}[thm]{Construction}

\usepackage{titlesec}
\titleformat{\section}
{\centering\bfseries\Large} 
{\thesection}{1em}{}

\usepackage{tocloft}

\hypersetup{
	linktoc=page 
}

\usepackage{titlesec}
\titleformat{\subsection}
{\normalfont\bfseries}        
{\thesubsection}                            
{0pt}                         
{}                            
\titlespacing*{\subsection}       
{0pt}                         
{5pt}                         
{0pt}                         

\titleformat{\subsubsection}
{\normalfont\bfseries}        
{\thesubsubsection}                            
{0pt}                         
{}                            
\titlespacing*{\subsubsection}       
{0pt}                         
{5pt}                         
{0pt}                         

\begin{document}
	\maketitle

	\section*{Abstract}
	Let $X$ be the special fiber of a unitary Shimura variety of hyperspecial level at a prime $p$ inert in the totally real field $F$. Let $Y\to X$ be the associated flag space. For every $L$-dominant weight $\lambda$, let $\mathcal{L}_Y(\lambda)$ denote the corresponding automorphic line bundle. We give an explicit necessary and sufficient criterion, in terms of the signature data and the coordinates of $\lambda$, for the ampleness of $\mathcal{L}_Y(\lambda)$. 
    The criterion generalizes the known ample cone for Hilbert modular and $U(2)$-Shimura varieties. The proof develops the machinery of the description of certain Ekedahl--Oort strata, a geometric Jacquet--Langlands correspondence between strata of unitary Shimura varieties with different signatures, and the construction of stratum Hasse invariants, and introduced a way to systematically deal with combinatorical data in the higher rank case.
	
	
	\addcontentsline{toc}{section}{Abstract}

	\medskip
	{\color{black} \tableofcontents}

	\medskip
	
	\section{Introduction}
	
	
	

\subsection{ Overview}

    Shimura varieties are higher dimensional generalizations of modular curves. The study of the cohomology of certain local systems or vector bundles on the Shimura variety relates Galois representations and automorphic forms, known as a typical approach towards the Langlands program. The past few decades have witnessed substantial developments in this theory (See, for example, \cite{CaraianiShin,HarrisSurvey} for some modern survey). The main focus of this paper is on the coherent cohomology of automorphic vector bundles. Knowledge of vanishing theorems of coherent cohomology is an important input for modularity lifting theorems (see \cite{HarrisTWMethod1,AtanasovHarrisTWMethod2} for the Taylor--Wiles method case; see also \cite{CalegariGeraghty} for beyond the Taylor--Wiles case; see also \cite{Lan-SuhGeneralVanishing,Lan-SuhCompactVanishing} for their vanishing theorems). Therefore, we believe that it is natural to ask for an ampleness criterion for automorphic vector bundles. By Kodaira type vanishing results, this directly implies vanishing theorems of coherent cohomology.

    On the other hand, the study of positivity for algebraic varieties is itself of great importance, especially in the characteristic-$p$ case. The existence of Frobenius structures leads to the emergence of many interesting —- and sometimes peculiar —- subvarieties. In the context of PEL-type Shimura varieties, one could use the moduli interpretation to obtain a very concrete description of these subvarieties. The study of these subvarieties leads us to establish a geometric Jacquet--Langlands correspondence in characteristic $p$. This correspondence helps us understand the intersection numbers over the special fiber, which finally result in a \emph{complete} answer of the ampleness criterion of automorphic line bundles.

    In the following, we first explain our main result. After that, we explain the strategy of our proof and the relations to arithmetic applications.

    Let $(G,X)$ be a Shimura datum. Upon fixing a neat compact open subgroup $K\subseteq G(\mathbb{A}_f)$, we obtain an 
	algebraic variety $Sh_K(G,X)$ that is defined over the reflex field $E=E(G,X)$:
	\begin{equation}
		Sh_K(G,X)(\mathbb{C})=G(\mathbb{Q})\backslash X\times(G(\mathbb{A}_f)/K).
	\end{equation}
	We further assume that $G_{\mathbb{Q}_p}$ is unramified, i.e., it is quasi-split and split over an unramified extension of $\mathbb{Q}_p$. We assume that $K\subseteq G(\mathbb{A}_f)$ is of the form $K^pK_p$, where $K^p\subseteq G(\mathbb{A}_f^{(p)})$ and $K_p\subseteq G(\mathbb{Q}_p)$ is a hyperspecial subgroup. We have the following results of \cite{KisinIntegralModel,VasiuIntegralModel}:
	
	\begin{thm} If $(G,X)$ is of abelian type, then $Sh_K(G,X)$ admits a smooth canonical integral model $\mathscr{S}_K(G,X)$ over $\mathcal{O}_{E,(p)}$.
	\end{thm}
	
	Let $S_{K}(G)=\mathscr{S}_K(G,X)_k$ be the geometric special fiber of the integral model. Here, $k=\overline{\mathbb{F}}_p$ and we omit $X$ from the notation for simplicity. For any $L$-dominant weight $\lambda$, where $L$ is the Levi subgroup of the parabolic subgroup of $G$ that corresponds to the Hodge filtration, there is an \emph{automorphic vector bundle} $\mathcal{V}(\lambda)$ on $Sh_K(G,X)$. The automorphic vector bundle $\mathcal{V}(\lambda)$ extends to $\mathscr{S}_K(G)$, and we still denote its pullback to the special fiber by $\mathcal{V}(\lambda)$. The purpose of this paper is to study the positivity of $\mathcal{V}(\lambda)$. To illustrate this, we first consider an example.
	
	\medskip\noindent\textbf{Example:} Let $F/\mathbb{Q}$ be a totally real field in which $p$ is inert. Let $\{\tau_1,\dots,\tau_f\}$ be the set of all real embeddings of $F$. If we fix an isomorphism $\bar{\mathbb{Q}}_p\simeq \mathbb{C}$, we can view the $\tau_i$'s as $p$-adic embeddings, and may order them such that $\sigma^{-1}\tau_i=\tau_{i+1}$. Here we view the subscripts as modulo $f$, and $\sigma$ is the Frobenius endomorphism. Let $G=\textnormal{Res}_{F/\mathbb{Q}}GL_{2,F}$ and $K_p=GL_2(\mathcal{O}_F)$. For a sufficiently small tame level $K^p\subseteq G(\mathbb{A}_f^{(p)})$, the corresponding algebraic variety $Sh_K(G)$ is known as the Hilbert modular variety. Its integral model $\mathscr{S}_K(G)/\mathcal{O}_F$ is the moduli space of abelian schemes with a certain polarization, an action of $\mathcal{O}_F$, and a level structure. (This is not entirely accurate. One needs to twist the moduli problem. See \cite{Tian-Xiao} for a detailed explanation.) 
	
	The universal abelian scheme carries an $\mathcal{O}_F$-action, and thus the Hodge filtration
	\begin{equation}
		0\longrightarrow \omega_{\mathcal{A}^\vee/S_K(G)}\longrightarrow\H(\mathcal{A}/S_K(G))\longrightarrow \textnormal{Lie}_{\mathcal{A}/S_K(G)}\longrightarrow 0
	\end{equation}
	decomposes into $\tau_i$-components for $1\le i\le f$. We denote by $\omega_i$ the subsheaf of $\omega_{\mathcal{A}^\vee/S_K(G)}$ on which $\mathcal{O}_F$ acts via $\tau_i$. In fact, $\omega_i$ is a line bundle. In the Hilbert case, given a tuple of integers $(k_i,l_i)_{1\le i\le f}$ such that $k_i+2l_i$ equals a constant $c$ for all $i$, the corresponding automorphic line bundle of weight $\lambda=(k_i,l_i)_i$ is given by:
	\begin{equation}
		\mathcal{L}(\lambda)=\bigotimes_{i=1}^f \big((\textnormal{det }\H(\mathcal{A}/S_K(G))_i)^{\otimes l_i}\otimes \omega_i^{\otimes k_i}\big).
	\end{equation}
	In \cite{YangAmpleness}, the author proved the following criterion (We omit the issue of compactification for simplicity)
	
	\begin{thm}\label{Ampleness of Hilbert Modular Variety}
		$\mathcal{L}(\lambda)$ is ample if and only if $pk_i> k_{i+1}$ for $1\le i\le f$.
	\end{thm}
	
	\begin{thm}\label{Nefness of Hilbert Modular Variety}
		$\mathcal{L}(\lambda)$ is nef if and only if $pk_i\ge k_{i+1}$ for $1\le i\le f$.
	\end{thm}
	
	As a corollary, we obtain a Kodaira type vanishing result in characteristic $p$:
	\begin{coro}
		Let $S_K(G)^{\textnormal{tor}}$ be the toroidal compactification (for a given smooth admissible polyhedral cone decomposition, as in \cite{LanCompactification}) of $S_K(G)$ and $D=S_K(G)^{\textnormal{tor}}\backslash S_K(G)$ be the boundary divisor. If $\lambda$ is an ample weight, then
		\begin{equation}
			H^i(S_K(G)^{\textnormal{tor}},\Omega^j(-D)\otimes \mathcal{L}(\lambda))=0, \qquad\textnormal{for $i+j>\textnormal{sup}\{f,2f-p\}$}.
		\end{equation}
	\end{coro}

	\medskip
	However, beyond the Hilbert case, $\mathcal{V}(\lambda)$ are usually vector bundles. It is less common to directly discuss the ampleness of $\mathcal{V}(\lambda)$. Instead, we consider a Borel--Weil construction (see, for example, \cite{Fulton-HarrisBook}). Roughly speaking, we define $Y$ as the moduli space over $k$ that parametrizes tuples $(s,\textnormal{Fil}^\bullet)$, where
	
	$\bullet$ $s$ is a point of $S_K(G)$;
	
	$\bullet$ $\textnormal{Fil}^\bullet$ is a full refinement of the Hodge filtration with respect to additional $G$-structures.
	
	The theory of $G$-zips and $G$-zipflags gives a more natural and more rigorous group-theoretic interpretation of the above construction. See Section 3 for an example. Admitting this moduli interpretation, it follows that the fibers of the natural forgetful map $Y\to S_K(G)$ are isomorphic to the flag variety $P/B$ over $k$. We refer to $Y$ as the flag space in this paper. 
	
	For any $\lambda\in X^\ast_{+,L}$, one can construct a natural line bundle $\mathcal{L}_Y(\lambda)$ over $Y$, which we will refer to as the \emph{automorphic line bundle} of weight $\lambda$ in this paper. 
	
	The following result, which is originally due to Borel--Weil and Bott (see \cite{Fulton-HarrisBook}), explains the relation between $\mathcal{L}_Y(\lambda)$ and $\mathcal{L}_X(\lambda)$.
	\begin{prop}
		Let $\pi:Y\to S_K(G)$ be the natural projection. For any $\lambda\in X^\ast_{+,L}$,
		\begin{equation}
			\pi_{\ast}\mathcal{L}_Y(\lambda)=\mathcal{V}(\lambda)
		\end{equation}
		and
		\begin{equation}
			H^i(Y,\mathcal{L}_Y(\lambda))=H^i(S_K(G),\mathcal{V}(\lambda)).
		\end{equation}
	\end{prop}
	
	Therefore, for our purpose, we should ask the following question instead.
	
	\medskip
	\noindent\textbf{Question: }When is the line bundle $\mathcal{L}_Y(\lambda)$ ample?
	
	\medskip
	We first review the history of this problem.
    \begin{enumerate}[itemsep=0pt,topsep=2pt, parsep=0pt]
        \item In \cite{Andreatta-Goren}, Andreatta and Goren computed the ample cone of Hilbert modular surface ($f=2$ case of Theorem \ref{Ampleness of Hilbert Modular Variety}), by constructing partial Hasse invariants and studying intersection theory on the surface directly.
        \item In \cite{Tian-Xiao}, Tian and Xiao conjectured the ampleness criterion for $U(2)$/Hilbert/quaternionic Shimura varieties. They studied the geometric description of Goren--Oort strata and proved the necessity part of Theorem \ref{Ampleness of Hilbert Modular Variety}.
        \item In \cite{YangAmpleness}, the author proved the ampleness criterion conjectured by \cite{Tian-Xiao}, by computing the intersection numbers using the description of strata.
        \item In the unpublished paper \cite{BGKS}, Brunebarbe--Goldring--Koskivirta--Stroh proposed a sufficiency criterion for ampleness for Hodge type Shimura varieties. However, their bound is known to be not sharp (even in the Hilbert case).
    \end{enumerate}

	Little is known about the ampleness criterion beyond the above cases, especially when the group $G$ is not of type $A_1$. The main difficulty is that the stratification of the flag space is usually very complicated, and computing intersection numbers $(\mathcal{L}_Y(\lambda)\cdot C)$ with a complete curve $C\subseteq Y$ is subtle.
	
	In this paper, we give a complete answer to the above question when $G_{\mathbb{Q}_p}^{\textnormal{sc}}$ is the unramified Weil restriction of $SL_n$.

	\subsection{ Statement of the main results.}
	
	Let $F$ be a totally real field of degree $N$ over $\mathbb{Q}$ in which $p$ is inert. We label the set of real embeddings of $F$ by $\{\tau_1,\dots,\tau_N\}$ as before. Let $E/F$ be an imaginary quadratic extension in which $p\mathcal{O}_E=\mathfrak{q}\mathfrak{q}^c$. Each $\tau_i$ extends to complex embeddings $\tilde\tau_i$ (over $\mathfrak{q}$) and $\tilde\tau_i^c$ (over $\mathfrak{q}^c$). For a tuple of integers $(m_i,n_i)_{1\le i\le N}$ such that $m_i+n_i=n$, let $G$ be the unitary similitude group over $\mathbb{Q}$, with signatures at infinite places given by $(m_i,n_i)_{1\le i\le N}$, i.e.,
	\begin{equation}
		G(\mathbb{R})=G\big(\prod_{i=1}^N U(m_i,n_i)\big).
	\end{equation}
	Let $X$ be the special fiber of the corresponding Shimura variety of signature $(m_i,n_i)$. We apologize for using the letter $X$ again (but for a different meaning) to maintain consistency with the rest of the paper. The signature condition implies that, in the Hodge filtration,
	\begin{equation}
		0\longrightarrow \omega_{\mathcal{A}^\vee/X,\tilde\tau_i}\longrightarrow\H(\mathcal{A}/X)_{\tilde\tau_i}\longrightarrow \textnormal{Lie}_{\mathcal{A}/X,\tilde\tau_i}\longrightarrow 0,
	\end{equation}
	$\omega_{\mathcal{A}^\vee/X,\tilde\tau_i}$ (resp. $\omega_{\mathcal{A}^\vee/X,\tilde\tau_{i}^c}$) is a subbundle of $\omega_{\mathcal{A}^\vee/X}$ of rank $m_i$ (resp. $n_i$). The polarization in the moduli problem induces a perfect pairing
	\begin{equation}
		\langle\cdot,\cdot\rangle:\H(\mathcal{A}/X)_{\tilde\tau_i}\times \H(\mathcal{A}/X)_{\tilde\tau_i^c}\longrightarrow \mathcal{O}_X
	\end{equation}
	such that $\omega_{\mathcal{A}^\vee/X,\tilde\tau_i^c}=\omega_{\mathcal{A}^\vee/X,\tilde\tau_i}^\perp$ under this pairing.
	
	We define the flag space $Y$ as the moduli space over $k$ whose $S$-valued points are $(s,\textnormal{Fil}^\bullet)$, where $s\in X(S)$, and $\textnormal{Fil}=\{\textnormal{Fil}_i,\textnormal{Fil}_i^c\}$ is a collection of filtrations at each place:
	\begin{equation}
		\textnormal{Fil}_i: 0=\mathcal{F}^i_0\subseteq \mathcal{F}^i_1\subseteq \cdots\subseteq \mathcal{F}^i_{m_i-1}\subseteq \mathcal{F}^i_{m_i}=\omega_{A^\vee/S,\tilde\tau_i}\subseteq \mathcal{F}^i_{m_i+1}\subseteq\cdots\subseteq\mathcal{F}^i_{n-1}\subseteq \mathcal{F}^i_n=\H(A/S)_{\tau_i}.
	\end{equation}
	Similarly for the $\tilde\tau_i^c$-part. In addition, we require that $\mathcal{F}^{i^c}_{u}=\mathcal{F}^{i,\perp}_{n-u}$ for $1\le i\le N$ and $0\le u\le n$. Here, $\mathcal{F}^{i^c}_{u}$ is a vector bundle in the $\tilde\tau_{i}^c$-component. Since there are dual relations between the $\tilde\tau_i$-and the $\tilde\tau_i^c$-parts, we can simplify the data by only recording the refinement of the Hodge filtration at all the $\tilde\tau_i$-components.
	
	In general, irreducible automorphic vector bundles are parametrized by a character in $X^\ast_{+,L}$, i.e., a character of $T$ that is dominant with respect to $L$. For our specific $G$, a character $\lambda\in X^\ast_{+,L}$ corresponds to a tuple of integers $\{k^i_j\}^{1\le i\le N}_{1\le j\le n}$. The automorphic line bundle of weight $\lambda$ is 
	\begin{equation}
		\mathcal{L}_Y(\lambda)=\bigotimes_{i=1}^{N}\bigotimes_{j=1}^n \big(\mathcal{F}^i_{j}/\mathcal{F}^i_{j-1}\big)^{\otimes k^i_j}.
	\end{equation}
	The main theorem of this paper is

	\begin{thm}\label{Main Theorem - Introduction}[Theorem \ref{The ampleness criterion for the flag space, case 1}, \ref{The ampleness criterion for the flag space, case 2}, \ref{The ampleness criterion for the flag space, case 2'}, \ref{The ampleness criterion for the flag space, case 3}]
		Let $T=\{i|m_i\notin \{0,n\}\}$ and $t=\# T$.
		\begin{enumerate}
			
			\item If $t=0$, then $X$ is a Shimura set, and $Y$ is a disjoint union of flag varieties in characteristic $p$. $\mathcal{L}_Y(\lambda)$ is ample if and only if
			\begin{equation}
				k^i_n>k^i_{n-1}>\cdots>k^i_1,\qquad\textnormal{for $1\le i\le N$}.
			\end{equation}
			
			\item If $t=1$, we assume $T=\{1\}$ for simplicity. If moreover $m_1=n-1$, then $\mathcal{L}_Y(\lambda)$ is ample if and only if
			\begin{equation}
				\begin{aligned}
					&k^1_{n-1}>k^1_{n-2}>\cdots>k^1_1,\\
					&k^i_{n}>k^i_{n-1}>\cdots>k^i_1,\qquad\textnormal{for $2\le i\le N$},\\
					&\big(\sum_{u=0}^{n-2}p^{Nu}\big)\sum_{j=0}^{N-1}\frac{1}{p^j}(k_1^{1+j}-k^{1+j}_n)>\sum_{r=2}^{n-1}p^{(r-2)N}\sum_{j=0}^{N-1}\frac{1}{p^j}(k^{1+j}_r-k^{1+j}_n). 
				\end{aligned}
			\end{equation}
			
			\item If $t=1$, we assume $T=\{1\}$ for simplicity. If moreover $m_1=1$, then $\mathcal{L}_Y(\lambda)$ is ample if and only if
			\begin{equation}
				\begin{aligned}
					&k^1_n>k^1_{n-1}>\cdots >k^1_2,\\
					&k^i_n>k^i_{n-1}>\cdots>k^i_1,\qquad\qquad \textnormal{for $i\ne 1$},\\
					&	\big(\sum_{u=0}^{n-2}p^{Nu}\big)\sum_{j=0}^{N-1}\frac{1}{p^j}(k_1^{1+j}-k^{1+j}_n)>\sum_{r=2}^{n-1}p^{(n-1-r)N}\sum_{j=0}^{N-1}\frac{1}{p^j}(k^{1+j}_1-k^{1+j}_r).
				\end{aligned}
			\end{equation}
			
			\item In all other cases, let $T=\{i_1,\dots,i_t\}$ with $i_1<i_2<\cdots<i_t$. Let $a(i_l)=i_{l+1}-i_l$ for $1\le l\le t-1$ and $a(i_t)=i_1+N-i_t$. Then $\mathcal{L}_Y(\lambda)$ is ample if and only if
			\begin{equation}
				\begin{aligned}
					&k^{i_l}_{m_{i_l}}>k^{i_l}_{m_{i_l}-1}>\cdots>k^{i_l}_{1},\quad k^{i_l}_{n}>k^{i_l}_{n-1}>\cdots>k^{i_l}_{m_{i_l}+1}, \qquad\textnormal{for all $1\le l\le t$},\\
					&k^j_n>k^j_{n-1}>\cdots>k^j_1,\qquad\textnormal{for $j\notin T$},\\
					&p^{a(i_l)}\sum_{j=0}^{a(i_l)-1}\frac{1}{p^j}(k^{i_l+j}_1-k^{i_1+j}_n)>k^{i_{l+1}}_{m_{i_{l+1}}}-k^{i_{l+1}}_{m_{i_{l+1}}+1}\qquad \textnormal{for all $1\le l\le t$.}
				\end{aligned}
			\end{equation}
			Moreover, $\mathcal{L}_Y(\lambda)$ is nef if and only if we substitute all the ">" by $\ge$ in the above inequalities.
		\end{enumerate}
	\end{thm}
	
    Note that the formulas are uniform in general, except the case $t=1$.
    
	We remark that the notion of flag space can be generalized to partial flag spaces, by recording a partial refinement of the Hodge filtration. There are analogous theorems for partial flag spaces over $X$. In fact, the proof of Theorem \ref{Main Theorem - Introduction} is essentially assembling these results over partial flag spaces (See Theorem \ref{The corresponding nefness criterion}, \ref{The ampleness criterion for Shimura variety}, \ref{The nefness criterion for partial flag space} in Section 3).
	
	Before sketching the proof of Theorem \ref{Main Theorem - Introduction}, we first discuss some applications of this ampleness criterion to other topics.

	\subsection{ Applications.}

	\subsubsection{ Vanishing of coherent cohomology.}	
	
	As we mentioned earlier, the study of different cohomology theories plays an important role in the world of Shimura varieties. There have been studies \cite{Caraiani-ScholzeCompactGenericPart,Caraiani-ScholzeGenericCohomology,Caraiani-Tamiozzo,Hamann-LeeVanishing}, etc, to control the (concentration of) \'etale cohomology of Shimura varieties, \cite{AlexandreVanishing,Lan-SuhCompactVanishing,Lan-SuhGeneralVanishing}, etc, to establish the vanishing results of coherent cohomology. Our ampleness criterion directly yields the Kodaira-type vanishing results.
	
	\begin{thm}
		Set $d=\textnormal{dim $Y$}$ and let $\mathcal{L}_Y(\lambda)$ be an ample line bundle over $Y$. Then
		\begin{equation}
			H^i(Y,\Omega^j_{\mathcal{A}/Y}\otimes\mathcal{L}_Y(\lambda))=0, \qquad \textnormal{for $i+j>\textnormal{max$\{d,2d-p\}$}$}.
		\end{equation}
	\end{thm}
	
	We can also establish new vanishing results on $X$. When $\lambda$ is a character of $L\subseteq P$, the irreducible representation of highest weight $\lambda$ is 1-dimensional. Consequently, the automorphic vector bundle $\mathcal{V}(\lambda)$ over $X$ is locally free of rank one. We sometimes refer to it as $\mathcal{L}_X(\lambda)$. In our unitary setup, such $\lambda$ corresponds to a tuple of integers $\{k_i\}_{i=1}^N$, and
	\begin{equation}
		\mathcal{L}_X(\lambda)=\mathcal{L}_X(\{k_i\})=\bigotimes_{i=1}^N \big(\textnormal{det }\omega_{\mathcal{A}^\vee/X,\tau_i}\big)^{\otimes k_i}.
	\end{equation}
	
	The proof of Theorem \ref{Main Theorem - Introduction} will also imply the following criterion for $\mathcal{L}_X(\lambda)$:
	
	\begin{thm}\label{Introduction-Vanishing of line bundles on X}[Theorem \ref{The ampleness criterion for Shimura variety}]
		If $T=\{i_1,\dots,i_t\}$. The line bundle $\mathcal{L}_X(\{k_i\})$ is ample if and only if
		\begin{equation}
			p^{a(i_l)}k_{i_l}>k_{i_{l+1}}, \qquad \textnormal{for $1\le l\le t$}.
		\end{equation}
	\end{thm}	
	
	Recall the Kodaira type vanishing result in \cite{Lan-SuhCompactVanishing}:
	\begin{thm}
		Suppose that $\mathcal{L}_X(\{k_i\})$ is ample, and $\mu\in X^{+,<_{\textnormal{re}}p}_{G}$. Then for any $w\in W^M$, we have
		\begin{equation}
			H^{i-l(w)}(X,\mathcal{L}_X(\lambda)\otimes \mathcal{V}^{\vee}({w\cdot \mu}))=0,\qquad \textnormal{for $i>d$}.
		\end{equation}
	\end{thm}
	
	We omit the detailed explanation of the notations here as they are beyond the scope of this paper. In Lan and Suh's original paper, they applied this theorem together with the ampleness of the Hodge line bundle $\textnormal{det }\omega$ (which corresponds to $\mathcal{L}_X(\{1,\dots,1\})$) to deduce their vanishing results. We could use the ampleness criterion in this paper to improve this vanishing result for unitary Shimura varieties.

    \subsubsection{ Applications to mod $p$ automorphic forms.}
    
	A direct application of the vanishing results of coherent cohomology is to understand the congruences of $p$-adic automorphic forms and liftability of mod $p$ automorphic forms, or more generally, cohomological classes. This approach has been used in \cite{Emerton-Reduzzi-Xiao,GeorgeBoxerThesis,GK-strata}, to construct Galois representations attached to cohomology classes. 
    
    In our joint work in progress \cite{YangYang}, we prove the following conjecture by Diamond and Sasaki in \cite{Diamond-Sasaki}.
    
    \begin{conjecture}
        Let $\overline\rho:G_F\to GL_2(\overline{\mathbb F}_p)$ be an irreducible mod $p$ Galois representation. Then $\overline{\rho}$ is \emph{algebraically modular} of weight $(\underline{k},\underline{l})$ if and only if $\overline{\rho}$ is \emph{geometrically modular} of weight $(\underline{k},\underline{l})$.
    \end{conjecture}
    
    For the definitions of two notions of modularity, see \cite{Buzzard-Diamond-Jarvis} and \cite{Diamond-Sasaki}. In fact, we generalize the statement to quaternionic Shimura varieties and prove the conjecture under some mild Taylor--Wiles hypotheses. The ampleness criterion is a first input to deduce the liftability of mod $p$ automorphic forms of certain weights. We also apply the geometric Jacquet--Langlands correspondence developed in this paper to relate automorphic forms for different quaternionic Shimura varieties. This convinces us that the ampleness criterion in this paper can provide a potential method in attacking the modularity of mod $p$ representations beyond the $GL_2$ case.


	\subsection{Ingredients of the proof.}
	\subsubsection{ Geometry of the special fiber and a geometric Jacquet--Langlands correspondence.}
	
	Let $X$ denote the special fiber of a Hodge type Shimura variety of hyperspecial level at $p$. There exists a structure morphism
	\begin{equation}
		\zeta:X\longrightarrow G\textnormal{-zip}^\mu.
	\end{equation}	
	The target of this morphism is called the stack of $G$-zips. This is a group-theoretical object first introduced by Moonen and Wedhorn in \cite{MoonenGroupsSchemes,Moonen-Wedhorn} for $G=GL_n$, and systematically studied by Pink, Wedhorn and Ziegler for a general reductive group $G$ in \cite{Pink-Wedhorn-ZieglerAlgebraicZipData,Pink-Wedhorn-ZieglerF-zips}. 
    The stack of $G$-zips captures the structure of the mod $p$ de Rham cohomology of the universal abelian scheme and its Frobenius structure.
    
	The stack of $G$-zips admits a stratification 
	$$G\textnormal{-zip}^\mu=\coprod_{w\in {^I}W} G\textnormal{-zip}^\mu_w.$$
	Zhang showed that the morphism $\zeta$ is smooth in \cite{ZhangEOStrata}, so we may pullback the stratification to X. The resulting stratification on $X$ is called the \emph{Ekedahl--Oort stratification}.
	
	There have been extensive studies on the geometry of Ekedahl--Oort stratification (See \cite{Goldring-Imai-Koskivirta,GK-stratification,KoskivirtaNormalization}, \cite{GK-strata,Viehmann-WedhornEONewton,Wedhorn-ZieglerTautological}, etc). We note that while the Ekedahl--Oort strata are smooth themselves, their closures are generally highly singular. We will also study the normalization of some Ekedahl--Oort strata. 
	
	The crucial geometric aspect of our proof, namely, a geometric mod $p$ Jacquet--Langlands correspondence, deserves discussion. In \cite{Helm}, Helm studied the sparse Goren--Oort strata of $U(2)$-Shimura varieties and proved that each sparse stratum is a $\mathbb{P}^1$-power bundle over another $U(2)$ Shimura variety with different signature. Subsequent works of \cite{Tian-Xiao,Helm-Tian-Xiao} employed Helm's trick to provide a comprehensive description of all Goren--Oort strata in Hilbert/quaternionic/$U(2)$ Shimura varieties. The work in \cite{Xiao-Zhu} significantly extended this description of strata to the basic locus of Hodge type Shimura varieties. They established a parametrization of the basic locus by another Shimura variety with different signature, where the fibers are isomorphic to iterated Deligne--Lusztig varieties.
	
	In this paper, rather than using this $\mathbb{P}^1$-bundle description, we conceptualize this phenomenon as a correspondence between the (partial) flag spaces of \emph{different} unitary Shimura varieties. The main theorem for this step is
	\begin{prop}[Theorem \ref{Description of strata}]
		Let $B\subseteq P_0\subseteq P$ be an intermediate parabolic subgroup such that the associated partial flag space $Y_{P_0}$ is characterized by recording precisely one subbundle $\mathcal{E}_{r_i}^i\subseteq \omega_i$ of rank $r_i$ for each $i$. Let $Z$ denote the closed subscheme of $Y_{P_0}$ cut out by the condition
		\begin{equation}
			\mathcal{E}^i_{r_i}\subseteq V^{-1}(\mathcal{E}_{r_{i+1}}^{i+1,(p)}),\qquad F(\mathcal{E}_{r_{i+1}}^{i+1,(p)})\subseteq \mathcal{E}^i_{r_i}
		\end{equation}
		for all $i$. Then $Z$ is, up to some Frobenius factor, isomorphic to a stratum $Z'$ in some partial flag space $Y'_{P_0'}$ over $X'$, where $X'$ is the special fiber of a unitary Shimura variety with hyperspecial level structure and signature condition
		\begin{equation}
			m_i'=m_{i}-r_{i}+r_{i-1},\qquad \textnormal{for $1\le i\le N$}.
		\end{equation}
		Furthermore, this correspondence respects extra structures.
	\end{prop}
	This correspondence can be understood through the following diagram
	
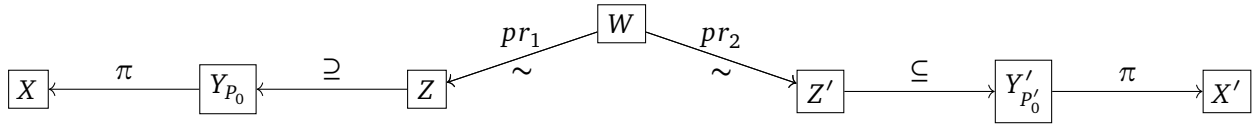
\begin{center}
		\begin{tikzpicture}[node distance=10pt and 1cm, auto, column sep=small]
			\node(1) [draw, align=center] {$W$};
			\node(2) [draw, align=center, below left=of 1] {$Z$};
			\node(3) [draw, align=center, below right=of 1] {$Z'$};
			\node(4) [draw, align=center, left=of 2] {$Y_{P_0}$};
			\node(5) [draw, align=center, right=of 3] {$Y'_{P_0'}$};
			\node(6) [draw, align=center, left=of 4] {$X$};
			\node(7) [draw, align=center, right=of 5] {$X'$};
			
			\draw [->] (1) -- node[midway, above] {$pr_1$} (2);
			\draw [->] (1) -- node[midway, below] {$\sim$} (2);
			\draw [->] (1) -- node[midway, above] {$pr_2$} (3);
			\draw [->] (1) -- node[midway, below] {$\sim$} (3);
			\draw [<-] (4) -- node[midway, above] {$\supseteq$} (2);
			\draw [->] (3) -- node[midway, above] {$\subseteq$} (5);
			\draw [->] (5) -- node[midway, above] {$\pi$} (7);
			\draw [->] (4) -- node[midway, above] {$\pi$} (6);
		\end{tikzpicture}
	\end{center}
	Here $X$ and $X'$ are unitary Shimura varieties with different signatures, $Y_{P_0}$ and $Y_{P_0'}'$ are partial flag spaces over $X$ and $X'$, $Z$ and $Z'$ are strata within $Y_{P_0}$ and $Y_{P_0'}'$ respectively. The symbol $\sim$ stands for an isomorphism up to Frobenius twists. This diagram allows us to relate automorphic line bundles over \emph{different} unitary Shimura varieties via pullback, allowing an inductive proof on possible signatures of the unitary group.

    A geometric Jacquet--Langlands correspondence of such kind has its application in various arithmetic problems: In \cite{Helm-Tian-Xiao,Xiao-Zhu} for the construction of Tate cycles; In \cite{Caraiani-Tamiozzo} to avoid trace formula computations; In our work in progress \cite{YangYang} to transform modularity of mod $p$ Galois representations.
    

	\subsubsection{ Strata Hasse invariants.}
	
	The classical Hasse invariant $h$ is a modular form of weight $p-1$ in characteristic $p$, induced by the Verschiebung map on $\omega$. It has $q$-expansion equal to 1. The zero locus of $h$ is the supersingular locus.
	
	It is natural to consider generalizations of the Hasse invariant to higher-dimensional Shimura varieties. By a \emph{strata Hasse invariant}, we mean a section $s\in H^0(\overline{X_w},\mathcal{V}(\lambda))$ for some $w\in {^I}W$ and $\lambda\in X^\ast_{+,L}$, such that $Z(s)$ is a union of strata. In \cite{Andreatta-Goren}, partial Hasse invariants were constructed for Hilbert modular varieties. In recent years, various generalizations have been developed (from different perspectives, though those related to the theory of $G$-zips are most relevant to our paper), see, for example, \cite{GeorgeBoxerThesis,Reduzzi-XiaoPartialHasse,Koskivirta-Wedhorn,IK-partial,Goldring-Imai-Koskivirta,GK-strata}. In the work of \cite{Emerton-Reduzzi-Xiao, GeorgeBoxerThesis, GK-strata}, they constructed strata Hasse invariants and used them to study congruences of automorphic forms.
	
	
	To prove the ampleness criterion, our naive idea is to express the class $[\mathcal{L}_Y(\lambda)]$ as a non-negative linear combination of strata Hasse invariants. Then the intersection number $(\mathcal{L}_Y(\lambda)\cdot C)$ would be non-negative as long as the test curve $C$ does not lie in any lower-dimensional strata. Then we can reduce the problem to a lower dimensional variety and continue this process. To avoid strata with bad geometric properties, we need a more careful construction of strata Hasse invariants whose zero locus does not contain these specific strata. This is a generalization of the previous work, which relies on the moduli interpretation.  

	\subsubsection{ Outline of the proof and structure of the paper.}
	
	In Section 2 we define unitary Shimura varieties. However, since our focus is on the local structure at $p$ and the signature at infinite places of the Shimura variety, the group $G$ can be replaced by any other unitary group over $\mathbb{Q}$ that differs from $G$ at finite places away from $p$. The result in this paper applies to these compact unitary Shimura varieties.
	
	In Section 3, we introduce the stack of $G$-zips and the Ekedahl--Oort stratification of $X$. We also review the construction of (partial) flag spaces $Y_{P'}$ over $X$ and their stratifications. We then provide an explicit construction of the automorphic line bundles $\mathcal{L}_{P'}(\lambda)$. We consider all possible line bundles over all partial flag spaces (including the case $P=P'$, in which $Y_{P'}=X$) simultaneously. We state our main theorem and its variants at the end of the section.
	
	In Section 4, we prove the first key theorem in this paper, namely, the description of strata mentioned above. The goal of this section is to establish the geometric Jacquet--Langlands correspondence for certain strata of the flag spaces associated to different unitary Shimura varieties.
	
	In Section 5, we prove the necessity part of Theorem \ref{Main Theorem - Introduction}. The idea is straightforward: Since $\mathcal{L}_Y(\lambda)$ is ample, its restriction to any closed subscheme of $Y$ must remain ample. For case (1) in Theorem \ref{Main Theorem - Introduction} (1), we deduce the inequalities from the theory of flag varieties; For case (2) and (3) of Theorem \ref{Main Theorem - Introduction}, we select a special curve $C\subseteq Y$ and restricting to $C$ yields the desired result; For the remaining case, we choose a stratum $Z\subseteq Y$ whose dimension is strictly larger than 1 and which admits a geometric description. The key is to construct a (highly nontrivial) morphism $Z\to Z'$ whose fibers are flag varieties. Restricting to the fibers of this morphism yields the desired inequality.
	
	In Section 6, we introduce the notion of slopes. 
	This definition provides a combinatorical understanding of different strata in the flag space.
	
	In Section 7, we prove the sufficiency part of Theorem \ref{Main Theorem - Introduction}. The proof is lengthy and involved, but the idea is very basic: First, we reduce the proof of the ampleness criterion to the nefness criterion, then we further reduce the problem to showing that the line bundles $\mathcal{L}_{Y^i_j}(\{k^i_j\};\alpha)$ are nef under certain numerical conditions (See Section 7 for the explicit definitions). 
    To achieve this, we express the class $[\mathcal{L}_{Y^i_j}(\{k_i\};\alpha)]$ in the rational Picard group as a non-negative linear combination
	\begin{equation}
		[\mathcal{L}_{Y^i_j}(\{k^i_j\};\alpha)]=\sum_{l=1}^{l_0}u_l[h_l]+[\mathcal{L}_0],
	\end{equation}
	where $h_l$'s are strata Hasse invariants, and $\mathcal{L}_0$ is a line bundle (obtained as a pullback under a proper morphism) that is nef by induction. This shows that
	\begin{equation}
		\big([\mathcal{L}_{Y^i_j}(\{k^i_j\};\alpha)]\cdot C\big)\ge0
	\end{equation}
	as long as $C$ is not contained in the zero locus of any $h_l$. We can thus reduce the problem to a stratum of lower dimension. However, the actual construction is highly combinatorical and requires induction on some variants. We will give a more detailed overview at the beginning of Section 7.

	\section{Unitary Shimura varieties}
	
	\begin{rmk}
		The main result of this paper is purely local, and the structure of the reductive group $G$ at places away from $p$ and $\infty$ is not really important for us. We use the following moduli interpretation for simplicity. The result also works for other moduli problems, as long as the unitary group satisfies the condition that $G_{\mathbb{Q}_p}$ is isomorphic to the Weil restriction of $GL_n$ of a totally real field and the resulting Shimura variety is \emph{compact}. For example, the moduli problem in \cite{Kottwitz}.
	\end{rmk}
	
	\subsection{ Notations.} Let $p$ be a prime number. We fix an isomorphism $\iota: \mathbb{C}\stackrel{\sim}{\longrightarrow} \overline{\mathbb{Q}}_p$. Let $F$ be a totally real field of degree $N$ over $\mathbb{Q}$ in which $p$ is inert and $\mathfrak{p}$ be the unique prime above $p$. We denote by $\Sigma$ the set of places of $F$ and $\Sigma_\infty$ the subset of real embeddings. Through the isomorphism $\iota$, we view the set $\Sigma_\infty$ also as the set of $p$-adic embeddings $F\hookrightarrow \overline{\mathbb{Q}}_p$, and we can label $\Sigma_{\infty}$ by $\tau_{1},\dots,\tau_{N}$, such that 
	\begin{equation}
		\sigma^{-1}\tau_{i}=\tau_{i+1},
	\end{equation}
	where $\sigma$ is the Frobenius endomorphism.
	
	Let $E/F$ be a CM extension such that $\mathfrak{p}\mathcal{O}_E=\mathfrak{q}\mathfrak{q}^c$. We fix the choice of a prime $\mathfrak{q}$ over $\mathfrak{p}$. Let $\Sigma_{E,\infty}$ be the set of complex embeddings of $E$ and view them as $p$-adic embeddings as above. For each $\tau_i\in\Sigma_{\infty}$, let $\tilde{\tau}_i$ and $\tilde{\tau}^c_i$ be the embeddings of $E$ that extend $\tau_i$ and correspond to $\mathfrak{q}$ and $\mathfrak{q}^c$ respectively.
	
	Set $\kappa=\mathbb{F}_{p^N}$ and $k=\overline{\mathbb{F}}_p$. The inclusion $\kappa\hookrightarrow k$ induces an embedding of the ring of Witt vectors $W(\kappa)$ into $\overline{\mathbb{Q}}_p$.
	
	\medskip
	\subsection{ The Shimura datum.} Let $n\ge 2$ be an integer and $\{(m_i,n_i\}_{1\le i\le N}$ be a tuple of non-negative integers such that $m_i+n_i=n$ for all $1\le i\le N$. We will define the compact unitary Shimura variety of good reduction of signature $\{(m_i,n_i)\}_{1\le i\le N}$ in the following.
	
	\begin{itemize}
		\item Choose an $2n$-dimensional $E$-vector space $V$ equipped with a nondegenerate alternating pairing $\langle\cdot,\cdot\rangle:V\times V\longrightarrow \mathbb{Q}$ satisfying
		\begin{equation}
			\langle \alpha x,y\rangle=\langle x,\bar\alpha y\rangle,\quad \alpha\in E,
		\end{equation}
		and the following signature condition: For each $\tau\in \Sigma_{\infty}$, the pairing $\langle\cdot,\cdot\rangle$ induces an alternating pairing on the complex vector space $V_\tau=V\otimes_{F,\tau}\mathbb{R}$. This pairing is always the ``imaginary'' part of a uniquely determined Hermitian form $[\cdot,\cdot]_{\tau_i}$ on $V_{\tau_i}$. We require the signature of $[\cdot,\cdot]_{\tau_i}$ to be $(m_i,n_i)$.
		\item We define $G$ to be the unitary similitude group of $V$. Concretely, for a $\mathbb{Q}$-algebra $R$,
		\begin{equation}
			G(R)=\{(g,c(g))\in \mathrm{Aut}_{E\otimes_\mathbb{Q}R}(V\otimes_{\mathbb{Q}}R)\times R^\times,\  \langle gx,gy \rangle =c(g)\langle x,y\rangle,\  \forall x,y\in V\otimes_{\mathbb{Q}}R \}.
		\end{equation}
		It is easy to see that
		\begin{equation}
			G(\mathbb{R})=\mathrm{G}\bigg(\prod_{1\le i\le N} \mathrm{U}(m_i,n_i)\bigg).
		\end{equation}
		\item We consider the homomorphism of $\mathbb{R}$-algebraic groups $h:\mathbb{S}=\rm{Res}_{\mathbb{C}/\mathbb{R}}\mathbb{G}_m\to G_\mathbb{R}$ defined by
		\begin{equation}
			h(z)=\prod_{1\le i\le N} 
			\left(\begin{matrix}
				z & & & & &\\
				&\ddots & & & &\\
				& &z & & & \\
				& & &\bar{z} & & \\
				& & & &\ddots & \\
				& & & & &\bar{z} 
			\end{matrix}\right),
		\end{equation}
		where on the diagonal for each $i$, $z$ appears with multiplicity $m_i$ and $\bar{z}$ appears with multiplicity $n_i$. Here, the embedding $\textnormal{U}(m_i,n_i)\subseteq \textnormal{GL}_n(\mathbb{C})$ is induced by $\tau_i$. 
		Let $X$ be the $G(\mathbb{R})$-conjugacy class of $h$.
		
	\end{itemize}
	The above defines a Shimura datum $(G,X)$.
	
	\medskip
	\subsection{ Moduli interpretation.}
	
	Assume $V_{\mathbb{Q}_p}$ admits a self dual lattice. Let $K=K^pK_p\subseteq G(\mathbb{A}_\mathbb{Q}^f)$ be an open subgroup with hyperspecial level structure at $p$. In other words, $G$ extends to a reductive group scheme $\mathcal{G}$ over $\mathbb{Z}_p$ with $\mathcal{G}\otimes_{\mathbb{Z}_p}\mathbb{Q}_p\simeq G_{\mathbb{Q}_p}$ and $K_p=\mathcal{G}(\mathbb{Z}_p)$. For $K^p$ sufficiently small, the Shimura variety $Sh_K(G)$ can be realized as the moduli scheme representating the functor:
	\begin{equation}
		Sch^{\textnormal{loc noe}}_{/W(k)}\longrightarrow Sets
	\end{equation}
	sending a locally noetherian scheme $S$ over $W(k)$ to the set $Sh_K(G_\mathcal{P})(S)$ consisting of equivalent classes of tuples $(A,\lambda,\eta)$ satisfying:
	\begin{enumerate}
		\item $A$ is an abelian scheme of dimension $2nN$ over $S$ equipped with an action of $\mathcal{O}_E$, such that the characteristic polynomial of $\alpha\in\mathcal{O}_E$ on $\omega_{A^\vee/S}$ is given by
		\begin{equation}
			\prod_{1\le i\le N}(x-\tilde{\tau_i}(\alpha))^{m_i}(x-\tilde{\tau_i}^c(\alpha))^{n_i}.
		\end{equation}
		In other words, if $\omega_{A^\vee/S,\tilde\tau_i}$ (resp. $\omega_{A^\vee/S,\tilde\tau_i^c}$) denotes the subsheaf of $\omega_{A^\vee/S}$ on which $\mathcal{O}_E$ acts via $\tilde\tau_i$ (resp. via $\tilde\tau_i^c$), then it is a subbundle of rank $m_i$ (resp. $n_i$).
		\item $\lambda:A\to A^\vee$ is a prime-to-$p$ quasi-polarization of $A$ such that the Rosati involution associated with $\lambda$ induces the complex conjugation on $\mathcal{O}_E$.
		\item $\eta$ is a $K^p$-level structure on $A$. Explicitly, this means that $\eta$ is a collection of $\pi_1(S_j,\bar{s}_j)$-stable $K^p$-orbit of isomorphisms $\eta_j:V\otimes_\mathbb{Q}\mathbb{A}^{(p)}_f\stackrel{\sim}{\longrightarrow}V^{(p)}(A_{\bar{s}_j})$ for each connected component $S_j$ of $S$ with a geometric point $\bar{s}_j$. Here $V^{(p)}(A_{\bar{s}_j})= \varprojlim\limits_{p\nmid N}(A_{\bar{s}_j}[N])\otimes_\mathbb{Z}\mathbb{Q}$ is the rational Tate module away from $p$. We further require that for some isomorphism $\nu(\eta_j)\in \text{Hom}({\mathbb{A}}_f^{(p)},{\mathbb{A}}_f^{(p)}(1))$, the following diagram is commutative:
		
		\centering{
			\begin{tikzcd}
				V\otimes_{\mathbb{Q}} \mathbb{A}_f^{(p)}\times V\otimes_{\mathbb{Q}} \mathbb{A}_f^{(p)} \arrow[rr, "{\langle\cdot,\cdot\rangle}"] \arrow[d, "\eta_j", shift left=10] \arrow[d, "\eta_j"', shift right=10] &  & {\mathbb{A}}_f^{(p)} \arrow[d, "\nu(\eta_j)"] \\
				V^{(p)}(A_{\bar{s}_j})\times V^{(p)}(A_{\bar{s}_j}) \arrow[rr, "\textnormal{Weil pairing}"]                                                                                    &  & {\mathbb{A}}_f^{(p)}(1).              
		\end{tikzcd}}
	\end{enumerate}
	Two triples $(A,\lambda,\eta)$ and $(A',\lambda',\eta')$ are equivalent if there is an $\mathcal{O}_E$-equivariant prime-to-$p$ quasi-isogeny $\phi:A\to A'$ such that $\phi^\vee\circ\lambda'\circ\phi=\lambda$ and $\eta'=\phi\circ \eta$.
	
	For the rest of the paper, we denote by $X=Sh_K(G)_k$ the geometric special fiber of the Shimura variety. We omit the level structure for simplicity.

	\section{$G$-zips and stratification of the special fiber}
	
	\indent In this section, we recall the notion of $G$-zips introduced by Pink, Wedhorn and Ziegler in \cite{Pink-Wedhorn-ZieglerAlgebraicZipData}. The theory of $G$-zips provides a group-theoretical description of the Ekedahl--Oort stratification on the special fiber of Shimura varieties. We also recall some relevant results by Goldring, Koskivirta, Imai, Wedhorn, etc.
	
	\medskip
	\subsection{ Notations.}
	
	Let $p$ be a prime number. Let $k=\overline{\mathbb{F}}_p$ be the algebraic closure of $\mathbb{F}_p$ as above. Let $G$ be a connected reductive group over $\mathbb{F}_p$. For a scheme $X$ over $\mathbb{F}_p$, we denote by $X^{(p)}$ its $p$-th power Frobenius twist and by $\varphi:X\to X^{(p)}$ the relative Frobenius morphism. Let $(B,T)$ be a Borel pair of $G$, meaning that $B\subseteq G$ is a Borel subgroup and $T\subseteq B$ is a maximal torus of $G$. We will use the following notations:
	
	$\bullet$ $X^\ast(T)$ (resp. $X_\ast(T)$) denotes the group of characters of $T$;
	
	$\bullet$ $W=W(G_k,T)=N_G(T)/T$ denotes the Weyl group of $G$;
	
	$\bullet$ $(\Phi,\Phi^+,\Delta)$ is the root system attached to the Borel pair $(B,T)$, i.e., $\Phi\subseteq X^\ast(T)$ is the set of $T$-roots in $G$, $\Phi^+\subseteq \Phi$ is the set of positive roots with respect to $B$, and $\Delta\subseteq \Phi^+$ is the set of simple roots;
	
	$\bullet$ For $\alpha\in\Phi$, let $s_\alpha\in W$ be the corresponding simple reflection. Let $\ell:W\to\mathbb{N}$ be the length function on $W$. We endow $W$ with the Bruhat order: $u\le v$ if and only if $v$ can be written in the form $v=s_{\alpha_1}\cdots s_{\alpha_{\ell(v)}}$ such that $u$ has a presentation of the form $s_{\alpha_{i_1}}\cdots s_{\alpha_{i_{\ell(u)}}}$ with $i_1<i_2<\cdots<i_{\ell(u)}$.
	
	$\bullet$ For a subset $I\subseteq \Delta$, let $W_I$ be the subgroup of $W$ generated by $\{s_\alpha|\ \alpha\in I\}$.
	
	$\bullet$ Let $w_0$ and $w_{0,I}$ denote the longest element in $W$ and $W_I$ respectively.
	
	$\bullet$ Let $^IW$ (resp. $W^I$) be the set of elements $w\in W$ of minimal length in the coset $W_Iw$ (resp. $wW_I$). Then $^IW$ (resp. $W^I$) is a set of representatives of the coset space $W_I\backslash W$ (resp. $W/W_I$). The longest element in $W^I$ is $w_0w_{0,I}$.
	
	$\bullet$ Let $X^\ast_+(T)$ be the set of dominant characters, i.e., characters	$\lambda\in X^\ast(T)$ that satisfies $\langle \lambda,\alpha^\vee\rangle\ge0$ for all $\alpha\in \Delta$.
	
	$\bullet$ Let $X^\ast_{+,I}(T)$ be the set of \emph{$I$-dominant characters}, i.e., characters $\lambda\in X^\ast(T)$ satisfying $\langle \lambda,\alpha^\vee\rangle\ge0$ for all $\alpha\in I$.
	
	$\bullet$ Let $P\subseteq G_k$ be a parabolic subgroup containing $T$. Then $P$ is conjugate in $G_k$ to a unique standard parabolic subgroup $P_I$ containing $B$. We say $P$ is of type $I_P:=I$. Let $\Delta^P=\Delta\backslash I$.

	\medskip
	\subsection{ The stack of $G$-zips.}
	
	\begin{definition}[\cite{Pink-Wedhorn-ZieglerF-zips,Pink-Wedhorn-ZieglerAlgebraicZipData}]
		Let $G$ be a connected reductive group over $\mathbb{F}_p$. A \emph{connected zip datum} is a tuple $\mathcal{Z}:=(G,P,L,Q,M,\varphi)$, where
		
		$\bullet$ $P,Q\subseteq G_k$ are parabolic subgroups.
		
		$\bullet$ $L\subseteq P$ and $M\subseteq Q$ are Levi subgroups and
		
		$\bullet$ $\varphi:L\to M$ is an isogeny.
	\end{definition}
	
	Given a choice of Levi subgroup $L$ (resp. $M$), every element $x\in P$ (resp. $y\in Q$) admits a unique decomposition $x=ul$ (resp. $y=vm$) where $u\in R_u(P)$ and $l\in L$ (resp. $v\in R_u(Q)$ and $m\in M$). We define the zip group as
	\begin{equation}
		E_{\mathcal Z}:=\{(x,y)\in P\times Q |\  \varphi(\bar{x})=\bar{y}\}.
	\end{equation}
	Here, $\overline{\bullet}$ means the Levi factor of $\bullet$. The zip group $E_\mathcal{Z}$ acts on $G$ via
	\begin{equation}
		(x,y)\cdot g:=x^{-1}gy.
	\end{equation}
	
	The following example provides a typical construction of zip data.
	
	\begin{example}
		By a cocharacter datum, we mean a pair $(G,\mu)$, where $G$ is a connected reductive group over $\mathbb{F}_p$ and $\mu:G_{m,k}\to G_k$ is a cocharacter. A cocharacter datum gives rise to a connected zip datum $\mathcal{Z}_\mu=(G,P,Q,L,M,\varphi)$ as follows:
		
		$\bullet$ $Q:=(P^+)^{(p)}$, where $P^+$ is the parabolic subgroup containing $B$ such that
		\begin{equation}
			P^+(k)=\{g\in G(k)\ |\ \lim_{t\to 0}\ \mu(t)g\mu(t)^{-1}\textnormal{ exists}\}.
		\end{equation}
		
		$\bullet$ $P:=P^-$ is the opposite parabolic subgroup of $P^+$.
		
		$\bullet$ $L:=P^-\cap P^+=\textnormal{Cent($\mu$)}$ and $M:=L^{(p)}$.
		
		$\bullet$ $\varphi:L\to M$ is the relative Frobenius morphism.
		
		We call this tuple a \emph{connected zip datum of type $(G,\mu)$}.
	\end{example}	
	
	\begin{definition}
		Let $S$ be a $k$-scheme and $\mathcal{Z}$ a connected zip datum. A \emph{$G$-zip of type $\mathcal{Z}$} over $S$ is a tuple $\underline{I}=(I,I_P,I_Q,\psi)$ where $I$ is a $G$-torsor over $S$, $I_P\subseteq I$ is a $P$-torsor, $I_Q\subseteq I$ is a $Q$-torsor, and $\psi:(I_P/R_uP)^{(p)}\stackrel{\sim}{\longrightarrow} I_Q/R_uQ$ is an isomorphism of $M$-torsors.
	\end{definition}
	
	The stack of $G$-zips, denoted by $G\textnormal{-zip}^{\mathcal{Z}}$, 
	is defined to be the classifying stack for $G$-zips of type $\mathcal{Z}$ over $k$. If $\mathcal{Z}$ arises from a cocharacter datum $(G,\mu)$, we simply denote this stack by $G\textnormal{-zip}^\mu$. Pink, Wedhorn and Ziegler proved the following in \cite{Pink-Wedhorn-ZieglerAlgebraicZipData}:
	\begin{thm}
		$G\textnormal{-zip}^\mathcal{Z}$ is a smooth algebraic stack of dimension 0 over $k$, and
		\begin{equation}
			G\textnormal{-zip}^\mathcal{Z}\simeq [E_\mathcal{Z}\backslash G].
		\end{equation}
	\end{thm}
	
	\begin{definition}
		Let $\mathcal{Z}$ be a connected zip datum. We call a triple $(B,T,z)$ a \emph{frame} for $\mathcal{Z}$ if $(B,T)$ is a Borel pair of $G_k$ defined over $\mathbb{F}_p$ and $z\in W$, such that
		\begin{equation}
			B\subseteq P,\qquad ^zB\subseteq Q,\qquad \varphi(B\cap L)=B\cap M= {^z}B\cap M.
		\end{equation}
	\end{definition}
	
	Note that this differs slightly from the original definition in \cite{Pink-Wedhorn-ZieglerAlgebraicZipData}, where they do not require the Borel pair to be defined over $\mathbb{F}_p$. Under our definition, a frame does not necessarily exist in general. However, in the case where $\mathcal{Z}=\mathcal{Z}_\mu$ comes from a cocharacter datum $(G,\mu)$, there exists $g\in G(k)$ such that $\mathcal{Z}_{\mu'}$ where $\mu'=\textnormal{ad}_g\circ\mu$ admits a frame\cite{IK-partial}. For our purpose (of studying stratification and automorphic vector bundles), it is harmless to assume the existence of frames.

	The orbit decomposition of $G$ with respect to the action of the zip group $E_{\mathcal{Z}}$ is actually a stratification.
	
	\begin{thm}[\cite{Pink-Wedhorn-ZieglerAlgebraicZipData}] Let $(B,T,z)$ be a frame for $\mathcal{Z}$ and $(\Phi,\Phi^+,\Delta)$ be the corresponding root system. Let $I$ and $J$ be the types of the parabolic subgroups $P$ and $Q$, respectively. Then
		
		1. (Orbit decomposition) There are bijections
		\begin{equation}
			\begin{aligned}
				& ^IW\longrightarrow \{E\textnormal{-orbits in $G$}\},\\
				& W^J\longrightarrow \{E\textnormal{-orbits in $G$}\}.
			\end{aligned}
		\end{equation}
		The first bijection is given by $w\mapsto o_{\mathcal{Z}}(g\dot{w})$. Here, $\dot{w}$ is a lift of $w$ to $N_G(T)\subseteq G$ and $o_{\mathcal{Z}}(\cdot)$ represents the orbit of an element. We denote the orbit corresponding to $w\in$$^IW$ or $W^J$ by $G_w$.
		
		2. (Dimension) $\textnormal{dim }G_w=\ell(w)+\textnormal{dim }P$.
		
		3. (Closure relation) For $w,w'\in$$^IW$ we write $w'\preccurlyeq w$ if and only if there exists $y\in W_I$ such that $y^{-1}w'\psi(y)\le w$. Here, $\psi:W_I\stackrel{\sim}{\rightarrow}{W_J}$ is the morphism induced by $\varphi\circ \textnormal{int}(z):P\to {^z}Q$ on the standard parabolics. We refer to this order as the twisted Bruhat order. Then, for any $w\in ^IW$,
		\begin{equation}
			\overline{G_w}=\coprod_{\substack{v\in ^IW\\v\preccurlyeq w}} G_v.
		\end{equation}
		
		4. (Geometry) Each orbit $G_w$ is a locally closed smooth subscheme of $G$.
	\end{thm}
	
	\begin{rmk}
		The closure $\overline{G_w}$ may be non-smooth. It is even non-normal in general. This lack of smoothness is one of the main difficulties in constructing strata Hasse invariants. In \cite{KoskivirtaNormalization}, the normalization of the closures of strata is studied, and we will use these results in our later discussion.
	\end{rmk}

	\medskip
	\subsection{ The Ekedahl--Oort stratification.}	The stack of $G$-zips provides a group-theoretic parametrization of the Ekedahl--Oort stratification of the special fiber of (Hodge type) Shimura varieties, which we now explain. We hope to present a general construction, so we temporarily forget the notations in Section 2 in this subsection.
	
	Let $(G,X)$ be a Shimura datum. We assume that the group $G$ is unramified at $p$. Then $G$ extends to a reductive group $\mathcal{G}$ over $\mathbb{Z}_{(p)}$ such that $\mathcal{G}_\mathbb{Q}=G$. Let $K_p:=\mathcal{G}(\mathbb{Z}_p)\subseteq G(\mathbb{Q}_p)$ be a hyperspecial subgroup. 
	
	Let $E=E(G,X)$ be the reflex field of $(G,X)$ and let $\mathcal{O}_E$ denote its ring of integers. Let $\mathfrak{p}$ be a prime of $\mathcal{O}_E$ above $p$. We denote by $E_\mathfrak{p}$ and $\mathcal{O}_{E,\mathfrak{p}}$ the completion of $E$ and $\mathcal{O}_E$ at $\mathfrak{p}$ respectively. If $(G,X)$ is of abelian type, by the work of Kisin \cite{KisinIntegralModel} and Vasiu \cite{VasiuIntegralModel}, for a sufficiently small open subgroup $K^p\subseteq G(\mathbb{A}^{p}_f)$ and $K=K^pK_p$, the Shimura variety $Sh(G,X)_K/E$ admits a canonical integral model $\mathcal{S}_K$ defined over $\mathcal{O}_{E,\mathfrak{p}}$. For another compact open subgroup $K'^p\subseteq G(\mathbb{A}^p_f)\subseteq K^p$, let $K'=K'^pK_p\subseteq K$. The inclusion of levels induces a finite etale projection
	\begin{equation}
		\pi_{K'/K}:\mathcal{S}_{K'}\longrightarrow \mathcal{S}_{K}.
	\end{equation}
	
	Let $X=(\mathcal{S}_K)_k$ be the geometric special fiber of the integral model, and we omit the level structure for simplicity. Since 
	finite pullbacks of ample line bundles remain ample, this simplification is harmless for the purpose of studying ampleness in this paper.
	
	When the Shimura datum $(G,X)$ is of Hodge type, we have the universal abelian scheme $\pi:\mathcal{A}\to X$. The Hodge-to-de Rham spectral sequence degenerates at the $E_1$-page, as shown in \cite{OdaDeRham}. We have the Hodge filtration
	\begin{equation}
		\textnormal{Fil}_H:\ 0\longrightarrow \omega_{\mathcal{A}/X}\longrightarrow \H(\mathcal{A}/X)\longrightarrow \textnormal{Lie}_{\mathcal{A}/X}\longrightarrow0
	\end{equation}
	and the conjugate filtration
	\begin{equation}
		\textnormal{Fil}_{\textnormal{conj}}:\ 0\longrightarrow \textnormal{Ker }V \longrightarrow \H(\mathcal{A}/X)\longrightarrow \omega_{\mathcal{A}/X}^{(p)}\longrightarrow 0,
	\end{equation}
	together with comparison isomorphisms between their subquotients
	\begin{equation}
		\begin{aligned}
			&\psi_1:\qquad \H(\mathcal{A}/X)/\textnormal{Ker $V$}\xlongrightarrow[\simeq]{V} \omega_{\mathcal{A}^\vee/X}^{(p)},\\
			&\psi_2:\qquad \textnormal{Ker $V$}\xlongrightarrow[\simeq]{F^{-1}}\big(\H(\mathcal{A}/X)/\omega_{\mathcal{A}^\vee/X}\big)^{(p)}.
		\end{aligned}
	\end{equation}
	Here, $\H(\mathcal{A}/X):=H^1_\textnormal{dR}(\mathcal{A}^\vee/X)$ denotes the de Rham homology of $\mathcal{A}$, and $\omega_{\mathcal{A}/X}=\pi_\ast \Omega_{\mathcal{A}/X}$. Let $\mathcal{Z}_\mu=(G,P,Q,L,M,\varphi)$ be the zip datum induced by $(G,\mu)$, where $\mu\in X$. The datum $(\H(\mathcal{A}^\vee/X),\textnormal{Fil}_H,\textnormal{Fil}_{\textnormal{conj}},\psi)$ induces a morphism
	\begin{equation}
		\zeta: X\longrightarrow G\textnormal{-zip}^\mu.
	\end{equation}
	For a more detailed explanation, see \cite{GK-strata}.
	
	\begin{thm}[\cite{ZhangEOStrata}] 
		The morphism $\zeta$ is smooth.
	\end{thm}
	
	\begin{definition}
		We define the Ekedahl--Oort stratification (or EO--stratification for short) of $X$ as
		\begin{equation}
			X=\coprod_{w\in {^I}W}\zeta^{-1}(G_w) \qquad \textnormal{or}\qquad X=\coprod_{w\in W^J}\zeta^{-1}(G_w)
		\end{equation}
	\end{definition}	
	
	We denote $\zeta^{-1}(G_w)$ by $X_w$. By the smoothness of $\zeta$, the EO-stratification has the following properties.
	
	\begin{prop} 
		1. (Dimension) $\textnormal{dim } X_w=\ell(w)$.
		
		2. (Closure relation)
		\begin{equation}
			\overline{X_w}=\coprod_{\substack{v\in ^IW\\v\preccurlyeq w}} X_v.
		\end{equation}
		
		3. (Geometry) Each EO--stratum $X_w$ is a locally closed subvariety of $X$.
	\end{prop}
	
	\noindent\textbf{Convention: }We will refer to $X_w$ and $\overline{X_w}$ as open and closed strata of $X$ respectively. By abuse of language, we will refer to both as strata.
	
	\begin{example}\label{Computation of Ekedahl--Oort stratum}
		We follow the conventions and definitions in Section 2, and explain how to determine the conditions on each EO--stratum in our unitary setup. Fixing the prime $\mathfrak{q}$ over $\mathfrak{p}$, the group $G$ over $\mathbb{F}_p$ is (non-canonically) isomorphic to $\textnormal{Res}_{\mathbb{F}_{p^N}/\mathbb{F}_p}GL_{n,\mathbb{F}_{p^N}}\times \mathbb{G}_m$. Thus, 
		the Weyl group $W$ is isomorphic to 
		$
		\prod_{\tau\in \Sigma_{\infty}}\mathcal{S}_n.
		$
		Here, each $\tau_i\in\Sigma_{\infty}$ induces an embedding of $\mathbb{F}_{p^N}$ into $k$, and $\mathcal{S}_n$ is the group of permutations of $n$-elements. The inverse Frobenius $\sigma^{-1}$ acts on $W$ via
		\begin{equation}
			\sigma^{-1}(w_1,\dots, w_{N})=(w_2,\dots, w_{N},w_1).
		\end{equation}
		The Levi subgroup $L_\mu$ of $P=P_{\mu}\subseteq G_k$ is isomorphic to
		\begin{equation}
			\left(\prod_{i=1}^N GL_{m_i}\times GL_{n_i} \right)\times \mathbb{G}_m.
		\end{equation}
		Thus, the Weyl group $W_I$ is isomorphic to the subgroup 
		\begin{equation}
			\prod_{i=1}^N \mathcal{S}_{m_i}\times \mathcal{S}_{n_i}
		\end{equation}
		of $W$. The coset
		\begin{equation}
			^IW=\Big\{(w_1,\dots,w_N)\in W | w_i^{-1}(1)<\cdots<w_i^{-1}(m_i),\ w_i^{-1}(m_i+1)<\cdots w_{i}^{-1}(n),\forall 1\le i\le N\Big\}.
		\end{equation}
		For $1\le i\le N$, Let $V^i$ be a standard $k$-vector space of dimension $n$ with basis $\{e^i_j\}_{1\le j\le n}$, and let $V=\oplus_{i=1}^N V^i$. For any $w=(w_1,\dots,w_N)\in$$^IW$, we define
		\begin{equation}
			\begin{aligned}
				&V^i_H:= ke^i_{w_i^{-1}(1)}\oplus \cdots \oplus ke^i_{w_i^{-1}(m_i)}, \qquad \textnormal{Fil}^{i,\textnormal{std}}_H:0\subseteq V^i_H\subseteq V^i\\
				&V^i_{\textnormal{conj}}:=ke^i_1\oplus \cdots\oplus ke^i_{n_{i+1}},\qquad\qquad\ \  \textnormal{Fil}^{i,\textnormal{std}}_\textnormal{conj}:0\subseteq V^i_{\textnormal{conj}}\subseteq V^i.\\
				&\psi^{i,\textnormal{std}}_1:V^i/V^i_{\textnormal{conj}}\stackrel{\simeq}{\longrightarrow} V^{i+1,(p)}_H \quad\qquad\textnormal{given by } \psi^{i,\textnormal{std}}_1(e^i_{n_{i+1}+j})=e^{i+1,(p)}_{w_{i+1}^{-1}(j)}, \quad\textnormal{for $1\le j\le m_{i+1}$}.\\
				&\psi^{i,\textnormal{std}}_2:V^i_{\textnormal{conj}}\stackrel{\simeq}{\longrightarrow} (V^{i+1}/V^{i+1,(p)}_H)\quad\  \textnormal{given by } \psi^{i,\textnormal{std}}_2(e^i_j)=e^{i+1,(p)}_{w_{i+1}^{-1}(m_{i+1}+j)},\quad \textnormal{for $1\le j\le n_{i+1}$}.
			\end{aligned}
		\end{equation}
	Here we view the superscripts as modulo $N$. We refer the tuple
	\begin{equation}
		\Big(V=\bigoplus_{i=1}^N V^i, \{\textnormal{Fil}^{i,\textnormal{std}}_H\}, \{\textnormal{Fil}^{i,\textnormal{std}}_{\textnormal{conj}}\}, \{\psi^{i,\textnormal{std}}_1,\psi^{i,\textnormal{std}}_2\}\Big)
	\end{equation}
	constructed above as a standard $G$-zip of type $\mathcal{Z}_\mu$ and shape $w$.
	
	By the moduli interpretation in Section 2, the Hodge and conjugate filtrations decompose into 
	\begin{equation}
		\begin{aligned}
			&\textnormal{Fil}^i_H: 0\subseteq \omega_{\mathcal{A}^\vee/X,\tilde\tau_i}\subseteq \H(\mathcal{A}/X)_{\tilde\tau_i},\\
			&\textnormal{Fil}^i_\textnormal{conj}:0\subseteq \textnormal{Ker $V$}_{\tilde\tau_i} \subseteq\H(\mathcal{A}/X)_{\tilde\tau_i},
		\end{aligned}
	\end{equation}
	and the comparison isomorphisms are given by
	\begin{equation}
		\begin{aligned}
			&\psi^i_1: \H(\mathcal{A}/X)_{\tilde\tau_i}/\textnormal{Ker $V$}_{\tilde\tau_i}\stackrel{\simeq}{\longrightarrow} \omega^{(p)}_{\mathcal{A}^\vee/X,\tilde\tau_{i+1}},\\
			&\psi^i_2: \textnormal{Ker $V$}_{\tilde\tau_i}\stackrel{\simeq}{\longrightarrow} (\H(\mathcal{A}/X)_{\tilde\tau_{i+1}}/\omega_{\mathcal{A}^\vee/X,\tilde\tau_{i+1}})^{(p)},
		\end{aligned}
	\end{equation}
	and similarly for the $\tilde{\tau}_i^c$-components. We can describe the conditions for the Ekedahl--Oort strata of $X$ by the proposition below.
	\begin{prop}
		A closed point $s\in X(k)$ lies in $X_w(k)$ if and only if there is a basis of the $k$-vector space $\H(\mathcal{A}_k/k)$ under which $(\H(\mathcal{A}_k/k),\allowbreak\{\textnormal{Fil}^i_H\},\allowbreak\{\textnormal{Fil}^i_{\textnormal{conj}}\},\{\psi^i_1,\psi^i_2\})$ is isomorphic to the above standard $G$-zip of type $\mathcal{Z}_\mu$ and shape $w$.
	\end{prop}
	\end{example}

	\begin{rmk}
		\begin{enumerate}
			\item A similar construction for $GL_n$-zips can be found in the work of Moonen and Wedhorn \cite{Moonen-Wedhorn}. In our paper, we take the conjugate filtration to be the standard filtration. This has benefits when we compute the tangent bundle of certain $\overline{X_w}$.
			\item The conditions defining the strata are typically difficult to express explicitly because we must account for the comparison isomorphism $\psi$. This also explains why the closed strata are usually non-normal, unlike Schubert cells, which are given by conditions of the type $\textnormal{dim }\mathcal{M}\cap \mathcal{N}\ge r$.
		\end{enumerate}
	\end{rmk}

	\medskip
	\subsection{ $G$-zipflags and automorphic bundles}
	
	The goal of this subsection is to introduce the flag space, or more generally, partial flag spaces, as well as their stratifications. Our proof of the ampleness criterion is based on a careful analysis of the geometry of these auxiliary varieties.
	
	Let $\mathcal{Z}=(G,P,L,Q,M,\varphi)$ be a connected zip datum and $(B,T,z)$ be a frame for $\mathcal{Z}$. 
	
	\begin{definition}
		Let $B\subseteq P_0\subseteq P$ be an intermediate parabolic subgroup. A \emph{$G$-zipflag of type $(\mathcal{Z},P_0)$} over a $k$-scheme $S$ is a pair $\underline{I}^{P_0}=(\underline{I},I_{P_0})$, where $\underline{I}=(I,I_P,I_Q,\psi)$ is a $G$-zip of type $\mathcal{Z}$ over $S$, and $I_{P_0}\subseteq I_P$ is a $P_0$-torsor.
	\end{definition}
	
	We define \emph{the stack of $G$-zipflag of type $(\mathcal{Z},P_0)$}, denoted by $G\textnormal{-zipflag}^{\mathcal{Z},P_0}$, to be the stack over $k$ such that for any $k$-scheme $S$, $G\textnormal{-zipflag}^{\mathcal{Z},P_0}(S)$ is the groupoid of $G$-zipflags of type $(\mathcal{Z},P_0)$ over $S$. If $P_0=P$, this reduces to the original stack of $G$-zips; if $P_0=B$, we typically omit the superscript $B$.
	
	Goldring and Koskivirta provide a characterization of this stack in \cite{GK-stratification,GK-strata}. 
	\begin{thm} Let
		\begin{equation}
			\hat{E}_{P_0}:=\{(p,q)\in P\times Q,\ p\in P_0\}\subseteq E_{\mathcal{Z}}.
		\end{equation}
		Let $\hat{E}_{P_0}$ act on $G$ via
		\begin{equation}
			(p,q)\cdot g:=p^{-1}gq.
		\end{equation}
		
		1. The stack $G\textnormal{-zipflag}^{\mathcal{Z},P_0}$ is smooth of dimension equal to $\textnormal{dim }P/P_0$ and is isomorphic to $[\hat{E}_{P_0}\backslash G]$.
		
		2. For all parabolic subgroups $B\subseteq P_1\subseteq P_0\subseteq P$, we have a commutative diagram
		\begin{equation}
			\begin{tikzcd}[column sep=large]
				G\textnormal{-zipflag}^{\mathcal{Z},P_1} \arrow[r, "{\pi_{P_1,P_0}}"] \arrow[d, "\simeq"'] & G\textnormal{-zipflag}^{\mathcal{Z},P_0} \arrow[d, "\simeq"] \\
				{[\hat{E}_{P_1}\backslash G]} \arrow[r]                                      & {[\hat{E}_{P_0}\backslash G]}. 
			\end{tikzcd}
		\end{equation}
		
		3. The projection $\pi_{P_1,P_0}$ is proper smooth, with fibers isomorphic to $P_0/P_1$.
	\end{thm}
	
	\begin{definition}
		We define \emph{the partial flag space of type $P_0$}, denoted by $Y_{P_0}$, to be the fiber product
		\begin{equation}
			\begin{tikzcd}[column sep=large]
				Y_{P_0} \arrow[r, "\zeta"] \arrow[d, "\pi_{P_0}"'] & G\textnormal{-zipflag}^{\mathcal{Z},P_0} \arrow[d, "\pi_{P_0}"] \\
				X \arrow[r, "\zeta"]   & G\textnormal{-zip}^\mathcal{Z}. 
			\end{tikzcd}
		\end{equation}
	\end{definition}
	Here, we denote both vertical projections $\pi_{P_0}$ and both horizontal arrows by $\zeta$ by abuse of notation. When $P_0=B$, we write $Y_B=Y$ and simply refer to it as the \emph{flag space}.
	
	\medskip
	\noindent\textbf{The coarse stratification. }In \cite{GK-stratification}, the authors introduced two stratifications of the stack $G\textnormal{-zipflag}^{\mathcal{Z},P_0}$. The strata that will appear in our proof of the ampleness criterion are defined by conditions on the intersections of vector bundles, i.e., conditions of the form $\textnormal{rank }(\mathcal{M}\cap\mathcal{N})\ge t$. This is referred to as the \emph{coarse stratification} by Goldring and Koskivirta.
	
	Put $Q_0:=(L\cap P_0)^{(p)}R_u(Q)$. This is a parabolic subgroup of $G$ such that
	\begin{equation}
		^zB\subseteq Q_0\subseteq Q,\qquad \hat{E}_{P_0}\subseteq P_0\times Q_0.
	\end{equation}
	Let $I_0\subseteq I$ and $J_0\subseteq J$ denote the types of $P_0$ and $Q_0$, respectively. The inclusion of groups induces a natural smooth surjective morphism of stacks
	\begin{equation}
		\iota_{P_0}:G\textnormal{-zipflag}^{\mathcal{Z},P_0}=[\hat{E}_{P_0}\backslash G]\longrightarrow [P_0\times Q_0\backslash G].
	\end{equation}
	The right hand side is bijective to $^{I_0}W^{J_0}$. We define the \emph{coarse flag strata} of $G\textnormal{-zipflag}^{\mathcal{Z},P_0}$ as the fibers of this morphism endowed with reduced structure. By further taking the inverse image under $\zeta$, we obtain the \emph{coarse flag strata} of the partial flag space $Y_{P_0}$. The closure relation of the coarse strata is given by the Bruhat order. Again, by abuse of language, we refer to the closure of a coarse stratum as a coarse stratum.
	

	\medskip
	\noindent\textbf{Vector bundles on quotient stacks.} Let $G$ be a $k$-algebraic group acting on a $k$-scheme $X$. To any finite dimensional representation $(\rho,V)$ of $G$, we can associate a vector bundle $\mathcal{V}(\rho)$ on the stack $[G\backslash X]$ such that
	\begin{equation}
		H^{0}([G\backslash X],\mathcal{V}(\rho))=\big\{f:X\longrightarrow V\big|\ f(g\cdot x)=\rho(g)f(x),\quad \forall g\in G,x\in X \big\}.
	\end{equation}

	\medskip
	\noindent\textbf{Automorphic bundles.} Let $H$ be any split connected reductive algebraic group over an arbitrary field $K$ with a Borel pair $(B_H,T_H)$. For any irreducible representation $(\rho,V)$ of $H$, there is a unique highest weight vector $v\in V$ with weight $\lambda$. There is a bijection between the set of irreducible representations of $H$ over $K$ and the set of dominant characters $X^\ast_{+}(T_H)$.
	
	In our case, let $\lambda\in X^\ast_{+,I}$ be an $I$-regular dominant weight. For any intermediate parabolic subgroup $B\subseteq P_0\subseteq P$, let $L_0$ be the Levi subgroup of $P_0$ containing $T$ and $I_0$ be the type of $P_0$. Then $\lambda$ is $I_0$-regular dominant. Let $B_{L_0}=B\cap L_0$, so that $(B_{L_0},T)$ is a Borel pair of $P_0$. The space
	\begin{equation}
		V(\lambda):=H^0(L_0/B_{L_0},\mathcal{L}_\lambda)
	\end{equation}
	is a representation of $L_0$ with highest weight $\lambda$. We view this as a representation of $\hat{E}_{P_0}$ via
	\begin{equation}
		\hat{E}_{P_0}\xlongrightarrow{\textnormal{pr}_1} P_0\xlongrightarrow{\pi} L_0.
	\end{equation}
	Here, the first map is the projection onto the first factor and the second is the natural quotient map. By the previous discussion, we obtain a vector bundle, denoted by $\mathcal{V}(\lambda)$, on the stack $G\textnormal{-zipflag}^{\mathcal{Z},P_0}=[\hat{E}_{P_0}\backslash G]$.
	\begin{definition}
		We call $\mathcal{V}(\lambda)$ as well as its pullback $\mathcal{V}_{Y_{P_0}}(\lambda)$ to $Y_{P_0}$ by $\zeta$ \emph{the automorphic vector bundle of weight $\lambda$}. If $P_0=B$, then $V(\lambda)$ is a 1-dimensional representation, we write $\mathcal{V}_Y(\lambda)$ as $\mathcal{L}_Y(\lambda)$ instead and refer to it as an \emph{automorphic line bundle of weight $\lambda$}.
	\end{definition}
	
	\begin{definition}
		Let $\chi:L_0\to k^\times$ be a character of $L_0$, again viewed as a character of $\hat{E}_{P_0}$. This defines a line bundle $\mathcal{L}(\chi)$ on $G\textnormal{-zipflag}^{\mathcal{Z},P_0}$ and its pullback $\mathcal{L}_{Y_{P_0}}(\chi)$ on $Y_{P_0}$. Note that if $P_0=B$, this construction agrees with the one described above.
	\end{definition} 
	
	\begin{rmk}
		For $P_0=P$ or $B$, the above definition agrees with the classical notion of automorphic vector bundles and automorphic line bundles. Our definition extends this notion to all intermediate partial flag spaces. We believe this extension of definition is natural for the following reasons.
		
		1. They also arise from irreducible representations of parabolic subgroups of $G$, and the pushforward of $V_{Y_{P_0}}(\lambda)$ to another partial flag space, whenever this is reasonable, corresponds to the parabolic induction. Therefore, the study of such bundles is closely related to representation theory.
		
		2. By the Borel--Weil--Bott theorem, on has
		\begin{equation}
			\pi_\ast(\mathcal{L}_{Y}(\lambda))=\mathcal{V}_X(\lambda).
		\end{equation}
		By the Kempf vanishing theorem (see \cite{JantzenRepresentationBook}), we moreover have
		\begin{equation}
			H^i(Y,\mathcal{L}_Y(\lambda))=H^i(X,\mathcal{V}_X(\lambda)).
		\end{equation}
		Indeed, for any intermediate parabolic subgroups $B\subseteq P_0\subseteq P_1\subseteq P$, if we denote by $\pi^{P_0}_{P_1}$ the projection $Y_{P_0}\to Y_{P_1}$, then we have
		\begin{equation}
			\pi_{P_1,\ast}^{P_0}(\mathcal{V}_{Y_{P_0}}(\lambda))=\mathcal{V}_{Y_{P_1}}(\lambda),\quad \textnormal{and}\quad H^i(Y_{P_0},\mathcal{V}_{Y_{P_0}}(\lambda))=H^i(Y_{P_1},\mathcal{V}_{Y_{P_1}}(\lambda)).
		\end{equation}
		Moreover, these isomorphisms are all equivariant for tame Hecke operators. Hence, the coherent cohomology of $\mathcal{V}_X(\lambda)$ on $X$, which is of significant interest to experts, is closely related to the coherent cohomology of the automorphic bundles on these partial flag spaces.
		
		3. The proof of the ampleness criterion in this paper heavily relies on these intermediate partial flag spaces and line bundles defined on them. Therefore, it makes sense to unify the notions.
	\end{rmk}
	
	It is natural to pose the following fundamental question.
	
	\medskip
	\noindent\textbf{Question: }For which weights $\lambda$, is the line bundle $\mathcal{L}_Y(\lambda)$ ample on $Y$?

    \medskip
	This paper is devoted to providing a complete answer to this question for unitary Shimura varieties. We believe that the machinery of the geometric description of certain strata (or geometric Jacquet--Langlands correspondence), as well as our method of handling the combinatorical data by induction on "slopes", could be generalized to broader cases. We hope to address this in our future research.
	

	\begin{example}
		In this example, we provide an explicit construction of the automorphic line bundles used in this paper. We restrict ourselves to the case $G=G\big(\prod_{i=1}^NU(m_i,n_i)\big)$. Let $X$ denote the special fiber of the unitary Shimura variety introduced in Section 2. Recall that the Hodge filtrations on the $\tilde\tau_i$-side and the $\tilde\tau_i^c$-side are dual to each other, so as long as we have made a refinement of the Hodge filtration on the $\tilde{\tau}_i$-component, the refinement on the $\tilde{\tau}_i^c$-component is immediately given by the dual to it. Hence, we can continue to focus on the $\tilde{\tau}_i$-components only. To simplify the notation, we write $\omega_i$ for $\omega_{\mathcal{A}^\vee/X,\tilde\tau_i}$ and $\H(\mathcal{A}/X)_{i}$ for $\H(\mathcal{A}/X)_{\tilde\tau_i}$. By the moduli interpretation, $\omega_i$ is a vector bundle of rank $m_i$. The flag space $Y$ parametrizes tuples $(A,\lambda,\eta,\iota,{\textnormal{Fil}^\bullet})$, where $(A,\lambda,\eta,\iota)\in X(S)$ and $\textnormal{Fil}^\bullet$ is a refinement of the Hodge filtration that respects the $G$-structure. By this, we mean that for each $1\le i\le N$, there is a refinement
		\begin{equation}
			\textnormal{Fil}^i:0=\mathcal{F}^i_0\subseteq \mathcal{F}^i_1\subseteq\cdots\subseteq \mathcal{F}^i_{m_i-1}\subseteq\mathcal{F}^i_{m_i}=\omega_i\subseteq \mathcal{F}^i_{m_i+1}\subseteq\cdots \mathcal{F}^i_{n-1}\subseteq\mathcal{F}^i_n=\H(A/X)_i.
		\end{equation}
	\end{example}
	Here, the superscript labels the component of the subbundle, and the subscript represents the rank of the corresponding subbundle. 
	
	\noindent\textbf{Convention: }In the rest of this paper, 
	we consider the superscripts and subscripts of $\mathcal{F}$ as modulo $N$ and $n$, respectively. We consider the subscripts of $m$ and $n$ as modulo $N$. In other words, we identify $\mathbb{Z}/N\mathbb{Z}$ with $\{1,\dots, N\}$. The convention has the ambiguity that $\mathcal{F}^i_n=\H(\mathcal{A}/X)_0$ but $\mathcal{F}^i_0=0$. However, by Lemma \ref{General Picard Group Relations} below, both are zero in the rational Picard group of $X$. Therefore, this ambiguity here is harmless to our purpose.
	
	By definition, in our unitary setup, a weight for $G$ corresponds to a tuple of integers
	\begin{equation}
		\underline{k}=(k^i_j)^{1\le i\le N}_{1\le j\le n}.
	\end{equation}
	The corresponding automorphic line bundle $\mathcal{L}_Y(\lambda)$ on $Y$ is
	\begin{equation}
		\mathcal{L}_Y(\lambda)=\bigotimes_{i=1}^N\bigotimes_{j=1}^n (\mathcal{F}^i_j/\mathcal{F}^i_{j-1})^{\otimes k^i_j}
	\end{equation}
	The Levi subgroup $L$ of $P\subseteq G_k$ is isomorphic to 
	\begin{equation}
		\bigg(\prod_{i=1}^N(GL_{m_i}\times GL_{n_i})\bigg)\times \mathbb{G}_m.
	\end{equation}
	Therefore, a character of $L$ corresponds to a tuple of integers $\chi=(\alpha_i,\beta_i)_{i=1}^N\in (\mathbb{Z}^2)^N$. The corresponding line bundle on $X$ is
	\begin{equation}
		\mathcal{L}_X(\chi)=\bigotimes_{i=1}^N \big(\textnormal{det }(\H(A/X)_i/\omega_i)\big)^{\otimes \alpha_i}\otimes \big(\textnormal{det }\omega_i\big)^{\otimes \beta_i}.
	\end{equation}
	Since we only care about the positivity of $\mathcal{L}_X(\chi)$, by Lemma \ref{General Picard Group Relations} below, we can replace $\mathcal{L}_X({\chi})$ with 
	\begin{equation}
		\bigotimes_{i=1}^N\big(\textnormal{det }\omega_i\big)^{\otimes (\beta_i-\alpha_i)}
	\end{equation}
	in $(\textnormal{Pic $X$})_{\mathbb{Q}}$. To simplify the notation, we put $k_i=\beta_i-\alpha_i$ and denote
	\begin{equation}
		\mathcal{L}_X(\chi)=\mathcal{L}_X(\{k_i\})=\bigotimes_{i=1}^N\big(\textnormal{det }\omega_i\big)^{\otimes k_i}.
	\end{equation}
	
	The goal of this paper is therefore to provide a numerical criterion for $[\mathcal{L}_Y(\lambda)]$ to be ample.

	\medskip
	\subsection{ Statement of the main theorem.}
	
	\medskip
	\noindent\textbf{Notation: }Let $\mathcal{M}$ be a vector bundle on an algebraic variety $S$. By $[\mathcal{M}]$ we mean the class of the line bundle $(\textnormal{det }\mathcal{M})\in (\textnormal{Pic $S$})_\mathbb{Q}$. Our proof of the ampleness (or nefness) criterion is based on computing the intersection numbers, so we will also write $[\mathcal{M}]$ for $c_1(\mathcal{M})=c_1(\textnormal{det }\mathcal{M})$. It follows easily from the definition that $[\mathcal{M}^{(p)}]=p[\mathcal{M}]$ when $S$ is an $\mathbb{F}_p$-scheme.

	\begin{lem}\label{General Picard Group Relations} Let $(\textnormal{Pic $X$})_{\mathbb{Q}}=(\textnormal{Pic $X$})\otimes_{\mathbb{Z}} \mathbb{Q}$ be the rational Picard group of $X$. We have the following relations in $(\textnormal{Pic $X$})_{\mathbb{Q}}$.
		\begin{enumerate}
			\item $[\omega_{\mathcal{A}^\vee/X,\tilde\tau_i}]=[\omega_{\mathcal{A}^\vee/X,\tilde\tau_i^c}]$. 
			
			\item $[\H(\mathcal{A}/X)_{\tilde\tau_i}]=[\H(\mathcal{A}/X)_{\tilde\tau_i^c}]=0$. 
		\end{enumerate}
		
	\end{lem}
	\noindent\emph{proof:} In the $U(2)$ or Hilbert case, the proof can be found in \cite{Tian-Xiao}. Our proof is a natural generalization to the higher rank case. 
	%
	
	Recall that for each $i$, the perfect pairing $\langle\cdot,\cdot\rangle$ restricts to a perfect pairing
	\begin{equation}
		\langle\cdot,\cdot\rangle_i:\H(\mathcal{A}/X)_{\tilde\tau_i}\times \H(\mathcal{A}/X)_{\tilde\tau_i^c}\longrightarrow \mathcal{O}_X
	\end{equation}
	induced by the polarization and the $\mathcal{O}_E$-action, such that $\omega_{\mathcal{A}^\vee/X,\tilde\tau_i}$ and $\omega_{\mathcal{A}^\vee/X,\tilde\tau_i^c}$ are annihilators of each other. Thus, $\omega^\vee_{\mathcal{A}^\vee/X,\tilde\tau_i}=\H(\mathcal{A}/X)_{\tilde\tau_i^c}/\omega_{\mathcal{A}^\vee/X,\tilde\tau_i^c}$ and we have
	\begin{equation}\label{Equation given by Omega}
		-[\omega_{\mathcal{A}^\vee/X,\tilde\tau_i}]=[\H(\mathcal{A}/X)_{\tilde\tau_i^c}/\omega_{\mathcal{A}^\vee/X,\tilde\tau_i^c}].
	\end{equation} 
	Moreover, this pairing satisfies
	\begin{equation}
		\langle Vx,y\rangle=\langle x, Fy\rangle^{(p)}
	\end{equation}
	for $x\in\H(\mathcal{A}/X)_{\tilde\tau_i}$ and $y\in (\H(\mathcal{A}/X)_{\tilde\tau_{i+1}^c})^{(p)}$. This shows that $(\textnormal{Ker }V_{\tilde\tau_i})^\perp=\textnormal{Ker }V_{\tilde\tau_i^c}$ (because $\textnormal{Im }F=\textnormal{Ker }V$). So we have $\H(\mathcal{A}/X)_{\tilde\tau_i}/\textnormal{Ker }V_{\tilde\tau_i} =(\textnormal{Ker }V_{\tilde\tau_i^c})^\vee$ and
	\begin{equation}	
		[\H(\mathcal{A}/X)_{\tilde\tau_i}/\textnormal{Ker }V_{\tilde\tau_i}]=-[\textnormal{Ker $V_{\tilde\tau_i^c}$}].
	\end{equation} 
	We have the short exact sequence
	\begin{equation}
		0\longrightarrow \textnormal{Ker $V_{\tilde\tau_i}$} \longrightarrow \H(\mathcal{A}/X)_{\tilde\tau_i} \stackrel{V}{\longrightarrow} \omega_{\mathcal{A}^\vee/X,\tilde\tau_{i+1}}^{(p)} \longrightarrow 0.
	\end{equation}
	So,
	\begin{equation}
		[\H(\mathcal{A}/X)_{\tilde\tau_i}]=[\textnormal{Ker }V_{\tilde\tau_i}]+p[\omega_{\mathcal{A}^\vee/X,\tilde\tau_{i+1}}].
	\end{equation}
	A similar computation holds for $\textnormal{Ker }V_{\tilde\tau_i^c}$. This shows
	\begin{equation}\label{Equation given by Ker}
		p[\omega_{\mathcal{A}^\vee/X,\tilde\tau_{i+1}}]=p[\omega_{\mathcal{A}^\vee/X,\tilde\tau_{i+1}^c}]-[\H(\mathcal{A}/X)_{\tilde\tau_i^c}].
	\end{equation}
	Combining the equations (\ref{Equation given by Omega}) (for subscript $i+1$) and (\ref{Equation given by Ker}), we deduce that
	\begin{equation}
		[\H(\mathcal{A}/X)_{\tilde\tau_i^c}]=p[\H(\mathcal{A}/X)_{\tilde\tau_{i+1}^c}].
	\end{equation}
	Letting $i$ run over $1,\dots, N$, we get
	\begin{equation}
		p^{N}[\H(\mathcal{A}/X)_{\tilde\tau_i^c}]=[\H(\mathcal{A}/X)_{\tilde\tau_i^c}].
	\end{equation}
	Thus, $[\H(\mathcal{A}/X)_{\tilde\tau_i^c}]=0$. Since $\H(\mathcal{A}/X)_{\tilde\tau_i}$ is dual to $\H(\mathcal{A}/X)_{\tilde\tau_i^c}$, the above equality also holds for $\H(\mathcal{A}/X)_{\tilde\tau_i}$. Besides, (\ref{Equation given by Omega}) becomes
	\begin{equation}
		[\omega_{\mathcal{A}^\vee/X,\tilde\tau_i}]=[\omega_{\mathcal{A}^\vee/X,\tilde\tau_i^c}].
	\end{equation}
	This completes the proof of the lemma. $\hfill\square$

	\begin{definition}\label{Essential set and essential degree}
		Let $X$ be the special fiber of a unitary Shimura variety of signature $(m_i,n_i)$ as in Section 2. Let $T=\{i|m_i\notin\{0,n\}\}$. We refer to $T$ as the \emph{essential set} of $X$, and $t=\# T$ as the \emph{essential degree} of $X$.
	\end{definition}

	\begin{definition}
		For any $i$, let $a(i)$ be the smallest positive integer such that $i+1,\dots,\allowbreak i+a(i)-1\notin T$, but $i+a(i)\in T$. Note that we do not assume $i\in T$ in this definition.
	\end{definition}

	Since the subscripts of $m$ and $n$ are viewed as modulo $N$, we should also view the elements in $T$ as modulo $N$.

	The following are the main theorems of this paper.

	\begin{thm}[The ampleness criterion for the flag space, case 1]\label{The ampleness criterion for the flag space, case 1}
		If $t=0$, i.e., $T$ is the empty set, then $\mathcal{L}_Y(\lambda)=\mathcal{L}_Y(\{k^i_j\})$ is ample if and only if $\lambda$ is regular dominant, i.e.,
		\begin{equation}
			\begin{aligned}
				k^i_n>k^i_{n-1}>\cdots >k^i_1,\quad \forall 1\le i\le N.
			\end{aligned}
		\end{equation}
	\end{thm}

	\begin{thm}[The ampleness criterion for the flag space, case 2]\label{The ampleness criterion for the flag space, case 2}
		If $t=1$, i.e., $T=\{i_1\}$ is a singleton. By definition $a(i_1)=N$. Assume $m_{i_1}=n-1$, then $\mathcal{L}_Y(\lambda)=\mathcal{L}_Y(\{k^i_j\})$ is ample if and only if
		\begin{equation}
			\begin{aligned}
				&k^{i_1}_{n-1}>k^{i_1}_{n-2}>\cdots>k_1^{i_1},\\
				&k^i_n>k^i_{n-1}>\cdots>k^i_1,\qquad\qquad \textnormal{for $i\ne i_1$},\\
				&\big(\sum_{u=0}^{n-2}p^{Nu}\big)\sum_{j=0}^{N-1}\frac{1}{p^j}(k_1^{i_1+j}-k^{i_1+j}_n)>\sum_{r=2}^{n-1}p^{(r-2)N}\sum_{j=0}^{N-1}\frac{1}{p^j}(k^{i_1+j}_r-k^{i_1+j}_n). 
			\end{aligned}
		\end{equation} 
	\end{thm}
	
	\begin{thm}[The ampleness criterion for the flag space, case 2']\label{The ampleness criterion for the flag space, case 2'}
		If $t=1$, i.e., $T=\{i_1\}$ is a singleton. By definition $a(i_1)=N$. Assume $m_{i_1}=1$, then $\mathcal{L}_Y(\lambda)=\mathcal{L}_Y(\{k^i_j\})$ is ample if and only if
		\begin{equation}
			\begin{aligned}
				&k^{i_1}_n>k^{i_1}_{n-1}>\cdots >k^{i_1}_2,\\
				&k^i_n>k^i_{n-1}>\cdots>k^i_1,\qquad\qquad \textnormal{for $i\ne i_1$},\\
				&	\big(\sum_{u=0}^{n-2}p^{Nu}\big)\sum_{j=0}^{N-1}\frac{1}{p^j}(k_1^{i_1+j}-k^{i_1+j}_n)>\sum_{r=2}^{n-1}p^{(n-1-r)N}\sum_{j=0}^{N-1}\frac{1}{p^j}(k^{i_1+j}_1-k^{i_1+j}_r).
			\end{aligned}
		\end{equation}
	\end{thm}	
	Note that this is indeed the dual version to the previous one.
	

	\begin{thm}[The ampleness criterion for the flag space, case 3]\label{The ampleness criterion for the flag space, case 3}
		Suppose that either $t=\# T\ge 2$, or $t=\# T=1$ but $m_{i_1}\notin\{1,n-1\}$. Let $T=\{i_1,i_2,\dots,i_t\}$ such that $i_l+a(i_l)=i_{l+1}$ for all $l$ (Here, we also view $i_{t+1}=i_1$). Then the line bundle $\mathcal{L}_Y(\lambda)=\mathcal{L}_Y(k^i_j)$ is ample if and only if
		\begin{equation}
			\begin{aligned}
				&k^{i_l}_{m_{i_l}}>k^{i_l}_{m_{i_l}-1}>\cdots>k^{i_l}_{1},\quad k^{i_l}_{n}>k^{i_l}_{n-1}>\cdots>k^{i_l}_{m_{i_l}+1}, \qquad\textnormal{for all $1\le l\le t$},\\
				&k^j_n>k^j_{n-1}>\cdots>k^j_1,\qquad\textnormal{for $j\notin T$},\\
				&p^{a(i_l)}\sum_{j=0}^{a(i_l)-1}\frac{1}{p^j}(k^{i_l+j}_1-k^{i_l+j}_n)>k^{i_{l+1}}_{m_{i_{l+1}}}-k^{i_{l+1}}_{m_{i_{l+1}}+1}\qquad \textnormal{for all $1\le l\le t$.}
			\end{aligned}
		\end{equation}
	\end{thm}
	
	\begin{thm}[The nefness criterion for the flag space]\label{The corresponding nefness criterion}
		In each of the above cases, $\mathcal{L}_Y(\lambda)$ is nef if and only if all ">" are replaced by $\ge$.
	\end{thm}
	

	\noindent\textbf{Variants of the main theorem.} Let $B\subseteq P'\subseteq P$ be an intermediate parabolic subgroup. Recall that we have introduced the partial flag space $Y_{P'}$ and the line bundle $\mathcal{L}_{Y_{P'}}(\lambda)$ associated with a character $\lambda$ of the Levi subgroup $L'$ of $P'$.
	
	If $P'=P$, then $Y_{P'}$ coincides with the special fiber $X$. We have the following.
	
	\begin{thm}[The ampleness criterion for the Shimura variety]\label{The ampleness criterion for Shimura variety} The line bundle
		$$
		\mathcal{L}_X(\{k_l\})=\bigotimes_{i=1}^N\big(\textnormal{det }\omega_i\big)^{\otimes k_i}=\bigotimes_{i\in T}\big(\textnormal{det }\omega_i\big)^{\otimes k_i}
		$$
		is ample (resp. nef) if and only if $p^{a(i)}k_i>k_{i+a(i)}$ (resp. $\ge$) for all $i\in T$.
	\end{thm}
	
	\begin{rmk}
		When $n=2$, the space $X$ corresponds to the special fiber of a $U(2)$ Shimura variety. The ampleness criterion in this case was previously established by the author in \cite{YangAmpleness}.
	\end{rmk}

	If $P'$ is a maximal proper parabolic subgroup of $P$, then $Y_{P'}$ corresponds to refining the Hodge filtration by introducing exactly one new vector bundle $\mathcal{F}^i_j\subseteq \H(\mathcal{A}/X)_i$ for some indices $i$ and $j$. We require that $\mathcal{F}^i_j$ either contains, or is contained in $\omega_i$, depending on the choice of $P'$. We denote this parabolic subgroup $P'$ by $P^i_j$ and the corresponding partial flag space $Y_{P^i_j}$ by $Y^i_j$. We refer to $Y^i_j$ as a \emph{minimal partial flag space}. More generally, the notation $Y^{i_1,\dots,i_r}_{j_1,\dots,j_r}$ refers to the moduli space that records $\{\mathcal{F}^{i_k}_{j_k}\}$'s (with inclusion conditions given by the parabolic subgroup) for all $1\le k\le r$. The type of $P^{i_1,\dots,i_r}_{j_1,\dots,j_r}$ is $I^{i_1,\dots,i_r}_{j_1,\dots,j_r}$. These intermediate partial flag spaces will be intensively used in the subsequent proofs.

	\begin{definition}\label{Flagged at some place}
		Let $P'\subseteq P$ be an intermediate parabolic subgroup and $Y_{P'}$ the corresponding partial flag space. If $Y_{P'}$ records a refinement of the Hodge filtration at $\tau_i$, we say that $Y_{P'}$ is \emph{flagged} at $i$.
	\end{definition}
	
	\begin{thm}[The nefness criterion for partial flag space]\label{The nefness criterion for partial flag space} Let $P'\subseteq P$ be an intermediate parabolic subgroup, and $\lambda$ be a character of $P'$. Then $\mathcal{L}_{Y_{P'}}(\lambda)$ is nef if and only if its pullback to $Y_B=Y$ is nef.	
	\end{thm}
	
	Our idea of proving the sufficiency part of the ampleness criterion is to reduce the problem to the nefness criterion for \emph{minimal} partial flag spaces. For convenience, we present the explicit inequalities below. Readers can skip this part first, and come back when reading the proof, especially the induction steps, in Section 7.
	
	Let $P^i_j$ be a maximal proper parabolic subgroup of $P$, and $L^i_j$ a Levi subgroup of $P^i_j$. A character of the Levi subgroup $L^i_j$ corresponds to a tuple of integers $(\{k_l\};\alpha)$. The associated automorphic line bundle on $Y^i_j$ can be written as:
	\begin{equation}
		\mathcal{L}_{Y^i_j}(\{k_l\};\alpha)=\mathcal{L}_{Y^i_j}(k_1,\dots,k_N;\alpha)= \big(\textnormal{det }(\omega_i/\mathcal{F}^i_j)\big)^{\otimes k_i}\otimes (\textnormal{det }\mathcal{F}^i_j)^{\otimes \alpha} \otimes \bigotimes_{\substack{s=1\\ s\ne i}}^N (\textnormal{det $\omega_s$})^{\otimes k_s}.
	\end{equation}
	Here and later, if $j>m_i$, then $\omega_i$ is contained in $\mathcal{F}^i_j$, yet we still write $\omega_i/\mathcal{F}^i_j$ to mean $(\mathcal{F}^i_j/\omega_i)^{\vee}$. 
	

	\begin{coro}[The nefness criterion for minimal partial space]\label{The nefness criterion for minimal partial flag space} We write $T=\{i_1,\dots,i_t\}$. Assumptions as above, we assume further that $P'=P^i_j$ is a maximal proper parabolic subgroup. 
		
		\begin{enumerate}
			\item  If $T$ is the empty set, then $\mathcal{L}_{Y^i_j}(\{k_l\};\alpha)$ is nef if and only if
			\begin{equation}
				k_i\ge\alpha.
			\end{equation}
			
			From now on we assume that $T$ is non-empty, then there is a unique $1\le l_0\le t$, such that $i\in \{i_{l_0},i_{l_0}+1,\dots, i_{l_0+1}-1\}$. For simplicity, we assume that $i_1\le i$.
			
			\item If $T$ is a singleton, namely, $T=\{i_1\}$, then $l_0=1$ by definition. If moreover $m_{i_1}=n-1$, then $\mathcal{L}_{Y^i_j}(\{k_l\};\alpha)$ is nef if and only if
			\begin{equation}
				\begin{aligned}
					&k_i\ge\alpha,\\
					&(p^{(n-j)N}-1)\big(k_{i_1}-\frac{1}{p^{i-i_1}}(k_i-\alpha)\big)\ge (p^{(n-j-1)N}-1)k_{i_1}.
				\end{aligned}
			\end{equation}
			
			\item If $T$ is a singleton, namely, $T=\{i_1\}$, then $l_0=0$ by definition. If moreover $m_{i_1}=1$, then $\mathcal{L}_{Y^i_j}(\{k_l\};\alpha)$ is nef if and only if
			\begin{equation}
				\begin{aligned}
					&k_i\ge \alpha,\\
					&(p^{jN}-1)\big(k_{i_1}-\frac{1}{p^{i-i_1}}(k_i-\alpha)\big)\ge (p^{(j-1)N}-1)k_{i_1}.
				\end{aligned}
			\end{equation}
			
			\item In the rest cases, i.e., either $t\ge2$, or $t=1$ but $m_{i_1}\notin\{1,n-1\}$. Then $\mathcal{L}_{Y^i_j}(\{k_l\};\alpha)$ is nef if and only if
			\begin{equation}
				\begin{aligned}
					&k_i\ge \alpha,\\
					&p^{a(i_j)}k_{i_j}\ge k_{i_{j+1}}, \qquad\textnormal{for $j\ne l_0$},\\
					&p^{a(i_{l_0})}\big(k_{i_{l_0}}-\frac{1}{p^{i-i_{l_0}}}(k_i-\alpha)\big)\ge k_{i_{l_0}+1}.
				\end{aligned}
			\end{equation}
		\end{enumerate}
	\end{coro}
	
	The author encounters a difficulty in proving the necessity part of the ampleness criterion on partial flag spaces, since we cannot pull ample line bundles back via proper morphisms. Nevertheless, the nefness criterion for partial flag spaces suffices for our purpose of understanding the ampleness criterion of $\mathcal{L}_Y(\lambda)$.

	\section{Calculus on the special fiber}
	
	The primary approach to the ampleness criterion involves doing intersection theory on the special fiber and the associated (partial) flag spaces. However, when dealing with certain low-dimensional strata, it turns out that neither the theory of $G$-zips nor the theory of Schubert calculus on flag varieties provide sufficiently many divisor classes. A crucial step in addressing this problem is to establish correspondences between the strata of different Shimura varieties. This enables us to relate automorphic line bundles for different unitary Shimura varieties and make suitable inductive arguments.
	
	This idea was initially introduced in Helm's work \cite{Helm} for sparse strata on $U(2)$ Shimura varieties. It was subsequently extended by \cite{Tian-Xiao} to $U(2)$ and quaternionic Shimura varieties, and further generalized by \cite{Helm-Tian-Xiao,Xiao-Zhu} to a broader framework. In this chapter, we present a simplified version tailored to our specific needs.

	\medskip
	\subsection{ Review of Dieudonn\'e theory}
	
	In this section, let $k$ be a perfect field of characteristic $p$, and let $W(k)$ denote its ring of Witt vectors. The Frobenius endomorphism $\sigma$ on $W(k)$ lifts the Frobenius map on the residue field $k$. The classical (covariant) Dieudonn\'e theory establishes an equivalence
	\begin{equation}
		\Bigg\{
		\begin{aligned}
			&p\textnormal{-divisible groups} \\
			&\textnormal{of finite height over $k$}
		\end{aligned}
		\Bigg\}
		\stackrel{\simeq}{\longrightarrow}
		\Bigg\{
		\begin{aligned}
			&W(k)\textnormal{-modules free of finite rank equipped with operators $F,V$}\\
			&\textnormal{such that $F$ is $\sigma$-linear, $V$ is $\sigma^{-1}$-linear, and $FV=VF=p$}.
		\end{aligned}
		\Bigg\}
	\end{equation}
	For $p$-torsion groups, we have
	\begin{equation}
		\Bigg\{
		\begin{aligned}
			\textnormal{finite } p\textnormal{-torsion groups over $k$}
		\end{aligned}
		\Bigg\}
		\stackrel{\simeq}{\longrightarrow}
		\Bigg\{
		\begin{aligned}
			&\textnormal{finite dimensional }k\textnormal{-vector spaces equipped with operators $F,V$}\\
			&\textnormal{such that $F$ is $\sigma$-linear, $V$ is $\sigma^{-1}$-linear, and $FV=VF=0$}.
		\end{aligned}
		\Bigg\}.
	\end{equation}
	
	Let $D$ denote the (covariant) functor sending a $p$-divisible group (resp. a $p$-torsion group) over $k$ to its associated $W(k)$-module (resp. $k$-vector space) in the target category. The functor $D$ satisfies the relation
	\begin{equation}
		D(G[p])=D(G)/pD(G).
	\end{equation}
	Here, $G[p]$ denotes the kernel of the multiplication-by-$p$ morphism on the $p$-divisible group $G$.
	\begin{rmk}
		In general, a finite $p$-torsion group over $k$ need not arise as the $p$-torsion subgroup of a $p$-divisible group. Specifically, the subcategory consisting of "finite groups which are the $p$-torsion of a $p$-divisible group" corresponds to $k$-vectors spaces in the target category which satisfies $\textnormal{Ker }F=\textnormal{Im }V$, $\textnormal{Ker} V=\textnormal{Im }F$. See \cite{MoonenGroupsSchemes} for an explanation.
	\end{rmk}
	
	Let $A/k$ be an abelian variety. Let $\tilde{D}(A)=D(A[p^{\infty}])$ denote the Dieudonn\'e module of the $p$-divisible group associated to $A$ and let $D(A)=D(A[p])$. There is a natural comparison isomorphism
	\begin{equation}
		\tilde{D}(A)/p\tilde{D}(A)\simeq D(A)\simeq \H(A/k).
	\end{equation}
	Recall that in our moduli problem, $A$ is equipped with an action of $\mathcal{O}_E$. This induces a corresponding action of $\mathcal{O}_E$ on the Dieudonn\'e modules $\tilde{D}(A)$ and ${D}(A)$, leading to the natural decomposition
	\begin{equation}
		\tilde{D}(A)=\bigoplus_{1\le i\le N}\big(\tilde{D}(A)_{\tilde\tau_i}\oplus \tilde{D}(A)_{\tilde\tau_i^c}\big).
	\end{equation}
	  The relative Frobenius and Verschiebung maps on $A/k$ induce $W(k)$-linear maps $V:\tilde{D}(A)\to \tilde{D}(A^{(p)})$ and $F:\tilde{D}(A^{(p)})\to \tilde{D}(A)$ respectively, which further induces maps $V$ and $F$ (note the order!) on each component:
	\begin{equation}
		V:\tilde{D}(A)_{\tilde\tau_i}\longrightarrow \tilde{D}(A)_{\tilde\tau_{i+1}},\qquad F:\tilde{D}(A)_{\tilde\tau_{i+1}}\longrightarrow \tilde{D}(A)_{\tilde\tau_i}.
	\end{equation}
	A similar decomposition holds for the $\tilde\tau_i^c$-components and for the Dieudonn\'e module $D(A)$. Furthermore, $A$ is equipped with a prime-to-$p$ quasi-polarization $\lambda$, which induces a perfect pairing
	\begin{equation}
		\langle \tilde{D}(A)_{\tilde\tau_i},\tilde{D}(A)_{\tilde\tau_i^c}\rangle_{\lambda}\longrightarrow W(k).
	\end{equation}
	This pairing satisfies $\langle Fx,y\rangle_{\lambda}^{\sigma}=\langle x,Vy\rangle_{\lambda}$ for $x\in \tilde{D}(A)_{\tilde\tau_i}$ and $y\in\tilde{D}(A)_{\tilde\tau_{i+1}^c}$. The subscript $\lambda$ is omitted when it is clear in the context. For any free $W(k)$-submodule $M\subseteq \tilde{D}(A)_{\tilde\tau_i}$, we define
	\begin{equation}
		M^\perp:=\big\{y\in \tilde{D}(A)_{\tilde\tau_i^c}|\ \langle x,y\rangle\in pW(k) \ \textnormal{for all }x\in M \big\}.
	\end{equation}
	By definition, $\omega_{\tilde\tau_i}$ and $\omega_{\tilde\tau_i^c}$ are orthogonal complements of each other. As in Section 3, we simplify our notations by writing $\tilde{D}(A)_i$ and ${D}(A)_i$ for $\tilde{D}(A)_{\tilde\tau_i}$ and $D(A)_{\tilde\tau_i}$ respectively.
	
	An isogeny $f:A\rightarrow B$ induces a homomorphism $\tilde{D}(A)\rightarrow \tilde{D}(B)$ and also for each $\tilde\tau_i$- and $\tilde\tau_i^c$-component, compatible with $F$ and $V$. When Ker $f\subseteq A[p]$, this induces $p\tilde{D}(B)\subseteq f(\tilde{D}(A))\subseteq \tilde{D}(B)$. Alternatively, this can be interpreted as an inclusion of lattices $\tilde{D}(A)\subseteq \tilde{D}(B)$ in the $W(k)[\frac{1}{p}]$-vector space $\tilde{D}(A)[\frac{1}{p}]\allowbreak =\tilde{D}(A)\otimes_{W(k)}W(k)[\frac{1}{p}]\allowbreak \simeq\tilde{D}(B)[\frac{1}{p}]$. The dual isogeny $g$ of $f$, then induces an inclusion $g:\tilde{D}(B)\hookrightarrow \frac{1}{p}\tilde{D}(A)$.

	\medskip
	\subsection{ Geometric description of strata}
	
	The following theorem, which corresponds to the $m=1$ case in \cite{Helm-Tian-Xiao}, will provide a tool for constructing auxiliary Shimura varieties.
	
	\begin{thm}\label{Construction via F,V-chain}
		Let $k$ be a perfect field containing the residue field $\mathcal{O}_F/(p)$, and let $(A,\lambda,\eta,\iota)\in X(k)$ be a $k$-point. Suppose that we are given a free $W(k)$-submodule $\tilde{\mathcal{E}}_i\subseteq \tilde{D}(A)_{i}$ for each $1\le i\le N$, which satisfies
		\begin{equation}\label{F,V-chain condition}
			p\tilde{D}(A)_{i}\subseteq \tilde{\mathcal{E}}_i\subseteq \tilde{D}(A)_{i},\quad V(\tilde{\mathcal{E}}_i)\subseteq \tilde{\mathcal{E}}_{i+1},\quad \textnormal{and}\quad F(\tilde{\mathcal{E}}_{i+1})\subseteq \tilde{\mathcal{E}}_i
		\end{equation}
		for all $i$. Then there exists a unique pair $(B,\phi)$ satisfying
		
		\medskip
		\textnormal{(1)} $B$ is an abelian variety over $k$, endowed with an $\mathcal{O}_E$-action $\iota_B:\mathcal{O}_E\to \textnormal{End}(B)$ and a prime-to-$p$ polarization $\lambda_B$, such that the Rosati involution associated with $\lambda_B$ induces the complex conjugation on $\mathcal{O}_E$.
		
		\textnormal{(2)} $\phi:B\to A$ is an $\mathcal{O}_E$-equivariant isogeny, such that for $1\le i\le N$, the inclusion $\tilde{\mathcal{E}}_i\hookrightarrow \tilde{D}(A)_i$ is identified with the map $\phi_{\ast,i}:\tilde{D}(B)_i\rightarrow \tilde{D}(A)_i$ induced by $\phi$, and $\phi^\vee\circ \lambda\circ \phi=p\lambda_B$. 
		
		\textnormal{(3)} If $\textnormal{length}_{W(k)}(\tilde{\mathcal{E}}_i/p\tilde{D}(A)_i)=l_i$ for $1\le i\le N$, then
		\begin{equation}
			\textnormal{dim }\omega_{B^\vee/k,i}=m_i-l_i+l_{i-1}.
		\end{equation}
		
		\textnormal{(4)} The abelian variety $B/k$ is equipped with a prime-to-$p$ level structure $\eta_B$ on $B$ such that $\eta=\phi\circ\eta_B$.
	\end{thm}
	\noindent\emph{Proof:} By our previous discussion, the polarization $\lambda$ induces a perfect pairing $\langle\cdot,\cdot\rangle_{\lambda}:\allowbreak\tilde{D}(A)_i\times\tilde{D}(A)_{i^c}\allowbreak\to W(k)$ for all $1\le i\le N$. Let $\tilde{E}_{i^c}:=\tilde{E}_i^{\perp}$ be the orthogonal complement of $\tilde{\mathcal{E}}_i$. One can verify that
	\begin{equation}
		p\tilde{D}(A)_{\tilde\tau_i^c}\subseteq \tilde{\mathcal{E}}_{i^c}\subseteq \tilde{D}(A)_{\tilde\tau_i^c},\quad V(\tilde{\mathcal{E}}_{i^c})\subseteq \tilde{\mathcal{E}}_{{(i+1)}^c},\quad \textnormal{and}\quad F(\tilde{\mathcal{E}}_{(i+1)^c})\subseteq \tilde{\mathcal{E}}_{i^c}.
	\end{equation}
	Thus, the submodule
	\begin{equation}
		\tilde{D}:=\bigoplus_{i=1}^N \big(\tilde{\mathcal{E}}_i\oplus \tilde{\mathcal{E}}_{i^c}\big)
	\end{equation}
	is a free $W(k)$-submodule of $\tilde{D}(A)$ equipped with an action of $F$ and $V$. By Dieudonn\'e theory, the quotient $\tilde{D}/p\tilde{D}(A)$ corresponds to a closed subgroup scheme $H\subseteq A[p]$. Moreover, it follows from the above construction that $H=H^\perp$, meaning that $H$ is maximal isotropic. Let $B:=A/H$, and let $\psi:A\to B$ be the natural quotient. We define $\phi:B\to A$ as the quotient with kernel $\psi(A[p])$. We now verify the properties of $(B,\phi)$.
	
	(1) The $\mathcal{O}_E$-action follows from the construction. Since $H\subseteq A[p]$ is maximal totally isotropic, it follows from \cite{MumfordAV} that there is a prime-to-$p$ polarization $\lambda_B$ on $B$ with $\phi^\vee\circ\lambda\circ\phi=p\lambda_B$.
	
	(2) By construction of $\phi$, The Dieudonn\'e module $\tilde{D}(B)$ is identified with $\tilde{D}$, and this identification holds componentwise for each $\tau_i$.
	
	(3) We have
	\begin{equation}
		\begin{aligned}
			\textnormal{dim}_k\ \omega_{B^\vee/k,i}&=\textnormal{dim}_k\ \frac{V\tilde{D}(B)_{i-1}}{p\tilde{D}(B)_i}=\textnormal{dim}_k\ \frac{V\tilde{\mathcal{E}}_{i-1}}{p\tilde{\mathcal{E}}_i}\\
			&=\textnormal{dim}_k\ \frac{V\tilde{\mathcal{E}}_{i-1}}{Vp\tilde{D}(A)_{i-1}}+\textnormal{dim}_k\ \frac{Vp\tilde{D}(A)_{i-1}}{p^2\tilde{D}(A)_i}-\textnormal{dim}_k\ \frac{p\tilde{\mathcal{E}}_{i}}{p^2\tilde{D}(A)_i}\\
			&=\textnormal{dim}_k\frac{\tilde{\mathcal{E}}_{i-1}}{p\tilde{D}(A)_{i-1}}-\textnormal{dim}_k\frac{\tilde{\mathcal{E}}_i}{p\tilde{D}(A)_i}+\textnormal{dim}_k\ \omega_i\\
			&=m_i-l_i+l_{i-1}.
		\end{aligned}
	\end{equation}
	
	(4) Since $\phi$ is a $p$-isogeny, the prime-to-$p$ level structure $\eta_B$ on $B$ is uniquely determined by $\eta$ on $A$.$\hfill\square$
	
	\begin{definition}
		Let $A$ be an abelian scheme over $k$ and $\{\tilde{\mathcal{E}}_i\}_{i=1}^N$ be a tuple of free $W(k)$-submodules of $\tilde{D}(A)$. We call it an \emph{$F,V$-chain in $\tilde{D}(A)$} if it satisfies condition (\ref{F,V-chain condition}).
		
		Similarly, let $\{\mathcal{E}_i\}_{i=1}^N$ be a tuple of subspaces of $D(A)$. We call it an \emph{$F,V$-chain in $D(A)$} if there exists an $F,V$-chain in $\tilde{D}(A)$, namely, $\{\tilde{\mathcal{E}}_i\}$, such that $\mathcal{E}_i=\allowbreak\tilde{\mathcal{E}}_i/p\tilde{D}(A)_i\allowbreak\subseteq\H(A/k)_i$ for $1\le i\le N$.
	\end{definition}
	
	
	\begin{definition}
		A morphism of $k$-schemes $f:Y\to Z$ is called a \emph{Frobenius factor} if there exists $r\in \mathbb{N}$ and a morphism $g:Z^{(p^r)}\to Y$, such that $g\circ f:Z^{(p^r)}\to Z$ is the $p^r$-th power Frobenius morphism.
	\end{definition}
	
	\begin{lem}\label{Frobenius twist morphisms}
		Let $f:Y\to Z$ be a proper morphism between $k$-schemes of finite type, such that $Y$ is reduced and $Z$ is normal. If $f$ is a bijection on closed points, then $f$ is a Frobenius factor.
	\end{lem}
	\noindent\emph{Proof.} This is \textup{\cite[Prop 4.8]{Helm}}. $\hfill\square$
	
	In the following theorem, we construct correspondences between (strata of partial flag spaces over) Shimura varieties.
	
	\begin{thm}\label{Description of strata}
		Let $B\subseteq P_0\subseteq P$ be an intermediate parabolic subgroup such that the associated partial flag space $Y_{P_0}$ is defined by recording at most one subbundle $\mathcal{E}_{r_i}^i\subseteq \omega_i$ of rank $r_i$ for each $i$ ($r_i$ can be $m_i$, in which case the moduli problem does not record any refinement at the $\tilde\tau_i$-component). Let $Z$ be the closed subscheme of $Y_{P_0}$ defined by the condition
		\begin{equation}
			\mathcal{E}^i_{r_i}\subseteq V^{-1}(\mathcal{E}_{r_{i+1}}^{i+1,(p)}),\qquad F(\mathcal{E}_{r_{i+1}}^{i+1,(p)})\subseteq \mathcal{E}^i_{r_i}
		\end{equation}
		for all $1\le i\le N$ (Here, the second condition is in fact empty, since $\mathcal{E}^{i+1,(p)}_{r_{i+1}}\subseteq \omega^{(p)}_{i+1}=\textnormal{Ker }F$, but we keep it for uniformity). Then $Z$ is, up to some Frobenius factor, isomorphic to a stratum $Z'$ in a partial flag space $Y'_{P_0'}$ over $X'$, where $X'$ is the special fiber of a unitary Shimura variety with hyperspecial level structure and signature $(m_i',n_i')_{i=1}^N$, where
		\begin{equation}
			m_i'=m_{i}-r_{i}+r_{i-1},\quad n_i'=n-m_i',\quad \textnormal{for $1\le i\le N$}.
		\end{equation}
		Moreover, this correspondence respects additional structures.
	\end{thm}
	\begin{rmk} The precise meaning of "respects additional structures" will be clarified in the proof. Specifically, if there is another intermediate parabolic subgroup $B\subseteq P_1\subseteq P_0\subseteq P$ such that $Y_{P_1}$ is given by further recording some $\mathcal{E}^u_v$ 's that refine the Hodge filtration in the moduli problem, and $Z_1$ is the pullback of $Z$ under the natural projection, then $Z_1$ is, up to some Frobenius factor, isomorphic to a stratum $Z_1'$ in $Y'_{P_1'}$ over $Z'\subseteq Y'_{P_0'}$ for some parabolic subgroup $B'\subseteq P_1'\subseteq P_0'\subseteq P'$.
	\end{rmk}

	\noindent\emph{Proof.} It is easy to check that $\sum_{i=1}^{N}m_i=\sum_{i=1}^{N}m_i'$.
	By \cite[Corollary 8.2]{Helm}, there exists an $n$-dimensional vector space $V'$ over $E$ equipped with a symplectic pairing $\langle\cdot,\cdot\rangle':V'\times V'\to \mathbb{Q}$ satisfying
	\begin{equation}
		\qquad \langle \alpha x,y\rangle' =\langle x,\bar{\alpha}y \rangle',\qquad \forall x,y\in V',\ \alpha\in E,
	\end{equation}
	such that $(V,\langle\cdot,\cdot\rangle)$ is isomorphic to $(V',\langle\cdot,\cdot\rangle')$ at all finite places, and for $1\le i\le N$, the signature of $\langle\cdot,\cdot\rangle'$ at the infinite place $\tau_i$ of $F$ is $(m_i',n_i')$. Let $G'$ be the similitude unitary group of $(V,\langle\cdot,\cdot\rangle')$. We have
	\begin{equation}
		G'(\mathbb{R})=G\big(\prod_{i=1}^{N} U(m_i',n_i')\big),
	\end{equation}
	and an isomorphism
	\begin{equation}
		G(\mathbb{A}_f)\simeq G'(\mathbb{A}_f).
	\end{equation}
	Let $K'^p$ be the compact open subgroup of $G'(\mathbb{A}_f^{(p)})$ corresponding to $K^p$ under this isomorphism, and $K'_p$ be a hyperspecial subgroup at $p$. Let $X'$ be the special fiber of the Shimura variety parametrizing tuples $(B,\lambda_B,\eta_B,\iota_B)\in X'(S)$, where
	
	(1) $B$ is an abelian scheme over $S$ of dimension $nN$ equipped with an $\mathcal{O}_E$-action $\iota_B:\mathcal{O}_E\hookrightarrow \textnormal{End}_SB$, satisfying the Kottwitz signature condition given by $(m_i',n_i')$.
	
	(2) $\lambda_B:B\to B^\vee$ is a prime-to-$p$ polarization, and the associated Rosati involution induces the complex conjugation on $\mathcal{O}_E$.
	
	(3) $\eta_B$ is an $K'^p$-level structure.
	
	We denote by $\mathcal{B}/X'$ the universal abelian scheme. Let $P'$ be the Hodge parabolic of $G'$ and $B'\subseteq P'$ be a Borel subgroup. We consider an intermediate parabolic subgroup $B'\subseteq P'_0\subseteq P'$ such that the associated partial flag space $Y'_{P'_0}$ over $X'$ is defined by specifying an additional vector subbundle
	\begin{equation}
		0\subseteq\mathcal{F}^i_{m_i-r_i}\subseteq\omega_{\mathcal{B}^\vee/X',i}
	\end{equation}
	of rank $m_i-r_i$ for each $1\le i\le N$. Let $Z'\subseteq Y'_{P_0'}$ denote the closed subscheme defined by the Schubert condition
	\begin{equation}
		V^{-1}(\mathcal{F}^{i+1,(p)}_{m_{i+1}-r_{i+1}})\subseteq\mathcal{F}^i_{m_i-r_i} ,\qquad \textnormal{for all $1\le i\le N$}.
	\end{equation}
	We further define a moduli scheme $W$, such that $W(S)$ parametrizes the equivalence classes of tuples $(A,\lambda,\eta,\allowbreak\iota,\allowbreak\textnormal{Fil}^\bullet,B,\allowbreak\lambda_B,\allowbreak\eta_B,\allowbreak\iota_B,\allowbreak\textnormal{Fil}'^\bullet,\phi)$, where
	
	(1) $(A,\lambda,\eta,\iota,\textnormal{Fil}^\bullet)\in Z(S)$;
	
	(2) $(B,\lambda_B,\eta_B,\iota_B,\textnormal{Fil}'^{\bullet})\in Z'(S)$;
	
	(3) $\phi:B\to A$ is an isogeny whose kernel is contained in $B[p]$, and the induced morphism $\phi_{\ast,i}$ on de Rham homology groups satisfies
	\begin{equation}
		\textnormal{Ker }\phi_{\ast,i}=V^{-1}(\mathcal{F}^{i+1,(p)}_{m_{i+1}-r_{i+1}}),\qquad \phi_{\ast,i}(\H(B/S)_i)=\mathcal{E}^{i}_{r_i},\qquad \textnormal{for all $1\le i\le N$.}
	\end{equation}
	Here $\phi_{\ast,i}$ stands for the induced morphism on $\tilde\tau_i$-components.
	
	The forgetful maps induce the following diagram
	\begin{center}
		\begin{tikzpicture}[node distance=10pt and 1cm, auto]
			\node(1) [align=center] {$W$};
			\node(2) [align=center, below left=of 1] {$Z$};
			\node(3) [align=center, below right=of 1] {$Z'$};
			\node(4) [align=center, left=of 2] {$Y_{P_0}$};
			\node(5) [align=center, right=of 3] {$Y'_{p_0'}.$};
			
			\draw [->] (1) -- node[midway, above] {$pr_1$} (2);
			\draw [->] (1) -- node[midway, above] {$pr_2$} (3);
			\draw [<-] (4) -- node[midway, above] {$\supseteq$} (2);
			\draw [->] (3) -- node[midway, above] {$\subseteq$} (5);
		\end{tikzpicture}
	\end{center}
	
	Next, we check the properties for $pr_1$ and $pr_2$.
	
	\medskip
	\noindent\emph{$pr_1$ is bijective on closed points:} It suffices to construct an inverse of $pr_1$ on closed points. Given $(A,\lambda,\eta,\iota,\textnormal{Fil}^\bullet)\in Z(k)$, by the condition of $\textnormal{Fil}^\bullet$ on $Z(k)$, the tuple $\{\mathcal{E}^i_{r_i}\}$ forms an $F,V$-chain. By Theorem \ref{Construction via F,V-chain}, there exists a unique tuple $(B,\lambda_B,\eta_B,\iota_B)\in X'(k)$ and an isogeny $\phi:B\to A$. For each $i$, let $\tilde{\mathcal{E}}^i$ be the inverse image of $\mathcal{E}^i_{r_i}\subseteq \tilde{D}(A)_i/p\tilde{D}(A)_i$ in $\tilde{D}(A)_i$. By the construction in the proof of Theorem \ref{Construction via F,V-chain}, we have $\tilde{\mathcal{E}}^{i}/p\tilde{\mathcal{E}}^{i}\simeq \H(B/k)_i$. Let $\tilde{\omega}_i$ denote the inverse image of $\omega_i\subseteq \tilde{D}(A)_i/p\tilde{D}(A)_i$ in $\tilde{D}(A)_i$. Then the image of $\tilde\omega_i$ in $\tilde{\mathcal{E}}^{i}/p\tilde{\mathcal{E}}^{i}$ is a vector space of dimension $m_i-r_i$, and we denote it by $\mathcal{F}^i_{m_i-r_i}$. Note that $\mathcal{F}^i_{m_i-r_i}$ is contained in $\omega_{B^\vee/k,i}$. Let 
	\begin{equation}
		\textnormal{Fil}'^i:0\subseteq \mathcal{F}^i_{m_i-r_i}\subseteq \omega_{B^\vee/k,i}\subseteq \H(B/k)_{i}
	\end{equation} 
	be a refinement of Hodge filtration for the $\tilde\tau_i$-component, and let $\textnormal{Fil}'^\bullet=\{\textnormal{Fil}'^i\}_{i=1}^N$. It is straightforward to check that $(A,\lambda,\eta,\iota,\textnormal{Fil}^\bullet,B,\lambda_B,\eta_B,\iota_B,\textnormal{Fil}'^\bullet,\phi)$ corresponds to a point in $W(k)$.
	
	\medskip
	\noindent\emph{$pr_2$ is bijective on closed points:} We can construct the inverse of $pr_2$ in a same manner as above, by considering the $F,V$-chain $\{V^{-1}(\mathcal{F}^{i+1,(p)}_{m_{i+1}-r_{i+1}})\}_{i=1}^N$. We omit the remaining detailed arguments.
	
	\medskip
	Next, we show the smoothness of $Z$ and $Z'$ by computing their tangent bundles.
	
	\medskip	
	\noindent\emph{Smoothness of $Z$:} Let $\mathbb{I}=k[\epsilon]/(\epsilon^2)$. Given a point $x\in Z(k)$ corresponding to $(A,\lambda,\eta,\iota,\textnormal{Fil}^\bullet)$, we must determine all possible liftings $x_{\mathbb{I}}\in Z(\mathbb{I})$ of $x$. Recall the following theorem \cite[pp.116-118]{GrothendieckCrystallineDieudonne}\cite[Chapter V, Theorem 2.3]{MessingCrystals}\cite{Mazur-Messing}.
	\begin{thm}[Grothendieck-Messing-Serre-Tate]\label{Deformation of Abelian Schemes}
		Let $S$ be a nilpotent thickening of a scheme $\overline{S}$ with a divided-power structure. Let $AV_{S}$ be the category of abelian schemes over $S'$. Let $AV_{\overline{S}}^+$ be the category of pairs $(\overline{A},\omega)$, where $\overline{A}$ is an abelian scheme over $\overline{S}$, and $\omega\subseteq H_1^{\textnormal{cris}}(\overline{A}/\overline{S})_S$ is a subbundle lifting $\omega_{\overline{A}^\vee/\overline{S}}\subseteq \H(\overline{A}/\overline{S})$. Then the natural functor
		\begin{equation}
			\begin{aligned}
				AV_{S}&\longrightarrow AV^+_{\overline{S}}\\
				A&\mapsto (\overline{A}=A\times_S\overline{S},\omega_{A^\vee/S})
			\end{aligned}
		\end{equation}
		is an equivalence of categories.
	\end{thm}
	
	Thus, to construct a lifting $A_\mathbb{I}$ of $A$, it suffices to construct a subbundle $\omega_{A^\vee_{\mathbb{I}}/\mathbb{I}}$ that lifts $\omega_{A^\vee/k}$. To lift the $\mathcal{O}_E$-action, we must take $\omega_{A^\vee_{\mathbb{I}}/\mathbb{I}}=\bigoplus_{i=1}^N(\omega_{A^\vee_{\mathbb{I}}/\mathbb{I},\tilde\tau_i}\oplus \omega_{A^\vee_{\mathbb{I}}/\mathbb{I},\tilde\tau_i^c})$ such that the lifting is componentwise. To lift the polarization $\lambda$, the lifting must satisfy $\omega_{A^\vee_{\mathbb{I}}/\mathbb{I},\tilde\tau_i^c}=\omega_{A^\vee_{\mathbb{I}}/\mathbb{I},\tilde\tau_i}^\perp$ under the pairing induced by $\lambda$. To lift $\textnormal{Fil}^\bullet$, it suffices to choose subbundles $0\subseteq \mathcal{E}^i_{r_i,\mathbb{I}}\subseteq \omega_{A^\vee_{\mathbb{I}}/\mathbb{I},\tilde\tau_i}$ lifting $\mathcal{E}_{r_i}^i$ for each $i$, such that $\mathcal{E}^i_{r_i,\mathbb{I}}\subseteq V^{-1}(\mathcal{E}^{i+1,(p)}_{r_{i+1},\mathbb{I}})$ for all $1\le i\le N$. 
	
	Note that the $p$-th power map on $k[\epsilon]/(\epsilon^2)$ takes $\epsilon$ to 0. Thus, $V^{-1}(\mathcal{E}^{i+1,(p)}_{r_{i+1},\mathbb{I}})$ is independent of the choice of lifting $\mathcal{E}^{i+1}_{r_{i+1}}$. In other words, for each $i$, it suffices to first lift $\mathcal{E}^{i}_{r_i}$ to a sub-bundle of the pre-fixed bundle $V^{-1}(\mathcal{E}^{i+1,(p)}_{r_{i+1}})$ of rank $r_i$, and then lift $\omega_{A^\vee/k,i}$ to a vector bundle containing $\mathcal{E}^i_{r_i}$ and contained in $H^{\textnormal{cris}}_1(A/K)_{\mathbb{I}}$. We obtain
	\begin{equation}
		0
		\longrightarrow \bigoplus_{i=1}^N \mathcal{H}om\big(\mathcal{E}^i_{r_i},V^{-1}(\mathcal{E}^{i+1,(p)}_{r_{i+1}})/\mathcal{E}^i_{r_i}\big)  \longrightarrow{\mathcal{T}}_{Z/k} 
		\longrightarrow \bigoplus_{i=1}^N \mathcal{H}om(\omega_{\mathcal{A}^\vee/Z,i}/\mathcal{E}^i_{r_i},\H(\mathcal{A}/Z)_{\tilde\tau_i}/\omega_{\mathcal{A}^\vee/Z,i}\big) \longrightarrow0.
	\end{equation}
	This calculation of the tangent bundle $\mathcal{T}_{Z/k}$ demonstrates that $Z$ is smooth.
	
	\medskip
	\noindent\emph{Smoothness of $Z'$:} The analysis of $\mathcal{T}_{Z'/k}$ is analogous.
	
	\medskip
	\noindent\emph{$pr_2$ is a bijection on tangent spaces:} Let $y=(A,\lambda,\eta,\iota,\textnormal{Fil}^\bullet,B,\lambda_B,\eta_B,\iota_B,\textnormal{Fil}'^\bullet,\phi)\in W(k)$ be the inverse image of $x$. Given a lift $x_{\mathbb{I}}$ of $x$, it suffices to construct a point $y_{\mathbb{I}}\in W(\mathbb{I})$ lifting $y$. Let $\psi:A\to B$ be the dual isogeny of $\phi$. Since $\phi\circ\psi=[p]_A$ and $\psi\circ\phi=[p]_B$, we have
	\begin{equation}
		\textnormal{Im }\psi_{\ast,i}=\textnormal{Ker }\phi_{\ast,i},\qquad \textnormal{Ker }\psi_{\ast,i}=\textnormal{Im }\phi_{\ast,i}.
	\end{equation}
	Here, by slight abuse of notation, we use $\phi_{\ast,i},\psi_{\ast,i}$ to denote the induced morphisms on the crystalline homology groups. We take the lifting $\omega_{A^\vee_{\mathbb{I}}/\mathbb{I}}$ to be $\psi_{\ast,i}^{-1}(\mathcal{F}^i_{m_i-r_i})$ and $\mathcal{E}^i_{r_i,\mathbb{I}}$ to be $\textnormal{Im }\phi_{\ast,i}$. Applying Theorem \ref{Deformation of Abelian Schemes}, we obtain $(A_\mathbb{I},\lambda_\mathbb{I},\eta_{\mathbb{I}},\iota_{\mathbb{I}},\textnormal{Fil}'^{\bullet}_{\mathbb{I}})$ and $\phi_{\mathbb{I}}:B_{\mathbb{I}}\to A_{\mathbb{I}}$. Combining these data gives a point $y_{\mathbb{I}}$. This construction gives the inverse of $pr_2$.
	
	\begin{rmk}
		The map $pr_1$ does not induce an isomorphism on tangent spaces due to Frobenius twists issues. However, we know that $pr_2$ is an isomorphism because it induces a bijection on both closed points and tangent spaces. Therefore, $W$ is smooth, and by Lemma \ref{Frobenius twist morphisms}, we deduce that these two strata $Z$ and $Z'$ in partial flag spaces of different Shimura varieties are isomorphism up to Frobenius twists. This completes the proof. $\hfill\square$
	\end{rmk}

	\begin{rmk}
		There is a dual version of this proposition: Notations and assumptions as above, except that for each $1\le i\le N$, $\mathcal{E}^i_{s_i}$ is a sub-bundle of $\H(\mathcal{A}/X)$ of rank $s_i$ containing $\omega_i$. Then $Z$ is, up to some Frobenius twist, isomorphic to a stratum $Z''$ in a partial flag space $Y''_{P_0''}$ over $X'$, where $X''$ is the special fiber of a unitary Shimura variety with hyperspecial level structure and signature $(m_i'',n_i'')_{i=1}^N$, where
		\begin{equation}
			m_i''=m_i-s_i+s_{i-1},\quad n_i''=n-m_i'',\quad\textnormal{for $1\le i\le N$}.
		\end{equation}
		The proof is analogous.
	\end{rmk}
	
	\medskip
	This proposition has the limitation that if the signature at some $\tau_i$ is $(0,n)$, then the vector bundle $\mathcal{E}^i_{r_i}$ can only be $0$. Therefore, we strengthen the proposition by introducing essential Frobenius and Verschiebung maps.
	
	\begin{definition}
		We define the essential Frobenius map $F_{es}:\H(\mathcal{A}^{(p)}/X)_{i+1}\to \H(\mathcal{A}/X)_i$ and the essential Verschiebung map $V_{es}:\H(\mathcal{A}/X)_i\to \H(\mathcal{A}^{(p)}/X)_{i+1}$ to be
		\begin{equation}
			F_{es}=\left\{
			\begin{aligned}
				&F, \qquad\quad\textnormal{if $m_{i+1}\ne n$},\\
				&V^{-1}, \qquad\textnormal{if $m_{i+1}=n$}.
			\end{aligned}\right.
			\qquad\textnormal{and}\qquad
			V_{es}=\left\{
			\begin{aligned}
				&V, \qquad\quad\textnormal{if $m_{i+1}\ne0$},\\
				&F^{-1},\qquad\textnormal{if $m_{i+1}=0$}.
			\end{aligned}\right.
		\end{equation}
		We define $\tilde{\omega}_{i}$ as the image of $V_{es}$. More precisely, $\tilde{\omega}_{i}=\omega_i$ if $m_i\ne0$, and $\tilde{\omega}_{i}=\H(\mathcal{A}/X)_i$ otherwise. Let $\tilde{m}_i$ denote the rank of $\tilde{\omega}_i$ and $\tilde{n}_i=n-\tilde{m}_i$. We define an $F_{es},V_{es}$-chain in a similar manner as before.
	\end{definition}
	
	\begin{prop}\label{Construction via F_{es},V_{es}-chain}
		Let $B\subseteq P_0\subseteq P$ be an intermediate parabolic subgroup such that the associated partial flag space $Y_{P_0}$ is defined by recording at most one subbundle $\mathcal{E}^i_{r_i}\subseteq \tilde\omega_{i}$ of rank $r_i$ for each $i$ ($r_i$ can be $m_i$, in which case we do not record anything at the $\tau_i$-component). Let $Z$ be the closed subscheme of $Y_{P_0}$ defined by the condition
		\begin{equation}
			\mathcal{E}^i_{r_i}\subseteq V_{es}^{-1}(\mathcal{E}^{i+1,(p)}_{r_{i+1}})
		\end{equation}
		for all $i$. If moreover $\mathcal{E}^i_{r_{i}}=V^{-1}_{es}(\mathcal{E}^{i+1,(p)}_{r_{i+1}})$ whenever $m_i=0$. Then $Z$ is, up to some Frobenius factor, isomorphic to a stratum $Z'$ in some partial flag space $Y'_{P_0'}$ over $X'$, where $X'$ is the special fiber of a unitary Shimura variety with hyperspecial level structure and signature condition
		\begin{equation}
			m_i'=\left\{
			\begin{aligned}
				&0, \qquad\qquad\qquad\quad\textnormal{if $m_{i-1}=0$},\\
				&n-r_i+r_{i-1}, \qquad\textnormal{ if $m_i=0$ and $m_{i-1}\ne0$},\\
				&m_i-r_i+r_{i-1}, \quad\ \ \textnormal{ otherwise.}
			\end{aligned}
			\right.
		\end{equation}
	\end{prop}	
	\noindent\emph{Proof:} Pointwise, let $\tilde{\mathcal{E}}^{i}_{r_i}$ be the inverse images of $\mathcal{E}^i_{r_i}$ in $\tilde{D}(A)$ via $\tilde{D}(A)/p\tilde{D}(A)\simeq D(A)$. We define $\tilde{\mathcal{E}}'^i_{r_i}$ as
	\begin{equation}
		\tilde{\mathcal{E}}'^{i}_{r_i}=\left\{
		\begin{aligned}
			&\tilde{\mathcal{E}}^i_{r_i}, \qquad\ \textnormal{ if $m_i\ne0$},\\
			&p\tilde{\mathcal{E}}^i_{r_i}, \qquad\textnormal{if $m_i=0$}.
		\end{aligned}
		\right.
	\end{equation}
	However, the set $\{\tilde{\mathcal{E}}'^i_{r_i}\}$, is not an $F,V$-chain, and we cannot directly construct the universal $p$-isogeny $\phi:B\to A$ using the method in the proof of Proposition \ref{Construction via F,V-chain}, since $p\tilde{\mathcal{E}}^i_{r_i}\subseteq pD(A)_i$. However, we can consider the $F,V$-chain $\{\tilde{\mathcal{J}}_i\}$ given by
	\begin{equation}
		\tilde{J}_i=\left\{
		\begin{aligned}
			&\tilde{\mathcal{E}}^i_{r_i}, \qquad\quad\textnormal{ if $m_i\ne0$},\\
			&p\tilde{D}(A), \qquad\textnormal{otherwise.}
		\end{aligned}
		\right.
	\end{equation}
	Thus, $\tilde{J}_i$ forms an $F,V$-chain in $\tilde{D}(A)$, and $\tilde{\mathcal{E}}'^i_{r_i}$ forms an $F,V$-chain in $\oplus_{i=1}^N\tilde{J}_i$. Applying the construction in the proof of Proposition \ref{Construction via F,V-chain} twice, one can construct universal abelian schemes $\mathcal{B},\mathcal{C}$ and universal $p$-isogenies
	\begin{equation}
		\phi_1:\mathcal{B}\to \mathcal{A},\qquad \phi_2:\mathcal{C}\to \mathcal{B}
	\end{equation} 
	together with the data of auxiliary vector subbundles (plus inclusion relations). In this manner, we establish the desired isomorphism (up to Frobenius twists). $\hfill\square$
	
	This refined version will be frequently used in our later proof of the sufficiency part of the ampleness criterion.

	\medskip
	\subsection{ Review of some results on positivity in algebraic geometry.}	
	
	We will frequently use the following list of results in \cite{Lazarsfeld} on positivity in algebraic geometry.

    \begin{definition}\textup{\cite[Definition 1.4.1]{Lazarsfeld}} Let $X$ be a complete variety. A line bundle $\mathcal{L}$ on $X$ is called \emph{numerically effective} if
    \begin{equation}
        \big(\mathcal{L}\cdot C\big)\ge0
    \end{equation}
    for every irreducible curve $C\subseteq X$.
    \end{definition}
    
	\begin{prop}\textup{\cite[Proposition 1.2.13]{Lazarsfeld}}\label{Ampleness under finite morphisms}
		Let $f:X\to Y$ be a finite morphism of proper schemes and let $\mathcal{L}$ be a line bundle over $Y$. Then $f^\ast\mathcal{L}$ is ample if $\mathcal{L}$ is ample. Furthermore, if $f$ is surjective, then the converse holds.
	\end{prop}
	
	\begin{prop}\textup{\cite[Example 1.4.4]{Lazarsfeld}}\label{Nefness under proper morphisms}
		Let $f:X\to Y$ be a proper morphism of proper schemes, and $\mathcal{L}$ a line bundle over $Y$. Then $f^\ast\mathcal{L}$ is nef if $\mathcal{L}$ is nef. Furthermore, if $f$ is surjective, then the converse holds.
	\end{prop}	
	
	\begin{prop}\label{Nef+Ample=Ample}
		Let $X$ be a proper scheme and $D_1,D_2\in \textnormal{Pic $X$}_{\mathbb{Q}}$. If $D_1$ is ample and $D_2$ is nef, then $D_1+D_2$ is ample.
	\end{prop}	
	
	\begin{definition}\cite[Definition 1.7.1]{Lazarsfeld}
		Let $f:X\to T$ be a proper morphism of algebraic varieties and $\mathcal{L}$ be a line bundle on $X$. 
        
		(1) $\mathcal{L}$ is very ample relative to $f$, or $f$-very ample if the canonical map
		\begin{equation}
			\rho:f^\ast f_\ast\mathcal{L}\longrightarrow \mathcal{L}
		\end{equation}
		is surjective and defines an embedding $X\hookrightarrow \mathbb{P}(f_\ast\mathcal{L})$ of schemes over $T$.
		
		(2) $\mathcal{L}$ is ample relative to $f$, or $f$-ample, if $f^{\otimes m}$ is $f$-very ample for some $m>0$.
	\end{definition}

    \begin{prop}\cite[Theorem 1.7.8]{Lazarsfeld} Let $f$ be a morphism of proper noetherian schemes. Let $\mathcal{L}$ be a line bundle on $X$. For $t\in T$ set 
    \begin{equation}
        X_t=f^{-1}(t),\qquad \mathcal{L}_t=\mathcal{L}|_{X_t}.
    \end{equation}
    Then $\mathcal{L}$ is $f$-ample if and only if $\mathcal{L}_t$ is ample on $X_t$ for all $t\in T$ (That is, fiberwise ample). 
    \end{prop}
    
	\begin{prop}\textup{\cite[Proposition 1.7.10]{Lazarsfeld}}\label{Relative ample class gives ample classes}
		Let $f:X\to T$ be a morphism of projective schemes. Let $\mathcal{L}$ be a line bundle on $X$, and $\mathcal{E}$ be an ample line bundle on $T$. Then $\mathcal{L}$ is $f$-ample if and only if $\mathcal{L}\otimes \mathcal{E}^{\otimes m}$ is an ample line bundle on $X$ for all $m\gg0.$
	\end{prop}

	\section{Proof of the ampleness criterion: The necessity part}

	In this section, we assume $\mathcal{L}_Y(\lambda)$ is an ample line bundle on the flag space $Y$. We will show that the weight $\lambda$ must satisfy the inequalities in Theorem \ref{The ampleness criterion for the flag space, case 1}, Theorem \ref{The ampleness criterion for the flag space, case 2}, Theorem \ref{The ampleness criterion for the flag space, case 2'}, and Theorem \ref{The ampleness criterion for the flag space, case 3}, depending on the signature condition of the underlying Shimura variety. Recall that $T=\{i|\ m_i\ne 0,n\}$ denotes the \emph{essential set} of $X$, and its cardinality $t=\#T$ is the \emph{essential degree} of $X$.

	The fibers of the projection $Y\to X$ are flag varieties, so the ampleness implies that the weight is regular and $L$-dominant. The remaining inequalities come from a case-by-case analysis. The basic idea is to restrict $\mathcal{L}_Y(\lambda)$ to "special" subvarieties.
	
	Let $\mathcal{V}$ be a vector bundle over a base scheme $S$. Throughout this section, we write $[\mathcal{V}]$ for the class of $\textnormal{det }\mathcal{V}$ in $(\textnormal{Pic }S)_{\mathbb{Q}}$.
	
	\begin{rmk}
		By Proposition \ref{Ampleness under finite morphisms}, $\mathcal{L}_Y(\lambda)$ is ample if and only if its pullback under finite surjective morphisms is ample. So the Frobenius twists in Proposition \ref{Construction via F,V-chain} and \ref{Construction via F_{es},V_{es}-chain} may be ignored for our purpose.
	\end{rmk}

	\begin{lem}\label{Computing Picard Classes}
		For any vector bundle $\mathcal{F}\subseteq\omega_i$ (resp. $\mathcal{E}\supseteq \omega_i$) on $X$, we have
		\begin{equation}
			[V^{-1}(\mathcal{F}^{(p)})]=p[\mathcal{F}]-p[\omega_i] \qquad(\textnormal{resp. }F(\mathcal{E}^{(p)})=p[\mathcal{E}]-p[{\omega}_i])
		\end{equation}
		in $(\textnormal{Pic }X)_{\mathbb{Q}}$. Moreover, this holds for their pullbacks via any proper morphism $Z\to X$.
	\end{lem}
	\noindent\emph{Proof: }Since $[\H(\mathcal{A}/X)_{i-1}]=0$, the short exact sequence
	\begin{equation}
		0\longrightarrow \textnormal{Ker }V_{i-1}\longrightarrow \H(\mathcal{A}/X)_{i-1}\longrightarrow \omega_i^{(p)}\longrightarrow 0
	\end{equation}
	implies that $[\textnormal{Ker $V_{i-1}$}]=-p[\omega_i]$ in $(\textnormal{Pic $X$})_{\mathbb{Q}}$. The short exact sequence
	\begin{equation}
		0\longrightarrow \textnormal{Ker $V_{i-1}$}\longrightarrow V^{-1}(\mathcal{F}^{(p)})\longrightarrow\mathcal{F}^{(p)}\longrightarrow 0
	\end{equation}
	yields $[V^{-1}(\mathcal{F}^{(p)})]=[\mathcal{F}^{(p)}]+[\textnormal{Ker }V_{i-1}]\allowbreak=p[\mathcal{F}]-p[\omega_i]$. The formula for $\mathcal{E}$ can be deduced similarly. $\hfill\square$

	We first deal with the case when $X$ has essential degree 1.
	
	
	

	\begin{lem}\label{Necessity of the ampleness criterion, case 1}
		If $t=1$, we assume $T=\{i_1\}=\{1\}$ without loss of generality. If moreover $m_1=n-1$ and 
		$$
		[\mathcal{L}_{Y}(\lambda)]=[\mathcal{L}_Y(\{k^i_j\})]=\sum_{i=1}^N\sum_{j=1}^{n}k^i_j[\mathcal{F}^i_j/\mathcal{F}^i_{j-1}]
		$$
		is ample, then $\lambda$ is $L$-regular dominant and
		\begin{equation}
			\bigg(\sum_{u=0}^{n-2}p^{Nu}\bigg)\sum_{j=0}^{N-1}\frac{1}{p^j}(k^{1+j}_1-k^{1+j}_n)> \sum_{r=2}^{n-1}\bigg(p^{(r-2)N}\sum_{j=0}^{N-1}\frac{1}{p^j}(k^{1+j}_r-k^{1+j}_n)\bigg).
		\end{equation}
		Note that when $n=2$ the right hand side is 0.
	\end{lem}
	
	\noindent\emph{Proof:} By definition, the vector bundle $\omega_1=\omega_{\mathcal{A}^\vee/X,\tilde\tau_1}$ is locally free of rank $1$, and the vector bundle $\omega_i=\omega_{\mathcal{A}^\vee/X,\tilde\tau_i}$ has rank $0$ or $n$ for $2\le i\le N$. This signature condition implies that the map
	\begin{equation}
		V_{es}: \H(\mathcal{A}/X)_i\longrightarrow \H(\mathcal{A}/X)_{i+1}
	\end{equation}
	is an isomorphism for $1\le i\le N-1$, and has kernel locally free of rank $1$ for $i=N$.
	
	Let $Z\subseteq Y$ be the closed subscheme defined by the condition
	\begin{equation}
		\begin{aligned}
			&\mathcal{F}^i_j=V^{-1}_{es}(\mathcal{F}_{j}^{i+1,(p)}), \qquad \textnormal{for $1\le i\le N-1,1\le j<n$},\\
			&\mathcal{F}^N_j=V^{-1}_{es}(\mathcal{F}_{j-1}^{1,(p)}), \qquad\  \textnormal{for $2\le j< n$}.
		\end{aligned}
	\end{equation}
	This condition corresponds to (the closure of) a one dimensional stratum on the flag space. By Lemma \ref{Computing Picard Classes}, we have 
	\begin{equation}
		\begin{aligned}
			&[\mathcal{F}^i_j]=p[\mathcal{F}^{i+1}_j],\qquad \textnormal{for $1\le i\le N-1,1\le j<n$},\\
			&[\mathcal{F}^N_j]=p([\mathcal{F}^1_{j-1}]-[\omega_1]),\qquad \textnormal{for $2\le j<n$}.
		\end{aligned}
	\end{equation}
	in $(\textnormal{Pic $Z$})_{\mathbb{Q}}$. We can express each class in terms of $[\omega_1]$:
	\begin{equation}
		[\mathcal{F}^i_j]=\frac{1}{p^{i-1}}\big(1+\frac{1}{p^N}+\cdots+\frac{1}{p^{(n-j-1)N}}\big)[\omega_1],\qquad \textnormal{for $1\le i\le N,1\le j<n$}.
	\end{equation}
	
	The class $[\omega_1]$ is nef on $Z$, since it is the pullback of the class of the Hodge line bundle $\textnormal{det }\omega_{\mathcal{A}^\vee/X}$ on $X$. We claim that $[\omega_1]$ is, moreover, ample. Notice that the image of $Z$ under the projection map $Y\to X$ is the closure of the unique 1-dimensional Ekedahl--Oort stratum (One could use Example \ref{Computation of Ekedahl--Oort stratum} to write down the conditions), and this projection map onto its image is actually a Frobenius factor. Thus, the pullback of the Hodge bundle to $Z$ is ample. There is an alternative way to show this: We construct a strata Hasse invariant by
    \begin{equation}
        \mathcal{F}^N_1\longrightarrow V^{-1}_{es}(\mathcal{F}^{1,(p)}_1)/\textnormal{Ker }V\stackrel{V}{\longrightarrow} \mathcal{F}^{1,(p)}_1.
    \end{equation}
    Here, the first map is induced by the inclusion map. This gives a section $h\in H^0(Z,\mathcal{F}^{1,(p)}_1\otimes \mathcal{F}^{N,-1}_1)$ whose zero locus is the zero dimensional Ekedahl--Oort stratum on $Y$. By the above computation of classes, we have
    \begin{equation}
        [h]=p(1-\frac{1}{p^{N(n-1)}})[\omega_1].
    \end{equation}
    This shows that $[\omega_1]$ has positive degree on $Z$, hence ample.
    
	Plugging the above equations into $[\mathcal{L}_Y(\lambda)]|_Z$, we obtain
	\begin{equation}
		\begin{aligned}
			[\mathcal{L}_Y(\{k^i_j\})]\big|_Z &=\sum_{i=1}^N\sum_{j=1}^n k^i_j[\mathcal{F}^i_j/\mathcal{F}^i_{j-1}]\big|_Z =\sum_{i=1}^{N}\sum_{j=1}^{n-1}k^i_j[\mathcal{F}^i_j/\mathcal{F}^i_{j-1}]\big|_Z-\sum_{i=1}^{N}k^i_n[\mathcal{F}^i_{n-1}]\\
			&=\sum_{i=1}^{N}\sum_{j=1}^{n-1}(k^i_j-k^i_n)[\mathcal{F}^i_j/\mathcal{F}^i_{j-1}]\big|_Z\\
			&=\bigg(\sum_{i=1}^N\frac{1}{p^{i-1}}\sum_{j=2}^{n-1} (-\frac{1}{p^{(n-j)N}})(k^i_j-k^i_n)+\sum_{i=1}^{N}\frac{1}{p^{i-1}}\big(1+\frac{1}{p^N}+\cdots+\frac{1}{p^{(n-2)N}}\big)(k^i_1-k^i_n)\bigg)[\omega_1]\big|_Z.
		\end{aligned}
	\end{equation}
	The positivity of the coefficient gives the desired inequality. $\hfill\square$

	\begin{lem}\label{Necessity of the ampleness criterion, case 3}
		If $t=\# T\ge 2$ and $[\mathcal{L}_Y(\lambda)]$ is ample, then $\lambda$ is $L$-regular dominant and
		\begin{equation}
			p^{a(i_l)}\sum_{j=0}^{a(i_l)-1}\frac{1}{p^{j}}(k_1^{i_l+j}-k_n^{i_l+j})>k^{i_{l+1}}_{m_{i_{l+1}}}-k^{i_{l+1}}_{m_{i_{l+1}}+1},\qquad \textnormal{for all $1\le l\le t$}.
		\end{equation}
		Here, the subscripts of $i_l$ are viewed as modulo $t$.
	\end{lem}
	
	\noindent\emph{Proof:} Let $Z$ be the stratum over which the Hodge filtration agrees with the conjugate filtration at $\tau_{i_l},\tau_{i_l+1},\dots,\tau_{i_{l+1}-1}$-components, i.e., $Z$ is the closed subscheme defined by
	\begin{equation}
		\begin{aligned}
			&\mathcal{F}^{i_l+j}_r=V_{es}^{-1}(\mathcal{F}^{i_l+j+1,(p)}_r),\qquad \textnormal{for $0\le j<a(i_l)-1,1\le r<n$},\\
			&\mathcal{F}^{i_{l+1}-1}_r=F_{es}(\mathcal{F}_{m_{i_{l+1}}+r}^{i_{l+1},(p)}),\qquad \ \textnormal{for $1\le r\le n_{i_{l+1}}$},\\
			&\mathcal{F}^{i_{l+1}-1}_r=V^{-1}_{es}(\mathcal{F}_{r-n_{i_{l+1}}}^{i_{l+1},(p)}),\qquad \textnormal{for $n_{i_{l+1}}<r<n$}.
		\end{aligned}
	\end{equation}
	In fact, we may write the last two conditions uniformly by viewing the subscripts as modulo $n$. By Lemma \ref{Computing Picard Classes}, we have
	\begin{equation}
		\begin{aligned}
			&[\mathcal{F}^{i_l+j}_r]=p[\mathcal{F}^{i_l+j+1}_r],\qquad \textnormal{for $0\le j<a(i_l)-1,1\le r<n$},\\
			&[\mathcal{F}^{i_{l+1}-1}_r]=p[\mathcal{F}^{i_{l+1}}_{m_{i_{l+1}}+r}/\omega_{i_{l+1}}],\qquad \textnormal{for $1\le r<n$},
		\end{aligned}
	\end{equation}
	in $(\textnormal{Pic $Z$})_{\mathbb{Q}}$. Assume $m_{i_l}+m_{i_{l+1}}\ge n$ for simplicity, the above condition implies the tuple
	\begin{equation}
		(\H(\mathcal{A}/X)_1,\dots,\H(\mathcal{A}/X)_{i_l-1},\omega_{i_l},\mathcal{F}^{i_l+1}_{m_{i_l}},\dots,\mathcal{F}^{i_{l+1}-1}_{m_{i_l}},\H(\mathcal{A}/X)_{i_{l+1}},\dots,\H(\mathcal{A}/X)_N)
	\end{equation}
	forms an $F_{es},V_{es}$-chain. By Proposition \ref{Construction via F_{es},V_{es}-chain}, we obtain a correspondence
	
	\begin{center}
		\begin{tikzpicture}[node distance=10pt and 1.5cm, auto]
			\node(1) [draw, align=center] {$W$};
			\node(2) [draw, align=center, below left=of 1] {$Z$};
			\node(3) [draw, align=center, below right=of 1] {$Z'$};
			\node(4) [draw, align=center, left=of 2] {$Y$};
			\node(5) [draw, align=center, right=of 3] {$Z''$};
			\node(6) [draw, align=center, right=of 5] {$X'$};
			
			\draw [->] (1) -- node[midway, above] {$pr_1$} (2);
			\draw [->] (1) -- node[midway, above] {$pr_2$} (3);
			\draw [<-] (4) -- node[midway, above] {$\supseteq$} (2);
			\draw [->] (3) -- node[midway, above] {$\pi$} (5);
			\draw [->] (5) -- node[midway, above] {$\pi$} (6);
		\end{tikzpicture}
	\end{center}
	To make the proof more readable, we recall the moduli problems of all the auxiliary varieties in this diagram.
	
	Here, $X'$ is a unitary Shimura variety with signature $(m'_{i},n'_i)_{i=1}^N$. We note that $m'_{i_{l+1}}=m_{i_l}+m_{i_{l+1}}-n$. Since $Z$ involves the data of a refinement of Hodge filtration, we must record the data for the auxiliary varieties correspondingly. More precisely, $Z'$ is the moduli space whose $S$ points are tuples $(B,\lambda_B,\eta_B,\iota_B,\textnormal{Fil}_B)$ where
	
	$\bullet$ $(B,\lambda_B,\eta_B,\iota_B)\in X'(S)$.
	
	$\bullet$ $\textnormal{Fil}_B$ is a refinement of the Hodge filtration, which agrees with $\textnormal{Fil}$ away from the $i_l,\dots,i_{l+1}$-places, and
	\begin{equation}
		\begin{aligned}
			\begin{aligned}
				&\textnormal{Fil}_B^{i_{l}+j}\textnormal{ is the Hodge filtration for $0\le j< a(i_l)$, i.e., we do not take any refinements},\\
				&\textnormal{Fil}_{B}^{i_{l+1}}: 0\subseteq \mathcal{E}^{i_{l+1}}_1\subseteq\cdots\subseteq\mathcal{E}^{i_{l+1}}_{m'_{i_{l+1}}-1}\subseteq\mathcal{E}^{i_{l+1}}_{m'_{i_{l+1}}}=\omega_{B,i_{l+1}}\subseteq\mathcal{E}^{i_{l+1}}_{m'_{i_{l+1}}+1}\subseteq\cdots\subseteq\mathcal{E}^{i_{l+1}}_{n-1}\subseteq\H(B/S)_{i_{l+1}}.		
			\end{aligned}
		\end{aligned}
	\end{equation}

	We take $W$ to be the moduli space whose $S$-points consist of tuples $(A,\lambda,\eta,\iota,\textnormal{Fil},B,\lambda_B,\eta_B,\iota_B,\textnormal{Fil}_B,\varphi)$, where
	
	$\bullet$ $(A,\lambda,\eta,\iota,\textnormal{Fil})\in Z(S)$.
	
	$\bullet$ $(B,\lambda_B,\eta_B,\iota_B,\textnormal{Fil}_B)\in Z'(S)$.
	
	$\bullet$ $\varphi:B\to A$ is an isogeny whose kernel is killed by $p$, such that the induced morphism $\varphi_\ast:\allowbreak\H({B}/S)\allowbreak\to \H({A}/S)$ is an isomorphism away from the $\tau_{i_l},\tau_{i_l+1},\dots,\tau_{i_{l+1}}$-components, and 
	\begin{equation}
		\varphi_{i,\ast}(\mathcal{E}^i_j)=\mathcal{F}^i_j,\qquad\textnormal{for $i\notin\{i_l,\dots,i_{l+1}\}$ and $1\le j<n$}.
	\end{equation}
	At $\tau_{i_l},\dots,\tau_{i_{l+1}-1}$, we require
	\begin{equation}
		\begin{aligned}
			&\textnormal{Im }\varphi_{i_{l+1}-j,\ast} =\mathcal{F}^{i_{l+1}-j}_{m_{i_l}},\qquad \textnormal{for $1\le j\le a(i_l)$}.
		\end{aligned}
	\end{equation}
	At $\tau_{i_{l+1}}$, we require $\varphi_{i_{l+1},\ast}$ to be an isomorphism, and
	\begin{equation}
		\varphi_{i_{l+1},\ast}(\mathcal{E}^{i_{l+1}}_{r})=\mathcal{F}^{i_{l+1}}_{r},\qquad \textnormal{for $1\le r<n$}.
	\end{equation}
	
	We define $Z''$ as the partial flag space over $X'$ which only records refinement of the Hodge filtration at places away from the $\tau_{i_l},\dots,\tau_{i_{l+1}}$-components. The forgetful map $Z'\to Z''$ is proper surjective, whose fibers are isomorphic to flag varieties
	\begin{equation}
		Fl_{GL_{2n-m_{i_l}-m_{i_{l+1}}}}\times Fl_{GL_{m_{i_l}+m_{i_{l+1}}-n}}.
	\end{equation}
	
	It follows from Proposition \ref{Description of strata} (with additional filtration data) that $pr_1$ and $pr_2$ are isomorphisms up to Frobenius factors. If $\mathcal{L}_Y(\lambda)$ is ample, it has to be ample on $Z$, hence $Z'$, hence the fibers of $Z'\to Z''$. 
	Pick any closed point $z''\in Z''$, let $C_{z''}$ denote the fiber of $Z'\to Z''$ at $z''$. We have
	\begin{equation}
		\begin{aligned}		
			[\mathcal{L}_Y(\lambda)]\big|_{C_{z''}}&=\sum_{i=1}^N\sum_{u=1}^{n}k^i_u[\mathcal{F}^i_u/\mathcal{F}^i_{u-1}]|_{C_{z''}}\\
			&=\sum_{j=0}^{a(i_l)}\sum_{u=1}^{n}k^{i_l+j}_{u}[\mathcal{F}^{i_l+j}_{u}/\mathcal{F}^{i_l+j}_{u-1}]\big|_{C_{z''}}\\
			&=\bigg(\sum_{j=0}^{a(i_l)-1}\sum_{u=1}^{n}k^{i_l+j}_{u}[\mathcal{F}^{i_{l}+j}_{u}/\mathcal{F}^{i_l+j}_{u-1}] +\sum_{v=1}^{n}k^{i_{l+1}}_{v}[\mathcal{F}^{i_{l+1}}_{v}/\mathcal{F}^{i_{l+1}}_{v-1}]\bigg)\big|_{C_{z''}}\\
			&=\bigg(\sum_{u=1}^{n}\big(\sum_{j=0}^{a(i_l)-1}p^{a(i_l)-j-1}k_u^{i_l+j}\big) [\mathcal{F}^{i_{l+1}-1}_u/\mathcal{F}^{i_{l+1}-1}_{u-1}] +\sum_{v=1}^{n}k^{i_{l+1}}_{v}[\mathcal{F}^{i_{l+1}}_{v}/\mathcal{F}^{i_{l+1}}_{v-1}]\bigg)\big|_{C_{z''}}\\
			&=\bigg(\sum_{u=1}^{n_{i_{l+1}}}\big(\sum_{j=0}^{a(i_l)-1}p^{a(i_l)-j}k_u^{i_l+j}\big) [\mathcal{F}^{i_{l+1}}_{m_{i_{l+1}}+u}/\mathcal{F}^{i_{l+1}}_{m_{i_{l+1}}+u-1}]\\
			&+\sum_{u=n_{i_{l+1}}+1}^{n}\big(\sum_{j=0}^{a(i_l)-1}p^{a(i_l)-j}k_u^{i_l+j}\big) [\mathcal{F}^{i_{l+1}}_{u-n_{i_{l+1}}}/\mathcal{F}^{i_{l+1}}_{u-n_{i_{l+1}}-1}] +\sum_{v=1}^{n}k^{i_{l+1}}_{v}[\mathcal{F}^{i_{l+1}}_{v}/\mathcal{F}^{i_{l+1}}_{v-1}]\big)\bigg)\big|_{C_{z''}}\\
			&=\bigg(\sum_{u=1}^{n_{i_{l+1}}}\big(\sum_{j=0}^{a(i_l)-1}p^{a(i_l)-j}k_u^{i_l+j}\big) [\mathcal{E}^{i_{l+1}}_{m_{i_{l+1}}+u}/\mathcal{E}^{i_{l+1}}_{m_{i_{l+1}}+u-1}]\\
			&+\sum_{u=n_{i_{l+1}}+1}^{n}\big(\sum_{j=0}^{a(i_l)-1}p^{a(i_l)-j}k_u^{i_l+j}\big) [\mathcal{E}^{i_{l+1}}_{u-n_{i_{l+1}}}/\mathcal{E}^{i_{l+1}}_{u-n_{i_{l+1}}-1}] +\sum_{v=1}^{n}k^{i_{l+1}}_{v}[\mathcal{E}^{i_{l+1}}_{v}/\mathcal{E}^{i_{l+1}}_{v-1}]\big)\bigg)\big|_{C_{z''}}.
		\end{aligned}
	\end{equation}
	Over each fiber $C_{z''}$, the weight is required to be regular dominant. This in particular implies that "the coefficient of $[\mathcal{E}^{i_{l+1}}_{m_{i_{l+1}}+1}/\mathcal{E}^{i_{l+1}}_{m_{i_{l+1}}}]$" is bigger than "the coefficient of $[\mathcal{E}^{i_{l+1}}_{m_{i_{l+1}}}/\mathcal{E}^{i_{l+1}}_{m_{i_{l+1}}-1}]$". That is,
	\begin{equation}
		\big(\sum_{j=0}^{a(i_l)-1}p^{a(i_l)-j}k^{i_l+j}_1\big)+k^{i_{l+1}}_{m_{i_{l+1}}+1} >\big(\sum_{j=0}^{a(i_l)-1}p^{a(i_l)-j}k^{i_l+j}_n\big)+k^{i_{l+1}}_{m_{i_{l+1}}}.
	\end{equation}
	In other words,
	\begin{equation}
		p^{a(i_l)}\bigg(\sum_{j=0}^{a(i_l)-1}\frac{1}{p^j}\big(k^{i_l+j}_1-k^{i_l+j}_n\big)\bigg) >k^{i_{l+1}}_{m_{i_{l+1}}}-k^{i_{l+1}}_{m_{i_{l+1}}+1}.
	\end{equation}
	We remark that the other inequalities derived from the regular dominant condition are redundant. 
	
	If we are in the case $m_{i_l}+m_{i_{l+1}}<n$, we can consider the $F_{es},V_{es}$-chain
	\begin{equation}
		(0,\dots,0,\omega_{i_l},\mathcal{F}^{i_l+1}_{m_{i_l}},\dots, \mathcal{F}^{i_{l+1}-1}_{m_{i_l}},0,\dots,0),
	\end{equation}
	and apply Proposition \ref{Construction via F_{es},V_{es}-chain} to show that $Z$ is, up to Frobenius factors, isomorphic to a stratum $Z'$ in the partial flag space of an auxiliary unitary Shimura variety with different signature. One can construct a projection $Z'\to Z''$ whose fibers are isomorphic to
	\begin{equation}
		Fl_{GL_{n-m_{i_l}-m_{i_{l+1}}}}\times Fl_{GL_{m_{i_l}+m_{i_{l+1}}}}
	\end{equation}
	in a completely same manner. Restricting to geometric fibers of this morphism, the regular dominant condition again implies "the coefficient of $[\mathcal{E}^{i_{l+1}}_{m_{i_{l+1}}+1}/\mathcal{E}^{i_{l+1}}_{m_{i_{l+1}}}]$" is bigger than "the coefficient of $[\mathcal{E}^{i_{l+1}}_{m_{i_{l+1}}}/\mathcal{E}^{i_{l+1}}_{m_{i_{l+1}}-1}]$". A direct but boring computation of Picard group relations yields the same inequality. An alternative approach is to work on the $\tilde\tau_i^c$-side and reduce to the previous case.
	
	Letting $l$ run over $1,2,\dots,t$, we obtain the desired necessary condition. $\hfill\square$

	\begin{lem}\label{Necessity of the ampleness criterion, case 2} If $t=1$, we assume $T=\{i_1\}=\{1\}$ for simplicity. If moreover $m_1\notin \{1,n-1\}$ and 
		\begin{equation}
			[\mathcal{L}_Y(\lambda)]=[\mathcal{L}_Y(\{k^i_j\})]=\sum_{i=1}^N\sum_{j=1}^nk^i_j[\mathcal{F}^i_j/\mathcal{F}^i_{j-1}]
		\end{equation}
		is ample, then $\lambda$ is $L$-regular dominant and
		\begin{equation}
			p^N\big(\sum_{j=0}^{N-1}\frac{1}{p^j}(k^{1+j}_1-k^{1+j}_n)\big)>k^1_{m_1}-k^1_{m_1+1}.
		\end{equation}
	\end{lem}
	\begin{rmk}
		The inequalities in Lemma \ref{Necessity of the ampleness criterion, case 3} and Lemma \ref{Necessity of the ampleness criterion, case 2} could be written uniformly. However, since there is only one essential place in this case, the choice of the $F_{es},V_{es}$-chain in the proof of Lemma \ref{Necessity of the ampleness criterion, case 3} is not valid in this case. It is somewhat surprising that the criterion arises from a relatively unusual construction.
	\end{rmk}
	
	\noindent\emph{Proof:} Without loss of generality, we may assume $m_2=n$, otherwise, we may consider the filtrations on $\tau^c$-side. This is a minor technical assumption to ensure that the $F_{es},V_{es}$-chain in the following is a genuine $F,V$-chain. Let $Z$ denote the stratum defined by the condition
	\begin{equation}
		\begin{aligned}
			&\mathcal{F}^i_j=V^{-1}_{es}(\mathcal{F}^{i+1,(p)}_{j}), \qquad&&\textnormal{for $2\le i\le N-1,1\le j\le n$},\\
			&\mathcal{F}^1_j=V^{-1}_{es}(\mathcal{F}^{2,(p)}_{j}), \qquad&&\textnormal{for $j\ne m_1$},\\
			&\mathcal{F}^N_j=V^{-1}_{es}(\mathcal{F}^{1,(p)}_{j-n_1}), \qquad&&\textnormal{for $j\notin \{n_1,m_1\}$}.
		\end{aligned}
	\end{equation}
	Here, the subscripts are again treated modulo $n$. The stratum $Z$ is 2-dimensional if $m_1= n_1$, and is 3-dimensional if $m_1\ne n_1$. Note that over $Z$, the tuple
	\begin{equation}
		(\mathcal{F}^1_{m_1-1},\mathcal{F}^2_{m_1},\dots,\mathcal{F}^N_{m_1})
	\end{equation}
	forms an $F_{es},V_{es}$-chain. It is indeed an $F,V$-chain by assumption. Applying Proposition \ref{Construction via F,V-chain}, we obtain a correspondence
	
	
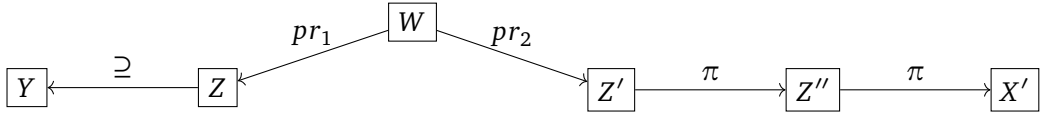
\begin{center}
		\begin{tikzpicture}[node distance=10pt and 1.5cm, auto]
			\node(1) [draw, align=center] {$W$};
			\node(2) [draw, align=center, below left=of 1] {$Z$};
			\node(3) [draw, align=center, below right=of 1] {$Z'$};
			\node(4) [draw, align=center, left=of 2] {$Y$};
			\node(5) [draw, align=center, right=of 3] {$Z''$};
			\node(6) [draw, align=center, right=of 5] {$X'$};
			
			\draw [->] (1) -- node[midway, above] {$pr_1$} (2);
			\draw [->] (1) -- node[midway, above] {$pr_2$} (3);
			\draw [<-] (4) -- node[midway, above] {$\supseteq$} (2);
			\draw [->] (3) -- node[midway, above] {$\pi$} (5);
			\draw [->] (5) -- node[midway, above] {$\pi$} (6);
		\end{tikzpicture}
	\end{center}
	For simplicity, we use the same letters as in the previous lemma, but they stand for entirely different objects. We now recall their moduli problems:
	
	$X'$ is the auxiliary Shimura variety of appropriate level structure, as given by Proposition \ref{Description of strata}. The signature of $X'$ changes only at the first two places, with $m_1'=m_1+1$ and $m_2=n-1$.
	
	$Y'$ is the full flag space over $X'$, which means that $Y'(S)$ consists of isomorphism classes of tuples $(B,\lambda_B,\eta_B,\iota_B,\textnormal{Fil}_B)$, where
	
	$\bullet$ $(B,\lambda_B,\eta_B,\iota_B)\in X'(S)$;
	
	$\bullet$ $\textnormal{Fil}_B$ is a refinement of the Hodge filtration:
	\begin{equation}
		\begin{aligned}
			&\textnormal{Fil}_B^{1}:0\subseteq \mathcal{E}^1_1\subseteq\cdots\subseteq \mathcal{E}^1_{m_1}\subseteq\mathcal{E}^1_{m_1+1}=\omega_{B,1}\subseteq\mathcal{E}^1_{m_1+2}\subseteq\cdots\subseteq\mathcal{E}^1_{n-1}\subseteq\H(B/S)_1, \\
			&\textnormal{Fil}_B^{2}:0\subseteq \mathcal{E}^2_1\subseteq\cdots\subseteq \mathcal{E}^2_{n-2}\subseteq\mathcal{E}^2_{n-1}=\omega_{B,2}\subseteq \H(B/S)_2, \\
			&\textnormal{Fil}_B^{1+j}:0\subseteq \mathcal{E}^{1+j}_1\subseteq\cdots\subseteq\mathcal{E}^{1+j}_{n-1}\subseteq\H(B/S)_{1+j}, \qquad\textnormal{for $2\le j\le N-1$}.
		\end{aligned}
	\end{equation}
	
	$Z'\subseteq Y'$ is the stratum defined by the condition (In the following, we always assume $m_1\ge n_1$ for simplicity. Otherwise, one changes the third equality to be "for $u\ge n_1$" and the fourth equality to be "for $u\le n_1-1$ and $u\ne n_1-m_1$".)
	\begin{equation}
		\begin{aligned}
			&\mathcal{E}^1_{u}=V^{-1}(\mathcal{E}^{2,(p)}_{u-1}), \quad\ \ \ \qquad\textnormal{for $u\notin\{1,m_1+1\}$},\\
			&\mathcal{E}^{1+j}_u=V^{-1}_{es}(\mathcal{E}^{2+j,(p)}_{u}), \qquad\textnormal{for $1\le j\le N-2,1\le u\le n-1$},\\
			&\mathcal{E}^N_{u}=V^{-1}(\mathcal{E}^{1,(p)}_{u-n_1+1}), \ \ \qquad\textnormal{for $u\ge n_1$ and $u\ne 2n_1$},\\
			&\mathcal{E}^N_{u}=F(\mathcal{E}^{1,(p)}_{u+m_1+1}), \ \quad \qquad\textnormal{for $u\le n_1-1$}.
		\end{aligned}
	\end{equation}
	We remark that if $m_1=n_1$, then the third condition above becomes "for $u\ge n_1$". This again justifies the fact that $Z$ is 3-dimensional if $m_1\ne n_1$, and 2-dimensional if $m_1=n_1$.
	
	$W$ is the moduli space whose $S$ points consists of tuples $(A,\lambda,\eta,\iota,B,\lambda_B,\eta_B,\iota_B,\varphi)$, where
	
	$\bullet$ $(A,\lambda,\eta,\iota)\in Z(S)$;
	
	$\bullet$ $(B,\lambda_B,\eta_B,\iota_B)\in Z'(S)$;
	
	$\bullet$ $\varphi:B\to A$ is an isogeny, whose kernel is killed by $p$, and the induced morphism $\varphi_{i,\ast}:\allowbreak\H(B/S)_i\allowbreak\to \H(A/S)_i$ satisfies (let $\psi$ be the dual isogeny of $\varphi$)
	\begin{equation}
		\begin{aligned}
			\textnormal{At $\tau_{1}$:}\qquad\qquad\qquad\qquad\quad &\textnormal{Ker }\varphi_{1,\ast}=\mathcal{E}^1_{n_1+1}, \\ 
			&\varphi_{1,\ast}(\mathcal{E}^1_{u})=\mathcal{F}^1_{u-n_1-1}, \ \qquad\textnormal{for $n_1+2\le u\le n$ and $u\ne m_1+1$},\\
			&\psi_{1,\ast}(\mathcal{F}^1_v)=\mathcal{E}^1_{v-m_1+1},\qquad\textnormal{for $m_1\le v\le n-1$}.\\
			\textnormal{At $\tau_{j}$ for $2\le j\le N$:}\qquad\quad\ &\textnormal{Ker }\varphi_{j,\ast}=\mathcal{E}^{j}_{n_1}, \\
			&\varphi_{j,\ast}(\mathcal{E}^{j}_{u})=\mathcal{F}^{j}_{u-n_1}, \qquad\  \textnormal{for $n_1+1\le u\le n$},\\
			&\psi_{j,\ast}(\mathcal{F}^{j}_{v})=\mathcal{E}^{j}_{v-m_1}, \qquad\textnormal{for $m_1+1\le v\le n-1$}.
		\end{aligned}
	\end{equation}
	We remark that if $m_1<n_1$, we change the condition at $\tau_1$ to be
		\begin{equation}
		\begin{aligned}
			\textnormal{At $\tau_{1}$:}\qquad\qquad\qquad\qquad\quad &\textnormal{Ker }\varphi_{1,\ast}=\mathcal{E}^1_{n_1+1}, \\ 
			&\varphi_{1,\ast}(\mathcal{E}^1_{u})=\mathcal{F}^1_{u-n_1-1}, \ \qquad\textnormal{for $n_1+2\le u\le n$},\\
			&\psi_{1,\ast}(\mathcal{F}^1_v)=\mathcal{E}^1_{v-m_1+1},\qquad\textnormal{for $m_1\le v\le n-1$ and $v\ne 2m_1$}.
		\end{aligned}
	\end{equation}

	$Y''$ is a partial flag space over $X'$ which records all the $\mathcal{E}^i_j$'s except that
	\begin{equation}
		\mathcal{E}^1_1,\qquad \mathcal{E}^1_{n_1+1},\qquad\textnormal{and }\mathcal{E}^{j}_{n_1}, \textnormal{ for $2\le j\le N$}.
	\end{equation}
	$Z''$ is the image of $Z'$ under the natural forgetful map $Y'\to Y''$. In fact, $Z$ is the stratum cut out by the same conditions as $Z'$ except for those involving $\mathcal{E}^1_{1}$, $\mathcal{E}^1_{n_1+1}$, and $\mathcal{E}^{j}_{n_1}$ for $2\le j\le N$. Giving a point on the fiber of $Z'\to Z''$ is equivalent to choosing $0\subseteq\mathcal{E}^1_1\subseteq\mathcal{E}^1_2$ and putting
	\begin{equation}
		\mathcal{E}^{N}_{n_1}=V^{-1}(\mathcal{E}^{1,(p)}_1),\ \mathcal{E}^{N-1}_{n_1}=V^{-1}_{es}(\mathcal{E}^{N,(p)}_{n_1}),\ \dots,\mathcal{E}^2_{n_1}=V^{-1}_{es}(\mathcal{E}^{3,(p)}_{n_1}),\ \mathcal{E}^1_{n_1+1}=V^{-1}(\mathcal{E}^{2,(p)}_{n_1}).
	\end{equation}
	Therefore, the fibers are $\mathbb{P}^1$, and the above vector bundles are the only nontrivial classes on the fibers.
	
	If $\mathcal{L}_Y(\lambda)$ is ample, then it is ample on $Z'$, and consequently, ample over the fibers of $Z'\to Z''$. Let $z\in Z''(k)$ be a closed point, then in $\textnormal{Pic $\mathbb{P}^1_z$}$, we have
	\begin{equation}
		\begin{aligned}
			&[\mathcal{F}^1_{m_1+1}/\mathcal{F}^1_{m_1}]=[\mathcal{E}^1_2/\mathcal{E}^1_1]=[\mathcal{O}(1)],\\
			&[\mathcal{F}^2_{m_1+1}/\mathcal{F}^2_{m_1}]=p[\mathcal{F}^3_{m_1+1}/\mathcal{F}^3_{m_1}]=\cdots=p^{N-2}[\mathcal{F}^N_{m_1+1}/\mathcal{F}^N_{m_1}]=p^{N-2}[\mathcal{E}^N_1]=0,\\
			&[\mathcal{F}^1_{m_1}/\mathcal{F}^1_{m_1-1}]=[\mathcal{E}^1_1]=[\mathcal{O}(-1)],\\
			&[\mathcal{F}^2_{m_1}/\mathcal{F}^2_{m_1-1}]=p[\mathcal{F}^3_{m_1}/\mathcal{F}^3_{m_1-1}]=\cdots=p^{N-2}[\mathcal{F}^N_{m_1}/\mathcal{F}^N_{m_1-1}]=[\H(B/X')_2/\omega_{B,2}]=0,\\
			&[\H(A/X)_1/\mathcal{F}^1_{n-1}]=p[\H(A/X)_2/\mathcal{F}^2_{n-1}]=p^2[\H(A/X)_3/\mathcal{F}^3_{n-1}]=\cdots=p^{N-1}[\H(A/X)_N/\mathcal{F}^N_{n-1}]\\
			&[\H(A/X)_1/\mathcal{F}^1_{n-1}]=[\mathcal{E}^1_{n_1+1}/\mathcal{E}^1_{n_1}]=\cdots=p^{N-1}[\mathcal{E}^N_{n_1}/\mathcal{E}^N_{n_1-1}]
			=p^N[\mathcal{E}^1_1]=p^N[\mathcal{O}(-1)],\\
			&[\mathcal{F}^1_1]=p[\mathcal{F}^2_1]=\cdots=p^{N-1}[\mathcal{F}^N_1],\\
			&[\mathcal{F}^N_1]=[\mathcal{E}^N_{n_1+1}/\mathcal{E}^N_{n_1}]=p[\mathcal{E}^1_2/\mathcal{E}^1_1]=p[\mathcal{O}(1)].
		\end{aligned}	
	\end{equation}
	For $[\mathcal{F}^i_j/\mathcal{F}^i_{j-1}]$, where $(i,j)$ are not in the above list, it is trivial on the fiber. Hence, we have
	\begin{equation}
		\begin{aligned}
			[\mathcal{L}_Y(\lambda)]\big|_{\mathbb{P}^1_z}&=\sum_{i=1}^N\sum_{j=1}^nk^i_j[\mathcal{F}^i_j/\mathcal{F}^i_{j-1}]\big|_{\mathbb{P}^1_z}=\sum_{i=1}^N\sum_{j\in\{1,m_1,m_1+1,n\}}k^i_j[\mathcal{F}^i_j/\mathcal{F}^i_{j-1}]\big|_{\mathbb{P}^1_z}\\
			&=\big(k^1_{m_1+1}-k^1_{m_1}-\sum_{i=1}^{N}p^{N+1-i}k^i_n+\sum_{i=1}^{N}p^{N+1-i}k^i_1\big)[\mathcal{O}(1)].
		\end{aligned}
	\end{equation}
	This gives the desired inequality
	\begin{equation}
		p^N\big(\sum_{j=0}^{N-1}\frac{1}{p^j}(k^{1+j}_1-k^{1+j}_n)\big)>k^1_{m_1}-k^1_{m_1+1}.
	\end{equation}
	This completes the proof. $\hfill\square$

	\medskip
	Combining these lemmas, we have established the necessity part of the ampleness criterion by a case-by-case argument.
	
	\begin{rmk}
		The proofs of the lemmas in this section remain valid if we replace "ample" by "nef", "regular dominant" by "dominant", and $>$ by $\ge$. Therefore, we also obtained the necessity part of the \emph{nefness} criterion.
	\end{rmk}
	
	\begin{rmk}
		A line bundle over a projective scheme is nef if and only if its pullback along a proper surjective morphism is nef. Applying this to the projection $Y=Y_B\to Y_{P'}$, we have also established the necessity part of the nefness criterion for $\mathcal{L}_{P'}(\lambda)$. Here, $B\subseteq P'\subseteq P$ is an intermediate parabolic subgroup, and $\lambda$ is a character of $P'$.
	\end{rmk}

	\section{Slopes}\label{Slope Chapter}
	
	The goal of this section is to introduce the notion that we call "slope" to measure different Ekedahl--Oort strata (on Shimura variety or any partial flag space). The proof of the sufficiency part of the ampleness (or nefness) criterion in the next chapter is based on induction with respect to slopes.
	
	\medskip
	\subsection{ Setup and definitions.}
	In the following, we continue to denote by $X$ the special fiber of the unitary Shimura variety with good reduction, where we have the vector bundles $\omega_i$. Let $\pi:T\to X$ be any proper morphism. As before, we modify the Frobenius and Verschiebung morphisms to be $F_{es}$ and $V_{es}$, respectively. Let $\tilde\omega_i$ be the image of $V_{es}$, i.e.,
	\begin{equation}
		\tilde\omega_i=\left\{
		\begin{aligned}
			&\omega_i,\qquad\qquad\quad\ \  \textnormal{if $m_i\ne0$},\\
			&\H(\mathcal{A}/X)_i,\qquad \textnormal{if $m_i=0$}.
		\end{aligned}	\right.
	\end{equation}
	Let $\tilde{m}_i$ denote the rank of $\tilde{\omega}_i$ and $\tilde{n}_i=n-\tilde{m}_i$. By abuse of notation, for any vector bundle $\mathcal{V}$ on $X$, we denote the pullback $\pi^\ast \mathcal{V}$ to $T$ still by $\mathcal{V}$ as long as the underlying scheme we are considering is obvious in the context.
	
	For $\mathcal{M}_1$ and $\mathcal{M}_2$ two subbundles of a vector bundle $\mathcal{M}$, the intersection sheaf $\mathcal{M}_1\cap\mathcal{M}_2$ is in general not a coherent sheaf, but the function
	\begin{equation}
		\varphi:s\mapsto \textnormal{dim}_{k(s)}\ \mathcal{M}_{1,s}\cap \mathcal{M}_{2,s}
	\end{equation} 
	takes a constant value $r$ on some dense open subset $U$ of $X$. 
	We will say that the intersection $\mathcal{M}_1\cap \mathcal{M}_2$ \emph{generically has rank $r$}.
	
	\begin{definition}[The slope of a vector bundle]\label{Slopes of vector bundles}
		Let $\mathcal{F}$ be any subbundle of rank $r$ contained in $\omega_i$. We define \emph{the slope of $\mathcal{F}$ on $T$ of order $s$} as the sequence of integers
		\begin{equation}
			\vec{r}_s(T,\mathcal{F}):=(r_1,r_2,\dots,r_{sN}).
		\end{equation}
		Here, the integers $r_a$ are defined inductively as follows.
		
		$\bullet$ We take $\mathcal{F}_0:=\mathcal{F}$, $\tilde{r}_0=r$, and $T_0:=T$.
		
		$\bullet$ For $1\le j\le sN-1$, assume that $\mathcal{F}_{j-1}\subseteq \tilde\omega_{i-j+1}$ and $T_{j-1}\subseteq T$ are defined, where $\mathcal{F}_{j-1}$ is a subbundle of $\tilde\omega_{i-j+1}$ of rank $\tilde{r}_{j-1}$ on $T_{j-1}$. Here, the subscript for $\tilde\omega$ is again viewed as modulo $N$. Then it is reasonable to talk about $V_{es}^{-1}(\mathcal{F}_{j-1}^{(p)})$, and we define
		$$
		\tilde{r}_{j}:=\textnormal{The generic intersection rank of $\tilde\omega_{i-j}\cap V_{es}^{-1}(\mathcal{F}_{j-1}^{(p)})$ on $T_{j-1}$}.
		$$
		We define
		\begin{equation}
			\begin{aligned}
				&T_{j}=\textnormal{A dense open subset of $T_{j-1}$ over which the above intersection has generic rank}.\\
				&\mathcal{F}_j:=\tilde\omega_{i-j}\cap V^{-1}_{es}(\mathcal{F}_{j-1}^{(p)}) \textnormal{ on $T_j$},\textnormal{ which is a subbundle of $\tilde\omega_{i-j}$},\\
				&r_j:=\tilde{n}_{i-j+1}+\tilde{r}_{j-1}-\tilde{r}_j.
			\end{aligned}
		\end{equation}
		
		$\bullet$ For $j=sN$, we define
        \begin{equation}
            \tilde{r}_{sN}:=\textnormal{The generic intersection rank of }\mathcal{F}_{(s-1)N}\cap  V^{-1}_{es}(\mathcal{F}^{(p)}_{sN-1}) \textnormal{ on } T_{sN-1}.
        \end{equation}
        We define
        \begin{equation}
            \begin{aligned}
                &T_{sN}=\textnormal{A dense open subset of } T_{sN-1} \textnormal{ over which the above intersection has generic rank}.\\
				&\mathcal{F}_{sN}:=\mathcal{F}_{(s-1)N}\cap V^{-1}_{es}(\mathcal{F}_{sN-1}^{(p)}) \textnormal{ on }T_{sN},\textnormal{ which is a subbundle of }\mathcal{F}_{(s-1)N},\\
				&r_{sN}:=\tilde{n}_{i+1}+\tilde{r}_{sN-1}-\tilde{r}_{sN}.
            \end{aligned}
        \end{equation}

		If $\mathcal{E}\supseteq \omega_i$, we define $\vec{r}_s(T,\mathcal{E})$ as $\vec{r}_s(T,\mathcal{E}^\perp)$ by working on the conjugate side.
	\end{definition}
	
	\begin{rmk}
		The choice of the open subsets $T_j$'s does not affect our definition. Hence, we can always take $T_j$ to be the largest locus over which the intersection attains its generic rank.
	\end{rmk}
	
	\begin{definition}[Total slope]\label{Total slope}
		Assumptions as above, we define \emph{the total slope of $\mathcal{F}$ on $T$ of order $s$} to be
		\begin{equation}
			r_{s,\textnormal{tot}}(T,\mathcal{F}):=\sum_{j=1}^{sN}r_j.
		\end{equation}
	\end{definition}
	
	\begin{definition}[Order of slopes]\label{Total order of slopes}
		For vector bundles $\mathcal{F},\mathcal{E}\subseteq \omega_i$, we define an order between $\vec{r}_s(T,\mathcal{F})$ and $\vec{r}_s(T,\mathcal{E})$ by \emph{lexicographical order}. Note that the slopes are comparable only when the vector bundles are both subbundles of the same $\omega_i$. We define a total order on the set of total slopes using the natural order.
	\end{definition}
	
	\begin{rmk}
		In this paper, we only need to use slopes of order 1 for the sufficiency part of the ampleness criterion. Hence, unless otherwise stated, we always assume $s=1$ and drop the subscript $s$ from our notions.
	\end{rmk}

	\medskip
	\subsection{ Illustration of slopes.} The definition of slopes may seem indirect. However, we introduce a method to appropriately represent the combinatorical data in a diagram, in which our notion of slopes is the slopes of the connecting lines between the vertices.
	
	\medskip
	\noindent\emph{Step 1: Draw the nodes.} We consider the Cartesian coordinate system.
	
	$\bullet$ Put $0$ at $(1,0)$, and $\mathcal{F}^1_j$ at $(1,j)$ for all $1\le j\le n$. We finish the construction of the first column.
	
	$\bullet$ Next, suppose that the first $u$ rows have already been drawn, and $\H(A/X)_u$ is placed at $(u,v)$, then we place $\tilde\omega_{u+1}$ at $(u+1,v)$. 
	
	$\bullet$ For all $l$, place $\mathcal{F}^{u+1}_l$ at $(u+1,v-m_{u+1}+l)$. We then finish the $(u+1)$-th column.
	
	$\bullet$ Conclude the inductive process at the $(N+1)$-th column.
	
	$\bullet$ Erase the points that correspond to the vector bundles $\mathcal{F}^{u}_l$ that contain $\tilde\omega_u$ for all $1\le u\le N+1$.
	
	$\bullet$ Draw segments connecting the top (resp. bottom) elements of adjacent columns.
	
	We remark that in the definition of slopes, we work on any scheme $T$ over $X$. We do not require $T$ to record the refinement of the Hodge filtration, i.e., $T$ is not necessarily a scheme over the flag space $Y$. At this moment, one should think of the $\mathcal{F}^i_j$'s that we put in the diagram simply as some symbols, not genuine vector bundles on $T$.

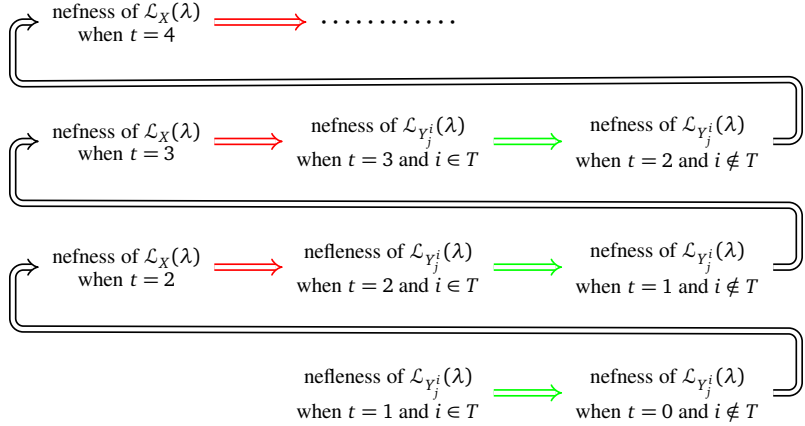
	\begin{example}	Let $G=G(U(2,1)\times U(3,0)\times U(1,2)\times U(0,3))$, the diagram becomes (Note that $\tilde\omega_i$ differs from $\omega_i$ only when $m_i=0$.)
		\begin{center}
			\begin{tikzpicture}[node distance=10pt and 1cm, auto]
				\begin{scope}
					\node (1) [draw, rectangle, minimum width=1cm, minimum height=0.75cm, align=center] {$\omega_1$};
					\node (2) [draw, rectangle, minimum width=1cm, minimum height=0.75cm, align=center, below=of 1] {$\mathcal{F}^1_1$};
					\node (3) [draw, rectangle, minimum width=1cm, minimum height=0.75cm, align=center, below=of 2] {$0$};
					\node (4) [draw, rectangle, minimum width=1cm, minimum height=0.75cm, align=center, above right=of 1] {$\omega_2$};
					\node (5) [draw, rectangle, minimum width=1cm, minimum height=0.75cm, align=center, below=of 4] {$\mathcal{F}^2_2$};
					\node (6) [draw, rectangle, minimum width=1cm, minimum height=0.75cm, align=center, below=of 5] {$\mathcal{F}^2_1$};
					\node (7) [draw, rectangle, minimum width=1cm, minimum height=0.75cm, align=center, below=of 6] {0};
					\node (8) [draw, rectangle, minimum width=1cm, minimum height=0.75cm, align=center, right=of 4] {$\omega_3$};
					\node (9) [draw, rectangle, minimum width=1cm, minimum height=0.75cm, align=center, below=of 8] {0};
					\node (10) [draw, rectangle, minimum width=1cm, minimum height=0.75cm, align=center, right=of 9] {0};
					\node (11) [draw, rectangle, minimum width=1cm, minimum height=0.75cm, align=center, above=of 10] {$\mathcal{F}^4_1$};
					\node (12) [draw, rectangle, minimum width=1cm, minimum height=0.75cm, align=center, above=of 11] {$\mathcal{F}^4_2$};
					\node (13) [draw, rectangle, minimum width=1cm, minimum height=0.75cm, align=center, above=of 12] {$\mathcal{F}^4_3$};
					\node (14) [draw, rectangle, minimum width=1cm, minimum height=0.75cm, align=center, above=of 13] {$\tilde\omega_4$};
					\node (15) [draw, rectangle, minimum width=1cm, minimum height=0.75cm, align=center, right=of 14] {$\omega_1$};
					\node (16) [draw, rectangle, minimum width=1cm, minimum height=0.75cm, align=center, below=of 15] {$\mathcal{F}^1_1$};
					\node (17) [draw, rectangle, minimum width=1cm, minimum height=0.75cm, align=center, below=of 16] {$0$};
				\end{scope}
				
				\draw[-, thick] (1)--(4);
				\draw[-, thick] (4)--(8);
				\draw[-, thick] (8)--(14);
				\draw[-, thick] (14)--(15);
				
				\draw[-, thick] (3)--(7);
				\draw[-, thick] (7)--(9);
				\draw[-, thick] (9)--(10);
				\draw[-, thick] (10)--(17);
			\end{tikzpicture}
		\end{center}
	\end{example}
	
	\medskip
	\noindent\emph{Step 2: Draw the lines.} To simplify the notations, we assume that $\mathcal{F}$ is a subbundle of $\omega_1$ of rank $l$ and denote it by $\mathcal{F}^1_l$. By the definition of slope, there is a dense open subset $T_N$ of $T$ on which all intersections of vector bundles that we considered have generic ranks. We perform the following operations to the above diagram:
	\begin{itemize}[itemsep=0pt,topsep=2pt, parsep=0pt]
	    \item For $0\le j\le N$, select all the vertices $\mathcal{F}^{1-j}_{\tilde{r}_j}$ in the diagram that correspond to the $\mathcal{F}_j$'s in our definition of slopes. Then we disregard all the other $\mathcal{F}^u_v$'s.
        \item For $0\le j\le N-1$, draw line segments connecting $\mathcal{F}^{1-j}_{\tilde{r}_j}$ and $\mathcal{F}^{-j}_{\tilde{r}_{j-1}}$. (The superscripts of $\mathcal{F}$ are viewed as modulo $N$)
	\end{itemize}

	
	Our notion of the slope of $\mathcal{F}=\mathcal{F}^1_l$ corresponds to the sequence of slopes of line segments from the right to the left of the diagram.

	\begin{example}
		Let $X$ be a 
        unitary Shimura variety associated with $G=G(U(4,1)\times U(2,3)\times U(3,2))$. Consider the partial flag space $Y^1_2$ over $X$, which records a subbundle $\mathcal{F}^1_2$ of rank 1 of $\omega_1$. The following diagram indicates that the slope of $\mathcal{F}^1_2$ on $Y^1_2$ is $\vec{r}(Y^1_2,\mathcal{F}^1_2)=(2,3,3)$.
		
		\begin{center}
			\begin{tikzpicture}[node distance=10pt and 1cm, auto]
				
				\begin{scope}
					\node (1) [draw, rectangle, minimum width=1cm, minimum height=0.75cm, align=center] {$\omega_1$};
					\node (2) [draw, rectangle, minimum width=1cm, minimum height=0.75cm, align=center, below=of 1] {$\cancel{\mathcal{F}^1_3}$};
					\node (3) [draw, rectangle, minimum width=1cm, minimum height=0.75cm, align=center, below=of 2] {$\mathcal{F}^1_2$};
					\node (4) [draw, rectangle, minimum width=1cm, minimum height=0.75cm, align=center, below=of 3] {$\cancel{\mathcal{F}^1_1}$};
					\node (5) [draw, rectangle, minimum width=1cm, minimum height=0.75cm, align=center, below=of 4] {$0$};
					\node (6) [draw, rectangle, minimum width=1cm, minimum height=0.75cm, align=center, above right=of 1] {$\omega_2$};
					\node (7) [draw, rectangle, minimum width=1cm, minimum height=0.75cm, align=center, below=of 6] {$\cancel{\mathcal{F}^2_1}$};
					\node (8) [draw, rectangle, minimum width=1cm, minimum height=0.75cm, align=center, below=of 7] {$0$};
					\node (9) [draw, rectangle, minimum width=1cm, minimum height=0.75cm, align=center, right=of 6] {$0$};
					\node (10) [draw, rectangle, minimum width=1cm, minimum height=0.75cm, align=center, above=of 9] {$\mathcal{F}^3_1$};
					\node (11) [draw, rectangle, minimum width=1cm, minimum height=0.75cm, align=center, above=of 10] {$\cancel{\mathcal{F}^3_2}$};
					\node (12) [draw, rectangle, minimum width=1cm, minimum height=0.75cm, align=center, above=of 11] {$\omega_3$};
					\node (13) [draw, rectangle, minimum width=1cm, minimum height=0.75cm, align=center, above right=of 12] {$\cancel{\mathcal{F}^1_3}$};
					\node (14) [draw, rectangle, minimum width=1cm, minimum height=0.75cm, align=center, above=of 13] {$\omega_1$};
					\node (15) [draw, rectangle, minimum width=1cm, minimum height=0.75cm, align=center, below=of 13] {$\mathcal{F}^1_2$};
					\node (16) [draw, rectangle, minimum width=1cm, minimum height=0.75cm, align=center, below=of 15] {$\cancel{\mathcal{F}^1_1}$};
					\node (17) [draw, rectangle, minimum width=1cm, minimum height=0.75cm, align=center, below=of 16] {$0$};

					\draw[-, line width=1.5pt] (1.north) -- (6.west);
					\draw[-, line width=1.5pt] (6.north) -- (12.west);
					\draw[-, line width=1.5pt] (12.north) -- (14.west);
					\draw[-, line width=1.5pt] (5.east) -- (8.south);
					\draw[-, line width=1.5pt] (8.east) -- (9.south);
					\draw[-, line width=1.5pt] (9.east) -- (17.south);
					
					\draw[-, thick, red] (15.west)--(10.north); 
					\draw[-, thick, red] (10.west)--(8.north); 
					\draw[-, thick, red] (8.west)--(5.north); 

					
				\end{scope}
			\end{tikzpicture}
		\end{center}
	\end{example}
	In other words, the intersection $V^{-1}(\mathcal{F}^{1,(p)}_2)\cap\omega_3$ has generic rank 1 on $T_N$, and we denote this intersection by $\mathcal{F}^3_1$ (One should pay attention that, $\mathcal{F}^3_1$ only represents a vector bundle defined on $T_N$. This differs from the one in the definition of partial flag space $Y^3_1$, although we are using the same symbols), and the intersection $V^{-1}(\mathcal{F}^{3,(p)}_1)\cap \omega_2$ is trivial on $T_N$, and finally the intersection $V^{-1}(0)\cap \mathcal{F}^1_2$ is trivial on $T_N$.
	
	\medskip
	The following properties of slopes follow directly from the definition.
	\begin{lem}\label{Properties of slope}
		Let $\mathcal{F}$ be a subbundle of $\omega_i$ over a scheme $T$, with $\vec{r}(T,\mathcal{F})=(r_1,\dots,r_N)$, then
		
		\textnormal{(1)} $r_j\ge0$ for all $j$.

        \textnormal{(2)} $r_j\le \tilde{n}_{i-j}$ for $1\le j\le N-1$.
        
		\textnormal{(3)} $r_{\textnormal{tot}}(T,\mathcal{F})\ge \sum_{j=1}^N\tilde{n}_j$, and the equality occurs if and only if the set of vector bundles $(\mathcal{F}_j)_{j=1}^N$ forms an $F_{es},V_{es}$-chain.
	\end{lem}
	\noindent\emph{Proof: }Trivial.$\hfill\square$

\section{Proof of the ampleness criterion: The sufficiency part}
	
\subsection{ Overview of the proof}

The goal of this section is to prove the sufficiency part of our main theorems: Theorem \ref{The ampleness criterion for the flag space, case 1}, Theorem \ref{The ampleness criterion for the flag space, case 2}, Theorem \ref{The ampleness criterion for the flag space, case 2'} and Theorem \ref{The ampleness criterion for the flag space, case 3}. The proof involves many inductive arguments as well as combinatorical choices. Before going into the proof, it is better to give an overview of the strategy of our proof in this subsection.

\begin{convention}
    In our proof, we will frequently pull line bundles back via proper (surjective) morphisms $\pi:Z_1\to Z_2$. For simplicity, we will only write $\pi^\ast\mathcal{L}$, $\mathcal{L}|_{Z_1}$, or even just $\mathcal{L}$ without specifying $\pi$ for these pullbacks, as long as the meaning is obvious in the paragraph.
\end{convention}
	
\begin{convention}
    To simplify our notation, we drop the subscripts of $V_{es}$ and $F_{es}$ and simply write $V$ and $F$. Similarly, we use $\omega_i$ to denote the image of \emph{essential} Verschiebung and $m_i$ to be the rank of $\omega_i$. (Originally we used $\tilde{\omega}_i$ and $\tilde{m}_i$ for these terms)
\end{convention}

	In Section 7.2, we reduce the sufficiency part of the ampleness criterion to the corresponding nefness criterion of the \emph{minimal partial flag spaces}. More explicitly, to prove the sufficiency part of Theorem \ref{Main Theorem - Introduction}, it suffices to prove the sufficiency part of the following theorems.
	
	\begin{thm}\label{Nefness criterion for minimal partial flag space}
		Let $P'=P^i_j$ be a minimal proper parabolic subgroup, and $Y_{P'}=Y^i_j$ be the corresponding partial flag space obtained by recording a rank $j$ vector bundle contained in (or containing, depending on $j$) $\omega_i$ in the moduli problem. Let 
		$$
			\mathcal{L}_{Y^i_j}(\lambda)=\mathcal{L}_{Y^i_j}(\{k_l\};\alpha) =\bigotimes_{\substack{l=1\\ l\ne i}}^N (\textnormal{det }\omega_l)^{\otimes k_l}\otimes \big(\textnormal{det }(\tilde\omega_i/\mathcal{F}^i_j)\big)^{\otimes k_i}\otimes (\textnormal{det }\mathcal{F}^i_j)^{\otimes \alpha}
		$$ 
		be the automorphic line bundle of weight $\lambda$ over $Y^i_j$.
		
		\begin{enumerate}
			\item  If $T$ is the empty set, then $\mathcal{L}_{Y^i_j}(\{k_l\};\alpha)$ is nef if and only if
			\begin{equation}
				k_i\ge\alpha.
			\end{equation}
			
			From now on we assume that $T$ is non-empty, then there is a unique $1\le t_0\le t$, such that $i\in \{i_{t_0},i_{t_0}+1,\dots, i_{t_0+1}-1\}$. To simplify the notation, we assume that $i_1=i_{t_0}\le i$.
			
			\item If $T$ is a singleton, namely, $T=\{i_1\}$. If moreover $m_{i_0}=n-1$, then $\mathcal{L}_{Y^i_j}(\{k_l\};\alpha)$ is nef if and only if
			\begin{equation}
				\begin{aligned}
					&k_i\ge\alpha,\\
					&(p^{(n-j)N}-1)\big(k_{i_1}-\frac{1}{p^{i-i_1}}(k_i-\alpha)\big)\ge (p^{(n-j-1)N}-1)k_{i_1}.
				\end{aligned}
			\end{equation}
			
			\item If $T$ is a singleton, namely, $T=\{i_1\}$. If moreover $m_{i_0}=1$, then $\mathcal{L}_{Y^i_j}(\{k_l\};\allowbreak\alpha)$ is nef if and only if
			\begin{equation}
				\begin{aligned}
					&k_i\ge \alpha,\\
					&(p^{jN}-1)\big(k_{i_1}-\frac{1}{p^{i-i_1}}(k_i-\alpha)\big)\ge (p^{(j-1)N}-1)k_{i_1}.
				\end{aligned}
			\end{equation}
			
			\item In the rest of the cases, we have $t\ge2$, or $t=1$ but $m_{i_1}\notin\{1,n-1\}$. 
            Then $\mathcal{L}_{Y^i_j}(\{k_l\};\allowbreak\alpha)$ is nef if and only if
			\begin{equation}
				\begin{aligned}
					&k_i\ge \alpha,\\
					&p^{a(i_{1})}\big(k_{i_1}-\frac{1}{p^{i-i_{1}}}(k_i-\alpha)\big)\ge k_{i_{1}+1},\\
					&p^{a(i_s)}k_{i_s}\ge k_{i_{s+1}}, \qquad\textnormal{for }2\le s\le t.
				\end{aligned}
			\end{equation}
		\end{enumerate}
	\end{thm}

    We will also prove the nefness criterion for certain line bundles over the Shimura variety $X$.
	\begin{thm}\label{Nefness criterion for Shimura variety}
		Write $T=\{i_1,\dots,i_t\}$. The line bundle $$\mathcal{L}_X(\lambda)=\mathcal{L}_X(\{k_l\})=\bigotimes_{l=1}^N (\textnormal{det }\omega_l)^{\otimes k_l}$$ over $X$ is nef if and only if
		\begin{equation}
			p^{a(i_s)}k_{i_s}\ge k_{i_{s+1}}
		\end{equation}
		holds for $1\le s\le t$.
	\end{thm}

    Notice that we simply change all the strict inequalities "$>$" to "$\ge$" for the nef cones. This agrees with the fact that the nef cone is the closure of the ample cone.
    
	\vspace{1em}
	We prove Theorem \ref{Nefness criterion for minimal partial flag space} and Theorem \ref{Nefness criterion for Shimura variety} simultaneously. The strategy, roughly speaking, is by induction with respect to the essential degree $t=\# T$. The following diagram illustrates how our inductive process works:

	\[
	\begin{tikzcd}
		\substack{\textnormal{nefness of $\mathcal{L}_X(\lambda)$}\\ \textnormal{when $t=4$}}  \arrow[r, Rightarrow, red] 
		&\cdots\cdots\cdots\cdots
		&\ \\
		\substack{\textnormal{nefness of $\mathcal{L}_X(\lambda)$}\\ \textnormal{when $t=3$}}  \arrow[r, Rightarrow, red] 
		&\substack{\textnormal{nefness of $\mathcal{L}_{Y^i_j}(\lambda)$}\\ \textnormal{when $t=3$ and $i\in T$}}  \arrow[r, Rightarrow, green] 
		&\substack{\textnormal{nefness of $\mathcal{L}_{Y^i_j}(\lambda)$}\\ \textnormal{when $t=2$ and $i\notin T$}}\arrow[ull, Rightarrow, rounded corners, to path={ -- ([xshift=2ex]\tikztostart.east)  
			-- ([xshift=2ex, yshift=4.75ex]\tikztostart.east) 
			-- ([xshift=-2ex, yshift=-4.75ex]\tikztotarget.west) 
			-- ([xshift=-2ex]\tikztotarget.west) 
			-- (\tikztotarget) }]\\
		\substack{\textnormal{nefness of $\mathcal{L}_X(\lambda)$}\\ \textnormal{when $t=2$}}  \arrow[r, Rightarrow, red] 
		&\substack{\textnormal{nefleness of $\mathcal{L}_{Y^i_j}(\lambda)$}\\ \textnormal{when $t=2$ and $i\in T$}}  \arrow[r, Rightarrow, green] 
		&\substack{\textnormal{nefness of $\mathcal{L}_{Y^i_j}(\lambda)$}\\ \textnormal{when $t=1$ and $i\notin T$}}
		\arrow[ull, Rightarrow, rounded corners, to path={ -- ([xshift=2ex]\tikztostart.east)  
			-- ([xshift=2ex, yshift=4.75ex]\tikztostart.east) 
			-- ([xshift=-2ex, yshift=-4.75ex]\tikztotarget.west) 
			-- ([xshift=-2ex]\tikztotarget.west) 
			-- (\tikztotarget) }]\\
		\textnormal{ }
		&\substack{\textnormal{nefleness of $\mathcal{L}_{Y^i_j}(\lambda)$}\\ \textnormal{when $t=1$ and $i\in T$}}  \arrow[r, Rightarrow, green] 
		&\substack{\textnormal{nefness of $\mathcal{L}_{Y^i_j}(\lambda)$}\\ \textnormal{when $t=0$ and $i\notin T$}} 
		\arrow[ull, Rightarrow, rounded corners, to path={ -- ([xshift=2ex]\tikztostart.east)  
			-- ([xshift=2ex, yshift=4.75ex]\tikztostart.east) 
			-- ([xshift=-2ex, yshift=-4.75ex]\tikztotarget.west) 
			-- ([xshift=-2ex]\tikztotarget.west) 
			-- (\tikztotarget) }] 
	\end{tikzcd}
	\]
	
	\medskip
	Here, each block means that "we prove the nefness of the line bundle whose weight satisfies the given the corresponding inequality in Theorem \ref{Nefness criterion for minimal partial flag space} and Theorem \ref{Nefness criterion for Shimura variety}, under the inductive hypothesis that the sufficiency part of the nefness criterion corresponding to the previous blocks holds". The weight $\lambda$ has different meanings for each block. 
    
    The case $t=0$ is the nefness criterion for line bundles over flag varieties (in characteristic $p$); The nefness criterion for $\mathcal{L}_X(\lambda)$ when $t=1$ is the ampleness of Hodge line bundle. We do not include these two trivial cases in the diagram. The reader should also pay attention to the essential degree $t$ on the rightmost column.

    \vspace{1em}
	In Sections 7.4, 7.5 and 7.6, we prove the nefness criterion corresponding to the left, middle and right columns respectively. This is the core of our proof.

\noindent\textbf{Strategy of the proof:} The basic idea of the proof is to analyze the intersection number
\begin{equation}
    \big(\mathcal{L}\cdot C\big)
\end{equation}
for $\mathcal{L}=\mathcal{L}_{Y^i_m}(\lambda)$ or $\mathcal{L}=\mathcal{L}_X(\lambda)$, where $\lambda$ satisfies the numerical conditions for the nef cone, and $C$ is any complete curve over $Y$ or $X$. We try to write the class $[\mathcal{L}]$ as a non-negative linear combination of certain divisors and nef classes.
\begin{equation}
    [\mathcal{L}]=\sum [h_i]+[\mathcal{L}_{\textnormal{nef}}].
\end{equation}
Here, each $h_i$ is a "generalized strata partial Hasse invariant", whose zero locus is an Eckdahl--Oort strata of lower dimension. $\mathcal{L}_{\textnormal{nef}}$ is (the pullback of) an automorphic line bundle of an appropriate choice of weight, which is known to be nef by induction hypothesis. Such an expression immediately implies
\begin{equation}
    \big(\mathcal{L}\cdot C\big)\ge0,
\end{equation}
as long as $C$ is not contained in the zero locus of any strata Hasse invariant $h_i$. If $C$ is contained in some $Z(h_i)$, then we repeat this process by a new choice of strata Hasse invariants and the nef line bundle $\mathcal{L}_{\textnormal{nef}}$. This helps us to reduce the dimension and eventually deduce the non-negativity of the intersection numbers.

\begin{rmk}
    As is explained above, in order to establish a desired decomposition of (the class of) a given automorphic line bundle, one has to overcome two technical difficulties: Construction of partial Hasse invariants over strata and resolution of singularities of certain Ekedahl--Oort strata. Both problems are crucial to understand the geometry and cohomology of the special fiber of Shimura varieties. There are only partial results for both questions. We refer to the introduction part of this paper for a survey of some relevant research. In our proof, we make use the moduli interpretation to construct strata Hasse invariants, and also resolve certain singular closed Ekedahl--Oort strata by working on a partial flag space given by more refined data.
\end{rmk}

	To prove the statements for the left column, it suffices to show that the weights on the boundary of the conjectured cone are nef weights. For example, we shall prove that $\lambda_0=(1,p^{a(i_1)},p^{a(i_1)+a(i_2)},\dots, p^{a_{i_1}+\cdots+a(i_{t-1})})$ is a nef weight. To attack this problem, we 
    study the so-called \emph{Schubert cells} (see Subsection 7.4 for the explicit definition), which are closures of certain Ekedahl--Oort strata. These closed subschemes are defined by explicit conditions in terms of the rank of intersections of $\omega_{i_j}$ and $\textnormal{Ker }V^{a(i_j)}$. Through these moduli description, it is easy to construct strata Hasse invariants which cut out lower dimensional Schubert cells. This enables us to establish the desired equalities as above.


	
	
	
	To prove the statement for the middle and right columns, we make induction on "slopes of $\mathcal{F}^i_j$". The rough picture for the proof is
	
	\[
	\begin{tikzcd}
		\substack{\textnormal{nefness of $\mathcal{L}_{Y^i_j}(\lambda)|_{Z_0}$}\\
			\textnormal{where $r_{\textnormal{tot}}(Z_0,\mathcal{F}^i_j)=\sum_{j=1}^Nn_j$}
		} \arrow[r, Rightarrow]
		&\substack{\textnormal{nefness of $\mathcal{L}_{Y^i_j}(\lambda)|_{Z_1}$}\\
			\textnormal{where $\vec{r}(Z_1,\mathcal{F}^i_j)>\vec{r}(Z_0,\mathcal{F}^i_j)$}
		}\arrow[r, Rightarrow]
		&\cdots\arrow[r, Rightarrow]
		&\substack{\textnormal{nefness of $\mathcal{L}_{Y^i_j}(\lambda)|_{Z}$}\\
			\textnormal{where $\vec{r}(Z,\mathcal{F}^i_j)$ is generic}
		}
	\end{tikzcd}
	\]
	In other words, we construct strata Hasse invariants, compute intersection numbers $(\mathcal{L}_{Y^i_j}(\lambda)\cdot C)$, and reduce the problem to lower dimensional subschemes over which the slope of $\mathcal{F}^i_j$ becomes smaller. A crucial remark is that, by each block, we mean proving the nefness of $\mathcal{L}_{Y^i_j}(\lambda)$ for all $j$'s simultaneously. For example, to prove $\mathcal{L}_{Y^1_m}(\lambda)$ is nef over $Z_1$, we have to use the inductive hypothesis that $\mathcal{L}_{Y^1_{m'}}(\lambda')|_{Z_0}$ is nef, provided that both $\lambda$ and $\lambda'$ belong to their respective nef cones.
	
	In subsection 7.3, we prove the induction base for all our future inductions.

    Now, let us state the induction hypothesis for this whole chapter.
	
\begin{definition}
    Let $X$ be the geometric special fiber of a unitary Shimura variety, and let $Y^i_j$ be a minimal partial flag space over $X$. Denote by $T$ the essential set of $X$ and $t=\# T$ the essential degree of $X$. For any proper morphism $\pi:W\to Y^i_j$, we define the tuple
    \begin{equation}
        \mathfrak{I}(W,\mathcal{L}_{Y^i_j}(\lambda)):=\left\{
        \begin{aligned}
            &\big(t,\vec{r}(W,\pi^\ast\mathcal{F}^i_j),j\big), \qquad && \textnormal{if }i\in T\\
            &\big(t+1,\vec{r}(W,\pi^\ast\mathcal{F}^i_j),j\big),\qquad  &&\textnormal{if }i\notin T.
        \end{aligned}\right.
    \end{equation}
    Similarly, for any proper morphism $W\to X$, we define
    \begin{equation}
        \mathfrak{I}(W,\mathcal{L}_X(\lambda)):=\big(t,(0,\dots,0),0\big).
    \end{equation}
    Finally, let $\mathfrak{I}$ be the set of all such tuples, and we define the order on $\mathfrak{I}$ as the lexicographical order.
\end{definition}
	
In the following subsections, we prove the ampleness criterion for \emph{all} unitary Shimura varieties simultaneously by induction with respect to the partial order on $\mathfrak{I}$. More explicitly,

\medskip
\noindent\textbf{Induction hypothesis:} To prove that $\mathcal{L}_{Y^i_j}(\lambda)|_W$ (or $\mathcal{L}_X(\mu)$) is nef if the weight $\lambda$ (or $\mu$) satisfies the corresponding inequalities in Theorem \ref{Nefness criterion for minimal partial flag space} and Theorem \ref{Nefness criterion for Shimura variety}, we assume that $\mathcal{L}_{{Y'}^{i'}_{j'}}(\lambda')|_Z$ and $\mathcal{L}_{X'}(\lambda'')|_Z$ is nef as long as the weight $\lambda'$ and $\lambda''$ satisfies the inequalities in the same theorems, and
\begin{equation}
\begin{aligned}
    &\mathfrak{I}(W,\mathcal{L}_{Y^i_j}(\lambda)) \textnormal{ or }\mathfrak{I}(W,\mathcal{L}_X(\mu))> \mathfrak{I}(Z,\mathcal{L}_{Y'^{i'}_{j'}}(\lambda')),\\
    &\mathfrak{I}(W,\mathcal{L}_{Y^i_j}(\lambda)) \textnormal{ or }\mathfrak{I}(W,\mathcal{L}_X(\mu))> \mathfrak{I}(Z,\mathcal{L}_{X'}(\lambda'')).
\end{aligned}
\end{equation}
Here, $X'$ runs over the special fiber of unitary Shimura varieties of all possible signatures and $Y'^{i'}_{j'}$ is a minimal partial flag space over $X'$. Please note that we prove the sufficiency part of the ampleness criterion for all unitary Shimura varieties simultaneously.

\medskip
Since the index set is rather complicated, we will write down the corresponding assumptions in each subsection accordingly.

	\medskip
	\subsection{ Step 1: Reduction to the sufficiency part of the nefness criterion.}
	
	\begin{lem}
		If the sufficiency part of the nefness criterion for $\mathcal{L}_Y(\lambda)$ (resp. $\mathcal{L}_X(\{k_l\})$) holds, then the sufficiency part of the ampleness criterion for $\mathcal{L}_Y(\lambda)$ (resp. $\mathcal{L}_X(\{k_l\})$) holds.
	\end{lem}	
	\noindent\emph{Proof: }The Hodge line bundle $\textnormal{det $\omega$}=\otimes_{i=1}^N\textnormal{det $\omega_i$}$ is ample, which corresponds to the weight $\lambda_0=\{k_l\}_{l=1}^N=\allowbreak\{1,1\dots,1\}$. Given any $\{k_l\}$ in the ample cone for $\mathcal{L}_X(\lambda)$, since the ample cone is an open subset of $\mathbb{R}^N$, there exists an $0<\varepsilon\ll 1$ such that
	\begin{equation}
		[\mathcal{L}_X(\{k_l\})]=[\mathcal{L}_X(\{k_l'\})]+ \varepsilon [\textnormal{det $\omega$}]
	\end{equation}
	in $\textnormal{(Pic $X$)}_{\mathbb{Q}}$ and $\{k_l'\}$ also lies in the ample cone. By Proposition \ref{Nef+Ample=Ample}, it suffices to prove $\mathcal{L}_X(\{k_l'\})$ is nef. Similarly, we can take an arbitrary $L$-regular dominant weight $\lambda$, the automorphic line bundle $\mathcal{L}_Y(\lambda)$ is therefore relatively ample (i.e., fiberwise ample) with respect to the projection $Y\to X$. Thus, for $M\gg1$, the line bundle $\mathcal{L}_Y(\lambda_0)\otimes \pi^\ast (\textnormal{det }\omega)^{\otimes M}$ is ample \cite{Lazarsfeld}. Again since the conjectured ample cone for $\mathcal{L}_Y(\lambda)$ is an open cone, there exists $0< \epsilon\ll 1$ such that
	\begin{equation}
		[\mathcal{L}_Y(\lambda)]=[\mathcal{L}_Y(\lambda')]+\epsilon [\mathcal{L}_Y(\lambda_0)\otimes \pi^\ast (\textnormal{det $\omega$})^{\otimes M}]
	\end{equation}
	in $(\textnormal{Pic $Y$})_{\mathbb{Q}}$ and $\lambda'$ lies in the conjectured ample cone. Tt suffices to prove $\mathcal{L}_Y(\lambda')$ is nef. $\hfill\square$
	
	\medskip
	Let $C_{\textnormal{nef}}$ denote the conjectured nef cone, i.e., the closure of the cone in Theorem \ref{Main Theorem - Introduction}. This is a convex cone in $\mathbb{R}^{nN}$. By Lemma \ref{General Picard Group Relations}, we have $[\textnormal{det }\H(\mathcal{A}/X)_i]=0$ in $(\textnormal{Pic $Y$})_\mathbb{Q}$, so we could view the nef cone as in $\mathbb{R}^{N(n-1)}$ by restricting to the intersection of hyperplanes
	\begin{equation}
		\sum_{j=1}^n k^i_j=0
	\end{equation}
	for all $1\le i\le N$. A convex cone is always generated by its boundary rays. For each case in Theorem \ref{Main Theorem - Introduction}, forcing $N(n-1)-1$ inequalities to be equalities and leaving exactly one inequality unchanged gives rise to a boundary ray of this cone. It is easy to see that any weight $\lambda$ lying in some boundary ray is a parallel weight for either $L$ or some $L^i_j$. This argument implies
	\begin{lem}\label{Reducing to partial flag spaces}
		In order to prove $\mathcal{L}_Y(\lambda)$ is nef for $\lambda\in C_{\textnormal{nef}}$, it suffices to prove for $\mathcal{L}_Y(\lambda)$ of either form:
		
		\medskip
		$\bullet$ $\lambda$ is a parallel weight for $L$, and $\mathcal{L}_Y(\lambda)=\pi^\ast \mathcal{L}_X(\lambda)=\pi^\ast\mathcal{L}_X(\{k_l\})$;
		
		$\bullet$ $\lambda$ is a parallel weight for $L^i_j$, and $\mathcal{L}_Y(\lambda)=\pi^\ast \mathcal{L}_{Y^i_j}(\lambda)=\pi^\ast\mathcal{L}_{Y^i_j}(\{k_l\};\alpha)$.$\hfill\square$
	\end{lem}
	
	By Lemma \ref{Nef+Ample=Ample}, we are reduced to proving the sufficency part of Theorem \ref{Nefness criterion for minimal partial flag space} and Theorem \ref{Nefness criterion for Shimura variety}.

	\medskip	
	\subsection{ Step 2: The induction base.}
	
	As described in subsection 7.1, the proof will heavily rely on inductions on the set $\mathfrak{I}$. In this part, we prove the base case for all induction steps:

\begin{lem} If the essential degree $t=0$, then the nefness criterion for $\mathcal{L}(\mathcal{\lambda})$ holds.
\end{lem}
\noindent\emph{Proof:} This follows directly from the regular dominance of $\lambda$ and the standard results of flag varieties in positive characteristics; see \cite{BrionFlagVarieties}, for example. $\hfill\square$

\begin{lem} If $T=\{i_1\}$ (that is, $t=1$) and $m_{i_1}=n-1$. For simplicity, we assume $i_1=1$. For any $1\le l\le n-1$ Let $Z\subseteq Y^1_l$ denote the subvariety defined by
\begin{equation}
    [\mathcal{F}^1_l\subseteq V^{-N}(\mathcal{F}^{1,(p^N)}_l)].
\end{equation}
Then $\mathcal{L}_{Y^1_l}(\lambda)|_Z=\mathcal{L}_{Y^1_l}(k_1;\alpha)|_Z$ is nef if
\begin{equation}
\begin{aligned}
    &k_{1}\ge \alpha,\\
    (&p^{N(n-l)}-1)\alpha\ge (p^{N(n-l-1)}-1)k_1.
\end{aligned}
\end{equation}
\end{lem}
\noindent\emph{Proof:} The goal is to prove that the intersection
\begin{equation}
    \big(k_1[\omega_1/\mathcal{F}^1_l]+\alpha[\mathcal{F}^1_l]\cdot C\big)\ge0
\end{equation}
for any complete curve $C\subseteq Z$ under the numercial condition above. We prove the statement of the proposition simultaneously for $1\le l \le n-1$.

	Consider a sequence of nested strata within $Z$:
	\begin{equation}
		Z_{n-l-1}\subseteq Z_{n-l-2}\subseteq \cdots\subseteq Z_1\subseteq Z_0:=Z,
	\end{equation}
	where $Z_j$ for $j\ge 1$ is defined by the condition
	\begin{equation}
		Z_j:=[\mathcal{F}_l^1\subseteq V^{-N}(\mathcal{F}_l^{1,(p^N)})\subseteq V^{-2N}(\mathcal{F}_l^{1,(p^{2N})})\subseteq\cdots \subseteq V^{-jN}(\mathcal{F}_l^{1,(p^{jN})})\subseteq\omega_1]\subseteq Z_{j-1}.
	\end{equation}

    \begin{rmk}
    In the above moduli interpretation, each $V^{-jN}(\mathcal{F}_l^{1,(p^{jN})})$ is not necessarily a vector bundle on $Z$, since $V^{-N}(\mathcal{F}_l^{1,(p^N)})$ fails to be contained in $\omega_1=\textnormal{Im }V^N$, so it is not well-defined to consider the inverse image under the (essential) Verschiebung map. However, on each stratum $Z_{j-1}$, we have $V^{-(j-1)N}(\mathcal{F}_l^{1,(p^{(j-1)N})})\subseteq \omega_1$, and $V^{-jN}(\mathcal{F}_l^{1,(p^{jN})})$ is a well-defined vector bundle of rank $l+j$ on $Z_{j-1}$. Therefore, the correct way to understand the condition is that the $Z_j$'s are consecutively defined by adding one new condition to $Z_{j-1}$.    
    \end{rmk}

	\begin{lem}
		Each $Z_j$ is a divisor in $Z_{j-1}$. The class of $Z_j$  in $\textnormal{Pic}(Z_{j-1})_{\mathbb{Q}}$ is equal to $$(p^{jN}-1)[\omega_1]-(p^{jN}-p^{(j-1)N})[\mathcal{F}_l^1].$$
	\end{lem}
	\emph{Proof: } $Z_j$ is the zero locus of the section
	\begin{equation}
		h_j:V^{-jN}(\mathcal{F}_l^{1,(p^{jN})})/V^{-(j-1)N}(\mathcal{F}_l^{1,(p^{(j-1)N})})\longrightarrow \H(A/X)_{\tilde\tau_1}/\omega_1.
	\end{equation}
    Here, the map is induced by the natural inclusion. Again, by Lemma \ref{Computing Picard Classes}, we have
	\begin{equation}
		\begin{aligned}
			c_1(h_j)&=c_1\big(\H(A/X)_{\tilde\tau_1}/\omega_1\big)-c_1\big(V^{-jN}(\mathcal{F}_l^{1,(p^{jN})})/V^{-{(j-1)N}}(\mathcal{F}^{1,(p^{(j-1)N})}_l)\big)\\
			&=-[\omega_1]-p^{(j-1)N}[V^{-1}(\mathcal{F}_l^{1,(p^N)})/\mathcal{F}_l]=(p^{jN}-1)[\omega_1]-(p^{jN}-p^{(j-1)N})[\mathcal{F}_l^1]
		\end{aligned}
	\end{equation}
	in Pic$(Z_{j-1})_\mathbb{Q}$.$\hfill\square$
	
	Now for any complete curve $C\subseteq Z$, let $j$ be the largest non-negative integer such that $C\subseteq Z_j$. If $j<n-l-1$, this means $C\nsubseteq Z(h_{j+1})$, and
	\begin{equation}
		\begin{aligned}
			[\mathcal{L}_{Y_l^1}(k_1;\alpha)]&=\frac{k_1-\alpha}{p^{(j+1)N}-p^{jN}}[h_{j+1}]+\big(k_1-\frac{p^{(j+1)N}-1}{p^{(j+1)N}-p^{jN}}(k_1-\alpha)\big)[\omega_1]\\
			&=\frac{k_1-\alpha}{p^{(j+1)N}-p^{jN}}[h_{j+1}]+\frac{(p^{(j+1)N}-1)\alpha-(p^{jN}-1)k_1}{p^{(j+1)N}-p^{jN}}[\omega_1].
		\end{aligned}
	\end{equation}
	By the assumption $(p^{(n-l)N}-1)k_l\ge (p^{(n-l-1)N}-1)k_{n-1}$, a direct computation implies that the coefficient of the second term is non-negative. Combining this with the nefness of $[\omega_1]$, we obtain
	\begin{equation}
		\begin{aligned}
			\big(\mathcal{L}_{Y_l}(k_{n-1}-k_l)\cdot C\big)&=\big(\frac{k_{n-1}-k_l}{p^{(j+1)N}-p^{jN}}[h_{j+1}]+\frac{(p^{(j+1)N}-1)k_l-(p^{jN}-1)k_{n-1}}{p^{(j+1)N}-p^{jN}}[\omega_1]\cdot C\big)\\
			&= \frac{k_{n-1}-k_l}{p^{(j+1)N}-p^{jN}}\big([h_{j+1}]\cdot C\big)+\frac{(p^{(j+1)N}-1)k_l-(p^{jN}-1)k_{n-1}}{p^{(j+1)N}-p^{jN}}\big([\omega_1]\cdot C\big)\ge0.
		\end{aligned}
	\end{equation}
	If $j=n-l-1$, i.e., $C\subseteq Z_{n-1-l}$. This stratum is defined by the condition
	\begin{equation}
		Z_{n-1-l}:=[\mathcal{F}^1_l\subseteq V^{-N}(\mathcal{F}_l^{1,(p^N)})\subseteq V^{-2N}(\mathcal{F}_l^{1,(p^{2N})})\subseteq\cdots \subseteq V^{-(n-1-l)N}(\mathcal{F}_l^{1,(p^{(n-l-1)N})})\subseteq\omega_1].
	\end{equation}
	Note that every time we take $V^{-N}$, the rank increases by 1, hence the last inclusion is in fact an equality of the vector bundles. By Lemma \ref{Computing Picard Classes}, we have the following in $(\textnormal{Pic $Z_{n-l-1}$})_{\mathbb{Q}}$
	\begin{equation}
		\begin{aligned}
			[\omega_1]&=[V^{-(n-1-l)N}(\mathcal{F}_l^{1,(p^{(n-l-1)N})})]=p^N[V^{-(n-2-l)N}(\mathcal{F}_l^{1,(p^{(n-l-2)N})})/\omega_1]\\&
			=\cdots\\
			&=p^{(n-1-l)N}[\mathcal{F}^1_l]-\big(p^N+p^{2N}+\cdots+p^{(n-1-l)N}[\omega_1]\big).
		\end{aligned}
	\end{equation}
	We can express $[\mathcal{F}^1_l]$ in terms of $[\omega_1]$ and plug into $[\mathcal{L}_{Y_l}(k_1;\alpha)]$:
	\begin{equation}
		\begin{aligned}
			[\mathcal{L}_{Y_l}(k_1;\alpha)]\big|_{Z_{n-l-1}}&=\bigg(k_1-\frac{p^{(n-l)N}-1}{p^{(n-1-l)N}(p^{N}-1)}(k_1-\alpha)\bigg)[\omega_1]\\
			&=\frac{1}{p^{(n-1-l)N}(p^{N}-1)}\Big((p^{(n-l)N}-1)\alpha-(p^{(n-1-l)N}-1)k_1\Big)[\omega_1].
		\end{aligned}
	\end{equation}
	The non-negativity of the coefficient is exactly our assumption. $\hfill\square$

We remark that in the case where $T=\{i_1\}$ and $m_{i_1}=1$, the proof of the corresponding statement is similar. Actually, by working on the $\tilde{\tau}_i^c$-side, the proof is exactly the same.

	\begin{lem} If $T=\{i_1\}$ (that is, $t=1$) and $m_{i_1}\ne \{1,n-1\}$. For simplicity, we assume $i_1=1$. For any $1\le l\le m_1$ Let $Z\subseteq Y^1_l$ denote the subvariety defined by
\begin{equation}
    [\mathcal{F}^1_l\subseteq V^{-N}(\mathcal{F}^{1,(p^N)}_l)].
\end{equation}
Then $\mathcal{L}_{Y^1_l}(\lambda)|_Z=\mathcal{L}_{Y^1_l}(k_1;\alpha)|_Z$ is nef if
\begin{equation}
\begin{aligned}
    &k_{1}\ge \alpha,\\
    &p^N\alpha\ge k_1.
\end{aligned}
\end{equation}
\end{lem}
	\noindent\emph{Proof:} By assumption, the set
	\begin{equation}
		\big\{\mathcal{F}^1_l,V^{-(N-1)}(\mathcal{F}_l^{1,(p^{N-1})}),\dots, V^{-1}(\mathcal{F}_l^{1,(p)})\big\}
	\end{equation}
	forms an $F,V$-chain. Hence, by Proposition \ref{Construction via F_{es},V_{es}-chain} $Z$, up to Frobenius twists, admits a proper morphism to an auxiliary Shimura variety $X'$. $X'$ also has essential degree $t'=1$, with $T'=\{2\}$ and $m_{2}'=m_1$. We can compute the first Chern class of $\omega_2'$ on $Z$ as follows:
	\begin{equation}
		\begin{aligned}
			p[\omega'_2]&=[\omega'^{(p)}_2]=[V(\H(\mathcal{A}'/X')_{\tilde{\tau}_1})]=[V{\mathcal{F}^1_l}/pV^{-(N-1)}({\mathcal{F}}_l^{1,(p^{N-1})})^{(p)}]\\
			&=\big[\frac{V\mathcal{F}^1_l}{Vp\H(\mathcal{A}/X)_{\tilde\tau_1}}\big]+\big[\frac{Vp\H(\mathcal{A}/X)_{\tilde\tau_1}}{pV^{-(N-1)}(\mathcal{F}_l^{1,(p^{N-1})})^{(p)}}\big]\\
			&=[\mathcal{F}^1_l]-p^{N-1}[V^{-1}(\mathcal{F}^{1,(p)}_l)] \textnormal{ (Because the Verschiebung map are all isomorphisms except one)}\\
            &=[\mathcal{F}^1_l]+p^N[\omega_1/\mathcal{F}^1_l].
		\end{aligned}
	\end{equation}
    The terms in the computation should be viewed as their pointwise liftings to $W(k)$-lattices in $\tilde{D}(\mathcal{A})[\frac{1}{p}]$. One could also deduce the class by a careful analysis of the (kernel and image of the) universal isogeny. Therefore, we could write
	\begin{equation}
		\begin{aligned}
			[\mathcal{L}_{Y^1_l}(k_1;\alpha)]\big|_Z&=k_1[\omega_1/\mathcal{F}^1_l]+\alpha[\mathcal{F}^1_l]=\frac{p}{p^N-1}(k_1-\alpha)[\omega_2']+\big(k_1-\frac{p^N}{p^N-1}(k_1-\alpha)\big)[\omega_1]\\
			&=\frac{p^N}{p^N-1}(k_1-\alpha)[\omega_2']+\frac{p^N\alpha-k_1}{p^N-1}[\omega_1].
		\end{aligned}
	\end{equation}
	Our assumption ensures that the coefficients are nonnegative, and since both $[\omega_1]$ and $[\omega_2']$ are classes of Hodge bundles on different Shimura varieties, we conclude that the restriction of $\mathcal{L}_{Y^1_l}(k_1;\alpha)$ to $Z$ is nef. The case $l>m_1$ follows from a similar argument.

	\medskip
	\subsection{ Step 3: The sufficiency part of $\mathcal{L}_X(\lambda)$}
	In this step, our aim is to prove the following class
	\begin{equation}
		[\mathcal{L}_X(\lambda)]=[\mathcal{L}_X(\{k_{s}\})]=\sum_{s=1}^t k_s[\omega_{i_s}]
	\end{equation}
	is a nef class under the hypothesis $p^{a(i_s)}k_s\ge k_{s+1}$ for all $1\le s\le t$. Here $T=\{{i_1},{i_2},\dots, {i_s}\}$ with $i_1<i_2<\cdots<i_t$, and the subscripts of $k$ are viewed as modulo $t$. Notice that if $t=1$, then we need to show $\mathcal{L}_X(\lambda)=k_1[\omega_{i_1}]$ is nef under the condition $k_1\ge0$. This follows directly from the ampleness of Hodge line bundle \cite{MadapusiCompactification}.
    
    From now on, we assume the essential degree $t\ge2$. The induction hypothesis in subsection 7.1, in particular, implies the following assumptions.
    
    \noindent\textbf{Induction hypothesis:} Theorem \ref{Nefness criterion for minimal partial flag space} holds for any minimal partial flag space over any Shimura variety with essential set $T'$ and essential degree $t'=\# T'$ such that either
    \begin{itemize}
        \item We have $i\in T'$ and $t'\le t-1$, or
        \item We have $i\notin T'$ and $t'\le t-2$.
    \end{itemize}
    

	For unitary Shimura varieties with different signatures, it is in general not possible to talk about the relations of their nef cones, since they have different Hodge parabolics $P$'s, and $\lambda$ corresponds to characters of different groups. However, for a fixed essential set $T$, note that $\lambda$ corresponds to a tuple of numbers $\{k_s\}$, we can viewed it as a vector in $\mathbb{R}^t$. Under this identification, it is natural to define the \emph{common} nef cone, denoted by $C_T$, as the intersection of the nef cones of all unitary Shimura varieties with essential set $T$.

	Let $C_{\textnormal{SV}}$ be the convex cone in $\mathbb{R}^t$ defined by the condition
	\begin{equation}
		p^{a(i_s)}k_s\ge k_{s+1},\qquad \textnormal{$1\le s\le t$}.
	\end{equation}
	By the necessity part of the ampleness criterion, we know $C_T\subseteq C_{\textnormal{SV}}$, and the goal of this section is to prove the opposite inclusion.
	
	\begin{lem}
		$C_{\textnormal{SV}}$ is spanned by the weights
		\begin{equation}
			\begin{aligned}
				&(1,p^{a(i_1)},p^{a(i_1)+a(i_2)},\dots,p^{a(i_1)+\cdots+a(i_{t-1})}),\\
				&(p^{a(i_2)+\cdots+a(i_t)},1,p^{a(i_2)},\dots,p^{a(i_2)+\cdots+a(i_{t-1})}),\\
				&\cdots\cdots\\
				&(p^{a(i_t)},p^{a(i_t)+a(i_1)},\dots,p^{a(i_t)+\cdots+a(i_{t-2})},1).
			\end{aligned}
		\end{equation}
	\end{lem}
	\noindent\emph{Proof: }This follows from a straightforward linear algebra.$\hfill\square$
	
	In order to prove $C_{SV}\subseteq C_T$, it suffices to show each boundary weight lies in $C_{T}$. Without loss of generality, we work with $(1,p^{a(i_1)},p^{a(i_1)+a(i_2)},\dots,p^{a(i_1)+\cdots+a(i_{t-1})})$ in the sequel. Let $\tilde{C}_{T},\tilde{C}_{\textnormal{SV}}\in \mathbb{R}^t$  be the inverse images of $C_T,C_{\textnormal{SV}}$ under the map
	\begin{equation}
		\begin{aligned}
			\xi:\mathbb{R}^t&\longrightarrow \mathbb{R}^t\\
			(k_1,\dots,k_n)&\longmapsto (k_1,k_2p^{a(i_1)},\dots,k_dp^{a(i_1)+\cdots +a(i_{t-1})}).
		\end{aligned}
	\end{equation}
	The goal is to show the weight $(1,\dots,1)\in \tilde{C}_{T}$. Note that the ampleness of the Hodge bundle implies 
	\begin{equation}
		(1,p^{-a(i_1)},\dots, p^{-(a({i_1})+\cdots+a(i_{t-1}))})\in \tilde{C}_T.
	\end{equation}
	Thus, $\tilde{C}_T$ is non-empty. The following lemma is the heart of our proof. It produces new nef weights from the known ones, hence enlarging $\tilde{C}_T$:
	
	\begin{lem}\label{Improvement of nef estimation}
		Let $\kappa=(k_1,\dots,k_t)\in \tilde{C}_{\textnormal{SV}}$. If $\kappa\in \tilde{C}_T$, then for $j=1,\dots,t-1$, $\kappa_j'=(k_1,\dots,\frac{k_j+k_{j+1}}{2},\frac{k_j+k_{j+1}}{2},\dots,k_t) \in\tilde{C}_T$.
	\end{lem}
	\noindent\emph{Proof: } Let $C$ be a complete curve, we aim to prove
	\begin{equation}
		([\mathcal{L}_{X}(\xi(\kappa_j'))]\cdot C)\ge0.
	\end{equation}
	
	\begin{definition}
		We call a closed subscheme of $X$ defined by conditions of the form 		
		$$
		\big\{\textnormal{dim }(\omega_{i_l}\cap \textnormal{Ker }V^{a(i_l)}_{es})\ge r_{i_l},1\le l\le t\big\}
		$$
		a Schubert cell of $X$.
	\end{definition}
	
	Note that $\omega_{i_l}$ has rank $m_{i_l}$ and $\textnormal{Ker }V_{i_l,es}^{a(i_l)}$ has rank $n_{i_{l+1}}$, so $r_{i_l}$ must be an integer within the interval
	\begin{equation}
		[\textnormal{Max}\{0,m_{i_l}-m_{i_{l+1}}\},\textnormal{Min}\{m_{i_l},n_{i_{l+1}}\}].
	\end{equation}
    Note that if $r_{i_l}$ is equal to the lower bound, then the condition at $i_l$ is actually trivial.
	
	\begin{rmk}
		The Scubert cells also satisfy properties such as non-emptiness, dimension formula, closure relations, etc, and forms a good stratification by considering the interior of each Schubert cell. On the other hand, since the Ekedahl--Oort stratification captures the relative position of the Hodge and conjugate filtrations, as well as a comparison isomorphism between the graded pieces, it is natural that Ekedahl--Oort stratification is finer than this Schubert-type stratification. In \cite{KoskivirtaNormalization}, the author studied when the closure of an Ekedahl--Oort stratum agrees with a Schubert cell. We will not delve deeply into this discussion in this paper. It is enough to only consider the Schubert-type stratification for our purpose in this step.
	\end{rmk}
    
	Let $Z$ be a Schubert cell of $X$. There are two possibilities for the conditions for $Z$: Either there exists some $l$, such that $r_{i_l}=\textnormal{Min}\{m_{i_l},m_{i_{l+1}}\}$, or there does not exist such an $l$. 

    \begin{lem}\label{Restriction to the Schubert Cell of the First Case}
        If $Z\subset X$ is a Schubert cell of the first case, then $\mathcal{L}_X(\lambda)|_Z$ is nef for any $\lambda\in C_{\textnormal{SV}}$.
    \end{lem}

    \noindent\emph{Proof}. In this first case, the condition implies that $\omega_{i_l}\subseteq \textnormal{Ker }V^{a({i_l})}$ or $\textnormal{Ker }V^{a(i_l)}\subseteq \omega_{i_l}$, depending on which vector bundle has the larger rank. Without loss of generality, suppose that we have $\omega_{i_l}\supseteq \textnormal{Ker }V^{a(i_l)}$ for some $l$. The sequence
	\begin{equation}
		\H(\mathcal{A}/X)_1,\dots,\H(\mathcal{A}/X)_{i_l-1},\omega_{i_l}, V\omega_{i_l},\dots,V^{a(i_l)-1}\omega_{i_l},\H(\mathcal{A}/X)_{i_{l+1}},\dots,\H(\mathcal{A}/X)_N
	\end{equation}
	forms an $F,V$-chain. By Theorem \ref{Description of strata}, the Schubert cell $Z$ is, up to appropriate Frobenius twists, isomorphic to the minimal partial flag space $\tilde{Y}^{i_{l+1}}_{m_{i_{l+1}}}$ over the auxiliary unitary Shimura variety $\tilde{X}$. (We are actually a bit of sloppy here: Each sheaf $V^r\omega_{i_l}$ is a subbundle of $\H(\mathcal{A}/X)_{i_{l+1}+r}^{(p^{r})}$, so the above chain is actually not an $F,V$-chain. Nevertheless, there are two possible ways to modify this issue: (1) One could take the chain
    \begin{equation}
        0,\dots,0, \textnormal{Ker }V^{a(i_l)},\dots, \textnormal{Ker }V,0,\dots,0.
    \end{equation}
    Then each $\textnormal{Ker }V^{r}\subseteq \H(\mathcal{A}/X)_{i_{l+1}-r}$ and one checks that this forms an $F,V$-chain. (2) We could generalize the notion of $F,V$-chain and allow each $\mathcal{E}_i$ to be a subbundle of $\H(\mathcal{A}/X)_i^{(p^{u_i})}$ for some $p^{u_i}$-twist. Since being nefness is an equivalent condition pullbacks under proper surjective morphisms, it is harmless even if we allow these Frobenius twists. The reason that we choose the above sequence here is that it is intuitively more natural)

    We collect some numerical information of $\tilde{Y}^{i_{l+1}}_{m_{i_{l+1}}}$ provided by Theorem \ref{Description of strata}: Let $0\subseteq\tilde{\omega}_{i_{l+1}}\subseteq\mathcal{J}\subseteq\H(\tilde{\mathcal{A}}/\tilde{X})_{i_{l+1}}$ denote the refinement of the Hodge filtration at $\tau_{i_{l+1}}$. Recall that each term in the above $F,V$-chain corresponds to the image of $\H(\tilde{\mathcal{A}}/\tilde{X})_j$ for some $j$, and $\omega_{i_{l+1}}$ corresponds to the image of $\mathcal{J}$, under the universal $p$-isogeny $\phi:\tilde{A}\to A$ in the moduli problem of the correspondence.
	
	\begin{lem}\label{Picard Classes Computation in Chapter 7} Let $(\tilde{m}_j,\tilde{n}_j)_{j\in\Sigma_\infty}$ be the signature of $\tilde{X}$. Let $\tilde\omega_j=\omega_{\tilde{A}^\vee/\tilde{X},\tilde\tau_j}$ for all $j$.
		\begin{enumerate}
			\item (Signature changes) $\tilde{m}_j=m_j$ for $j\notin\{i_l,i_{l+1}\}$, $\tilde{m}_{i_l}=n$, and $\tilde{m}_{i_{l+1}}=m_{i_{l+1}}-n_{i_l}$. In particular, the essential set $\tilde{T}=T\backslash \{i_l\}$ or $\tilde{T}=T\backslash\{i_l,i_{l+1}\}$, depending on whether $m_{i_{l+1}}=n_{i_l}$.
			\item (Classes in Picard group) In $(\textnormal{Pic }\tilde{X})_{\mathbb{Q}}$, we have  $[\tilde{\omega}_{i_l}]=0$, $[\tilde{\omega}_{i_{l+1}}]=\frac{1}{p^{a(i_l)}}[\omega_{i_l}]+[\omega_{i_{l+1}}]$, and $[\omega_j]=[\tilde{\omega}_{j}]$ for $j\notin \{i_l,i_{l+1}\}$. Moreover, in $(\textnormal{Pic }\tilde{Y}^{i_{l+1}}_{m_{i_{l+1}}})_{\mathbb{Q}}$, we have $[\mathcal{J}]=[\omega_{i_{l+1}}]$.
		\end{enumerate}
	\end{lem}
    \emph{Proof.} This is a direct application of Lemma \ref{Computing Picard Classes}. We only explain the computation for $[\tilde{\omega}_{i_{l+1}}]$ here. Recall that
    \begin{equation}
        \begin{aligned}
            &\phi_{i_l,\ast}(\H(\tilde{\mathcal{A}}/X)_i)=\omega_{i_l},\\
            &\phi_{i_{l+1},\ast} \textnormal{ is an isomorphism}.
        \end{aligned}
    \end{equation}
    We have 
    \begin{equation}
        [\tilde{\omega}_{i_{l+1}}]=[\phi_{i_{l+1},\ast}(\tilde{\omega}_{i_{l+1}})]=\frac{1}{p^{a(i_l)}}[V^{a(i_l)}\omega_i].
    \end{equation}
    A direct application of Lemma \ref{General Picard Group Relations} to the short exact sequence
    \begin{equation}
        0\longrightarrow \textnormal{Ker }V^{a(i_l)}\longrightarrow \omega_i\longrightarrow V^{a(i_l)}\omega_i\longrightarrow0
    \end{equation}
    yields the result. The rest of the computations are similar.    $\hfill\square$

	Now we have
	\begin{equation}
		\begin{aligned}
			[\mathcal{L}_X(\lambda)]\big|_Z&=\sum_{j=1}^{t}k_j[\tilde{\omega}_{i_j}]=p^{a(i_l)}k_{l}[\tilde{\omega}_{i_{l+1}}]+(k_{l+1}-p^{a(i_l)}k_{l})[\mathcal{J}]+\sum_{\substack{j=1 \\ j\ne l,l+1}}^{t} k_j[\tilde{\omega}_{i_j}]\\
			&=p^{a(i_l)}k_l[\tilde{\omega}_{i_{l+1}}/\mathcal{J}]+k_{l+1}[\mathcal{J}]+\sum_{\substack{j=1 \\ j\ne l,l+1}}^{t} k_j[\tilde{\omega}_{i_j}].
		\end{aligned}
	\end{equation}
	By the above lemma, either
    \begin{itemize}
        \item The auxiliary Shimura variety $\tilde{X}$ has essential degree $t'=t-1$, and $\tilde{Y}_{m_{i_{l+1}}}^{i_{l+1}}$ is flagged at an place in the essential set, or
        \item The auxiliary Shimura variety $\tilde{X}$ has essential degree $t'=t-2$, and $\tilde{Y}_{m_{i_{l+1}}}^{i_{l+1}}$ is flagged at an place not in the essential set.
    \end{itemize}
    Both satisfies the Induction hypothesis at the beginning of this subsection. By induction, it suffices to check $\mathcal{L}_X(\lambda)|_Z$ satisfies the nef condition for $\tilde{Y}^{i_{l+1}}_{m_{i_{l+1}}}$. The inequalities away from $i_l,i_{l+1}$ are the same as the original ones. The dominance condition $-k_{l+1}\ge -p^{a(i_l)}k_l$ for $\mathcal{J}$ follows from the assumption. At $i_l$ and $i_{l+1}$, the rest of the inequalities follow from a direct discussion:
	
	\medskip
	\noindent (1) If $t\ge 3$ and $m_{i_{l+1}}\ne n_{i_l}$. Then $t'=t-1\ge2$ and
	\begin{equation}
		p^{a(i_{l-1})+a(i_l)}k_{l-1}\ge p^{a(i_l)}k_l,\qquad p^{a(i_l)}k_l\ge k_{l+1}, \qquad p^{a(i_{l+1})}k_{l+1}\ge k_{l+2}.
	\end{equation}

    \noindent (2 )If $t=2$ and $m_{i_{l+1}}\ne n_{i_l}$, then $\tilde{T}=\{i_{l+1}\}$ has cardinality 1 and $\tilde{Y}^{i_{l+1}}_{m_{i_l}}$ is flagged at the place $\tilde\tau_{i_{l+1}}$. If moreover $\tilde{m}_{i_{l+1}}\in\{1,n-1\}$, then one checks 
	\begin{equation}
		p^{a(i_l)}k_l\ge k_{l+1},\qquad ((p^N)^r-1)k_{l+1}\ge ((p^N)^{r-1}-1)p^{a(i_l)}k_l
	\end{equation}
	holds for any $r\ge 1$. The latter inequality follows from (note that $a(i_l)+a(i_{l+1})=N$ in this case)
	\begin{equation}
		\frac{p^{Nr}-1}{(p^{N(r-1)}-1)p^{a(i_l)}}\ge p^{a(i_{l+1})}.
	\end{equation}

    \noindent (3) If $t=2$ and $m_{i_{l+1}}\ne n_{i_l}$, then $\tilde{T}=\{i_{l+1}\}$ has cardinality 1 and $\tilde{Y}^{i_{l+1}}_{m_{i_l}}$ is flagged at the place $\tilde\tau_{i_{l+1}}$. If moreover $\tilde{m}_{i_l}\notin\{1,n-1\}$, then one checks (note that $N=a(i_l)+a(i_{l+1})$)
    \begin{equation}
		p^Nk_{l+1}=p^{a(i_l)+a(i_{l+1})}k_{l+1}\ge p^{a(i_l)}k_l.
	\end{equation}

    \noindent (4) If $t\ge 4$ and $m_{i_{l+1}}=n_{i_l}$, then $\tilde{T}$ has cardinality $\ge 2$. We have $[\tilde{\omega}_{i_{l+1}}]=[\H(\tilde{\mathcal{A}}/X)_{\tilde\tau_{i_{l+1}}}]=0$. The inequality
    \begin{equation}
		p^{a(i_{l-1})+a(i_l)+a(i_{l+1})}\big(k_{l-1}+\frac{1}{p^{a(i_{l-1})+a(i_l)}}(k_{l+1}-p^{a(i_l)}k_l)\big)\ge k_{l+2}
	\end{equation}
    holds because
	\begin{equation}
		p^{a(i_{l-1})}k_{l-1}\ge k_l\quad \textnormal{and}\quad p^{a(i_{l+1})}k_{l+1}\ge k_{l+2}.
	\end{equation}

	

    \noindent (5) If $t=3$ and $m_{i_{l+1}}=n_{i_l}$, then $\tilde{T}=\{i_{l+2}\}$ has cardinality 1. Assume further that $\tilde{m}_{i_{l+2}}\notin  \{1,n-1\}$, then the inequality
    \begin{equation}
		p^N\big(k_{l+2}+\frac{1}{p^{a(i_l)+a(i_{l-1})}}(k_{l+1}-p^{a(i_l)}k_l)\big)\ge k_{l+2}
	\end{equation}
	holds because $N=a(i_l)+a(i_{l+1})+a(i_{l+2})$ and the same reason as (4).

	

    \noindent (6) If $t=3$ and $m_{i_{l+1}}=n_{i_l}$, then $\tilde{T}=\{i_{l+2}\}$ has cardinality 1. Assume further that $\tilde{m}_{i_{l+2}}\notin  \{1,n-1\}$, then it suffices to check
    \begin{equation}
		((p^N)^r-1)\big(k_{l+2}+\frac{1}{p^{a(i_{l+2})+a(i_l)}}(k_{l+1}-p^{a(i_l)}k_l)\big)\ge ((p^N)^{r-1}-1)k_{l+2}
	\end{equation}
    for any $1\le r\le n-1$. This is equivalent to (note that $N=a(i_l)+a(i_{l+1})+a(i_{l+2})$ in this case)
	\begin{equation}
		p^{a(i_{l+2})+a(i_l)}(p^{Nr}-p^{N(r-1)})k_{l+2}+(p^{Nr}-1)k_{l+1}\ge (p^{Nr}-1)p^{a(i_l)}k_l.
	\end{equation}
    The left hand side
    \begin{equation}
		\begin{aligned}
			p^{a(i_{l+2})+a(i_l)}(p^{Nr}-p^{N(r-1)})k_{l+2}+(p^{Nr}-1)k_{l+1}\ge p^{a(i_l)}(p^{Nr}-p^{N(r-1)})k_l+p^{-a(i_{l+1})-a(i_{l+2})}(p^{Nr}-1)k_{l}.
		\end{aligned}
	\end{equation}
    One checks that it is larger than the right hand side by comparing the coefficients.

    \noindent (7) If $t=2$ and $m_{i_{l+1}}=n_{i_l}$, then $\tilde{T}=\emptyset$ and $\tilde{X}$ is set of discrete points. The nefness corresponds to 
    \begin{equation}
        p^{a(i_l)}k_l\ge k_{l+1}.
    \end{equation}
    
	Conclusively, we have shown that for any Schubert cell $Z$ of the first case, the restriction of $\mathcal{L}_X(\lambda)$ to $Z$ is nef for any $\lambda\in C_{\textnormal{SV}}$. In particular, $\mathcal{L}_X(\xi(\kappa_j'))|_Z$ is nef. This completes the proof of Lemma \ref{Restriction to the Schubert Cell of the First Case}. $\hfill\square$
	
	\medskip
	Next, assume $Z$ is a Schubert cell given by the conditions
    \begin{equation}
        \textnormal{dim }(\omega_{i_l}\cap \textnormal{Ker }V_{es}^{a(i_l)})\ge r_{i_l}, \qquad 1\le l\le t.
    \end{equation} We moreover assume that it is of the second case. We first construct an auxliary scheme.
	
	\medskip
	\noindent\textbf{Construction:} Let $i_l\in T$, 
	
	$\bullet$ Take $\mathcal{F}^{i_l}$ to be a subbundle of $\H(A/X)_{i_l}$ of rank $r_{i_l}$ contained in both $\omega_{i_l}$ and $\textnormal{Ker }V_{es}^{a(i_l)}$. If $r_{i_l}=0$, then we make the trivial choice $\mathcal{F}^{i_l}=0$;
	
	$\bullet$ Take $\mathcal{E}^{i_l}$ to be a subbundle of $\H(A/X)_{i_l}$ of rank $m_{i_l}+n_{i_{l+1}}-r_{i_l}$ that contains both $\omega_{i_l}$ and $\textnormal{Ker }V_{es}^{a(i_l)}$. If $r_{i_l}=m_{i_l}-m_{i_{l+1}}$, then we make the trivial choice $\mathcal{E}^{i_l}=\H(A/X)_{i_l}$.
	
	If fact, we should write these auxiliary vector bundles as $\mathcal{F}^{i_l}_{r_{i_l}}$ and $\mathcal{F}^{i_l}_{m_{i_l}+n_{i_{l+1}}-r_{i_l}}$ by our standard convention. However, since we only pick out one extra subbundle for each essential place, we prefer to record them in a simpler form. 
	
	For each $1\le l\le t$, the datum $(A,\lambda,\eta,\iota,\mathcal{F}^{i_l}\subseteq \omega_{i_l}\subseteq \mathcal{E}^{i_l})$ corresponds to a stratum $Z_{i_l}$ in the partial flag space $P^{i_l,i_l}_{r_{i_l},m_{i_l}+n_{i_{l+1}}-r_{i_l}}$. Furthermore, the forgetful map $Z_{i_l}\to Z$ is projective and surjective. We are reduced to prove that the pullback of $\mathcal{L}_X(\lambda)$ to $Z_{i_l}$ is nef.
	
	The sequences
	\begin{equation}
		\begin{aligned}
        &0,\dots,0,\mathcal{F}^{i_l},V_{es}\mathcal{F}^{i_l},\dots,V_{es}^{a(i_l)-1}\mathcal{F}^{i_l},0,\dots,0,\\
			&\H(\mathcal{A}/X)_1,\dots,\H(\mathcal{A}/X)_{i_l-1},\mathcal{E}^{i_l},V\mathcal{E}^{i_l},\dots,V^{a(i_l)-1}\mathcal{E}^{i_l},\H(\mathcal{A}/X)_{i_{l+1}},\dots,\H(\mathcal{A}/X)_N		
		\end{aligned}
	\end{equation}
	are both $F,V$-chains (Here again, we need to allow Frobenius twists in the definition of $F,V$-chains). By Proposition \ref*{Description of strata}, we obtain correspondences of strata in partial flag spaces of unitary Shimura varieties for both $F,V$-chains.
    \begin{center}
		\begin{tikzpicture}[node distance=12pt and 1.5cm, auto]
			\node(1) [align=center] {$W'$};
			\node(2) [align=center, below left=of 1] {$Z_{i_l}$};
			\node(3) [align=center, below right=of 1] {$Z'$};
			\node(4) [align=center, left=of 2] {$X$};
			\node(5) [align=center, right=of 3] {$X'.$};
			
			\draw [->] (1) -- node[midway, above] {$pr_1$} (2);
			\draw [->] (1) -- node[midway, above] {$pr_2$} (3);
			\draw [<-] (4) -- node[midway, above] {} (2);
			\draw [->] (3) -- node[midway, above] {} (5);
		\end{tikzpicture}
	\end{center}
        \begin{center}
		\begin{tikzpicture}[node distance=12pt and 1.5cm, auto]
			\node(1) [align=center] {$W''$};
			\node(2) [align=center, below left=of 1] {$Z_{i_l}$};
			\node(3) [align=center, below right=of 1] {$Z''$};
			\node(4) [align=center, left=of 2] {$X$};
			\node(5) [align=center, right=of 3] {$X''.$};
			
			\draw [->] (1) -- node[midway, above] {$pr_1$} (2);
			\draw [->] (1) -- node[midway, above] {$pr_2$} (3);
			\draw [<-] (4) -- node[midway, above] {} (2);
			\draw [->] (3) -- node[midway, above] {} (5);
		\end{tikzpicture}
	\end{center}
    We denote the auxiliary unitary Shimura varieties by $X'$ and $X''$, and the universal abelian schemes by $\mathcal{A}'/X'$ and $\mathcal{A}''/X''$, respectively. We recall the signature conditions for $X'$ and $X''$ as follows.
	
	\begin{lem} $\mathcal{A}'/X'$ has signature condition $m_j'=m_j$ for $j\notin\{i_l,i_{l+1}\}$, $m_{i_l}'=m_{i_l}-r_{i_l}$ and $m'_{i_{l+1}}=m_{i_{l+1}}+r_{i_{l}}$; Similarly, $\mathcal{A}''/X''$ has signature condition $m_j''=m_j$ for $j\notin\{i_l,i_{l+1}\}$, $m_{i_l}''=m_{i_l+1}+r_{i_l}$ and $m_{i_{l+1}}''=m_{i_l}-r_{i_l}$. In particular, both $X'$ and $X''$ have essential set $T$.
	\end{lem}
	\noindent\emph{Proof: }This follows directly from the construction via Dieudonn\'e modules in Theorem \ref{Description of strata}. $\hfill\square$
	
	\begin{lem} In $(\textnormal{Pic $X'$})_{\mathbb{Q}}$, we have
		\begin{equation}
			\begin{aligned}
				&[\omega'_{j}]=[\omega_j]=0,\quad \textnormal{if }j\notin T,\\
				&[\omega'_{j}]=[\omega_j],\quad \quad\ \ \ \textnormal{if }j\in T \textnormal{ but }j\ne i_l,i_{l+1},\\
				&[\omega'_{i_l}]=[\omega_i]-[\mathcal{F}^{i_l}],\\
				&[\omega'_{i_{l+1}}]=\frac{1}{p^{a(i_l)}}[\mathcal{F}^{i_l}]+[\omega_{i_{l+1}}].
			\end{aligned}
		\end{equation}
		The statement for $[\omega''_j]$ in $(\textnormal{Pic $X''$})_{\mathbb{Q}}$ is similar, with $\mathcal{F}^{i_l}$ replaced by $\mathcal{E}^{i_l}$ in each equation.
	\end{lem}
	\noindent\emph{Proof:} These lemmas follow from a direct computation similar to Lemma \ref{Picard Classes Computation in Chapter 7}. $\hfill\square$
	
	For each $1\le l\le t$, we can construct a strata Hasse invariant on $Z$:
	$$
	h_{i_l}:\textnormal{det }(\omega_{i_l}/\mathcal{F}^{i_l})\longrightarrow \textnormal{det } (\mathcal{E}^{i_l}/\textnormal{Ker }V_{es}^{a(i_l)}).
	$$
	The zero locus of $h_{i_l}$ is a lower-dimensional Schubert cell by the moduli description. In $(\textnormal{Pic $Z$})_{\mathbb{Q}}$, we have
	\begin{equation}
		[h_{i_l}]=p^{a(i_l)}[\omega_{i_{l+1}}]-[\omega_{i_l}]+[\mathcal{F}^{i_l}]+[\mathcal{E}^{i_l}],
	\end{equation}
	which follows from a direct application of Lemma \ref{Computing Picard Classes}.

	Let $C$ be a complete curve in $X$. Let $Z$ be the Schubert stratum of lowest dimension containing $C$. We can perform the above constructions for this particular $Z$. 
	Note that the auxiliary unitary Shimura variety $X'$ has the same essential set $T$ as $X$, we deduce that for any $\kappa\in \tilde{C}_T$ and $1\le l\le t-1$, the class
	\begin{equation}
		\begin{aligned}
			[\mathcal{L}_{X'}(\xi(\kappa))]&=\sum_{j\in T} k_jp^{a(i_1)+\cdots+a(i_{j-1})} [\omega'_{i_j}]\\
			&=\sum_{j\in T\backslash\{l,l+1\}}k_jp^{a(i_1)+\cdots+a(i_{j-1})} [\omega_{i_j}]+k_lp^{a(i_1)+\cdots+a(i_{l-1})}([\omega_{i_l}]-[\mathcal{F}^{i_l}])\\
			&\qquad\qquad\qquad\qquad\ \ \qquad\qquad\qquad +k_{l+1}p^{a(i_1)+\cdots+a(i_l)}\big(\frac{1}{p^{a(i_l)}}[\mathcal{F}^{i_l}]+[\omega_{i_{l+1}}]\big)
		\end{aligned}
	\end{equation}
	is nef. Similarly for $X''$ with $[\mathcal{F}^{i_j}]$ replaced by $[\mathcal{E}^{i_j}]$. As before, we have the Hasse invariant 
	$$
	[h_{i_l}]=p^{a(i_l)}[\omega_{i_{l+1}}]-[\omega_{i_l}]+[\mathcal{F}^{i_l}]+[\mathcal{E}^{i_l}].
	$$ 
	Recall that in the statement of Lemma \ref{Improvement of nef estimation}, we expect to prove $\kappa_l'\in \tilde{C}_T$. We can eliminate the classes $[\mathcal{F}^{i_l}]$ and $[\mathcal{E}^{i_l}]$ by
	\begin{equation}
		\begin{aligned}
			[\mathcal{L}_X(\xi(\kappa_l'))]=\frac{1}{2}[\mathcal{L}_{X'}(\xi(\kappa))]+\frac{1}{2}[\mathcal{L}_{X''}(\xi(\kappa))]+\frac{1}{2}(k_l-k_{l+1}) p^{a(i_1)+\cdots+a(i_{l-1})} [h_{i_l}].
		\end{aligned}
	\end{equation}
	By the choice of $Z$, we know that the curve $C\nsubseteq Z(h_{i_l})$. Thus,
	\begin{equation}
		\begin{aligned}
			&\big([\mathcal{L}_X(\xi(\kappa_j'))]\cdot C\big)=\frac{1}{2}\big([\mathcal{L}_{X'}(\xi(\kappa))]\cdot C\big)+\frac{1}{2}\big([\mathcal{L}_{X''}(\xi(\kappa))]\cdot C\big) +\frac{1}{2}(k_j-k_{j+1})p^{a(i_1)+\cdots+a(i_{j-1})}\big([h_{i_j}]\cdot C\big)\ge0.
		\end{aligned}
	\end{equation}
    Here, from the definition of $\xi$, we see that the cone $\tilde{C}_T$ is given by the condition
    \begin{equation}
        k_1\ge k_2\ge\cdots\ge k_t,\qquad p^N k_t\ge k_1.
    \end{equation}
    This explains the non-negativity of the intersection number.
    
	Combining the two cases together, we have shown that for any complete curve $C\subseteq X$,
	\begin{equation}
		\big(\mathcal{L}_X(\xi(\kappa_l))\cdot C\big)\ge 0,
	\end{equation}
	i.e., $\kappa_l\in \tilde{C}_{T}$, for any $1\le l\le t-1$. This completes the proof of Lemma \ref{Improvement of nef estimation}. $\hfill\square$
	
	\medskip
	Recall that we have the ampleness of the Hodge bundle, this is equivalent to
	\begin{equation}
		(1,p^{-a(i_1)},\dots,p^{-(a(i_1)+a(i_2)+\cdots +a(i_{l-1}))})\in \tilde{C}_T.
	\end{equation}
	Let $\tilde{C}_0$ be the ray corresponding to the above class. We utilize the above lemma to improve our estimation of the cone of nef weights.
	
    \medskip
	\noindent\emph{Proof of $\tilde{C}_T=\tilde{C}_{\textnormal{SV}}$:} It suffices to prove that the left hand side contains the right hand side. Starting from the ray $\tilde{C}_0$, we apply Lemma \ref{Improvement of nef estimation} to all the tuples in $\tilde{C}_0$, and take the convex hull of $\tilde{C}_0$ and the resulting set. We denote this new cone by $\tilde{C}_1$. By construction $\tilde{C}_1\subseteq \tilde{C}_{T}\subseteq \tilde{C}_{\textnormal{SV}}$. We then repeat this procedure for $\tilde{C}_1$ to obtain $\tilde{C}_2$, and continue this process iteratively. Consequently, we get a sequence of increasing convex subcones of $\tilde{C}_{T}$:
	\begin{equation}
		\tilde{C}_0\subseteq \tilde{C}_1\subseteq \tilde{C}_2\subseteq\cdots \subseteq \tilde{C}_T.
	\end{equation}
	
	The $k_1=1$ slicing of the $\tilde{C}_i$ forms an increasing sequence of closed subsets of the compact set $\tilde{C}_{\textnormal{SV}}\cap \{k_1=1\}$, hence the limit exists. Let $\tilde{C}$ denote the convex cone of the closure of this limit. Since the nef condition is a closed condition, it follows that $\tilde{C}\subseteq \tilde{C}_{T}\subseteq \tilde{C}_{\textnormal{SV}}$. Moreover, $\tilde{C}$ has the property that for any tuple $(k_1,\dots,k_t)\in \tilde{C}$, the new tuple given by taking average of adjacent entries also lies in $\tilde{C}$.
	
	\medskip
	\noindent\textbf{Claim: }$(1,\dots,1)\in\tilde{C}$.
	
	\noindent\emph{Proof: }Let $\tilde{D}$ denote the $k_1=1$ slicing of $\tilde{C}$. The continuous function
	\begin{equation}
		f(1,k_2,\dots,k_t)=1-k_t.
	\end{equation}
	attains its minimal value $f_{\textnormal{min}}$ on the compact set $\tilde{D}$. Let $$\tilde{E}=\{(1,\alpha_2,\dots,\alpha_t)\in\tilde{D}\big|\ f(1,\alpha_2,\dots,\alpha_t)=f_{\textnormal{min}}>0 \}.$$ 
	This is another compact set. Assume $(1,\dots,1)\notin \tilde{D}$. We define another function only on $\tilde{E}$
	\begin{equation}
		g(1,\alpha_2,\dots,\alpha_t)=j,\quad\textnormal{if $\alpha_2=\cdots=\alpha_{j-1}=1$ but $\alpha_j<1$}.
	\end{equation}
	This bounded function attains its minimal value $s$ at some point $(1,\dots,1,\beta_s,\dots,\beta_t)$ in $\tilde{E}$. However, 
	\begin{itemize}
		\item If $2<s<t$, then by Lemma \ref{Improvement of nef estimation}, 
		\begin{equation}
			(1,\dots,\frac{1+\beta_s}{2},\frac{1+\beta_s}{2},\beta_{s+1},\dots,\beta_t)\in \tilde{D}.
		\end{equation}
		This tuple also lies in $\tilde{E}$, but $g$ has values $s-1$ on it, which contradicts the assumption on the minimality of $s$.
		\item  If $s=t$, again by Lemma \ref{Improvement of nef estimation},
		\begin{equation}
			(1,\dots,\frac{1+\beta_t}{2})\in \tilde{D}.
		\end{equation}
		However, $f$ takes the value $\frac{1-k_t}{2}<1-k_t$ on this tuple, contradicting the minimality of $\tilde{E}$.
		\item If $s=2$, again by Lemma \ref{Improvement of nef estimation},
		\begin{equation}
			(\frac{1+\beta_2}{2},\frac{1+\beta_2}{2},\dots,\beta_t)\in \tilde{C}.
		\end{equation}
		Since $\tilde{C}$ is a cone, dividing by $\frac{1+\beta_2}{2}$ yields
		\begin{equation}
			(1,1,\dots,\frac{2\beta_n}{1+\beta_2})\in \tilde{D}.
		\end{equation}
		However, $f$ takes the value $1-\frac{\beta_t}{1+\beta_2}<1-\beta_t$, contradicting the minimality of $\tilde{E}$.
	\end{itemize}
	
	We conclude that $(1,\dots,1)\in \tilde{C}_{\textnormal{nef}}$. In other words, we have shown that the weight
	\begin{equation}
		(1,p^{a(i_1)},\cdots,p^{a(i_1)+a(i_2)+\cdots+a(i_{N-1})})
	\end{equation}
	is nef. Since we could present the argument in this subsection to all boundary weights of $\tilde{C}_{\textnormal{SV}}$, we have shown that $\tilde{C}_{T}=\tilde{C}_{\textnormal{SV}}$, and hence $C_T=C_{\textnormal{SV}}$. Consequently, this establishes the sufficiency part of the nefness criterion for $\mathcal{L}_X(\lambda)$. $\hfill\square$

	\medskip
	\subsection{ Step 4: The sufficiency part of the nefness of $\mathcal{L}_{Y^i_m}(\lambda)$ with $i\in T$.}
	
	For simplicity, we assume that $i=i_1=1\in T$ and $m<m_{1}$. The proof is again based on analyzing the intersection number
	\begin{equation}
		\big([\mathcal{L}_{Y^1_m}(\lambda)]\cdot C\big)
	\end{equation}
	for any complete curve $C\subseteq Y^1_m$, or its pullback along any proper surjective morphism. 
    The induction hypothesis in subsection 7.1, in particular, implies the following assumption.
	
	\medskip
	\noindent\textbf{Induction Hypothesis: }
	
	$\bullet$ The sufficiency part of the nefness criterion for $\mathcal{L}_{Y^1_l}(\lambda)$, where $1\le l< m$, holds for any proper variety over $Y^1_l$ where $\mathcal{F}^1_l$ has slope $\le\vec{r}_0$. More precisely, if $\pi:W\to Y^1_l$ is a proper morphism such that $\vec{r}(Z,\pi^\ast\mathcal{F}^1_l)\le \vec{r}_0$, then $\pi^\ast\mathcal{L}_{Y^1_l}(\lambda)$ is nef provided that $\lambda$ satisfies the condition given by Theorem \ref{Nefness criterion for minimal partial flag space}.
	
	$\bullet$ The sufficiency part of the nefness criterion for $\mathcal{L}_{Y^1_l}(\lambda)$, where $m\le l\le m_1$, holds for any proper variety over $Y^1_l$ where $\mathcal{F}^1_l$ has slope $<\vec{r}_0$. More precisely, if $\pi:W\to Y^1_l$ is a proper morphism such that $\vec{r}(Z,\pi^\ast\mathcal{F}^1_l)< \vec{r}_0$, then $\pi^\ast\mathcal{L}_{Y^1_l}(\lambda)$ is nef provided that $\lambda$ satisfies the condition given by Theorem \ref{Nefness criterion for minimal partial flag space}.
	
	Under these hypotheses, we prove that for any projective morphism $Z\to Y^1_j$, where $\vec{r}(Z,\pi^\ast\mathcal{F}^1_j)=\vec{r}_0\allowbreak=(r_1,r_2,\dots,r_N)$, the pullback $\pi^\ast\mathcal{L}_{Y^1_j}(\lambda)$ to $Z$ is nef if the weight $\lambda$ satisfies the condition given by Theorem \ref{Nefness criterion for minimal partial flag space}.


    For the proof, we first notice that there are two possibilities of $r_{\textnormal{tot}}$, according to Lemma \ref{Properties of slope}:
    \begin{enumerate}[itemsep=0pt,topsep=2pt, parsep=0pt]
        \item Either $r_{\textnormal{tot}}>n_1+n_2+\cdots +n_N$, or
        \item $r_{\textnormal{tot}}=n_1+\cdots+n_N$.
    \end{enumerate}

	We first consider the case $r_{\textnormal{tot}}>n_1+\cdots+n_N$.
	
	\begin{construction}\label{Construction when flagged at essential set} We construct a variety $\tilde{Z}$ through a sequence of blow-ups via the following steps:

    \vspace{0.25em}
    Step 1. Blow-ups in the direction of $V^{-1}$.
    \begin{enumerate}[label=\Alph*.,
    leftmargin=4em,   
    labelsep=0.6em,   
    itemsep=0pt,      
    topsep=2pt,       
    parsep=0pt        
    ]
        \item Let $m^1_{N+1}(=m^1_1)=m$. For $2\le j\le N$, assuming that $m^1_{j+1}$ is defined, we set $m_j^1:=\allowbreak n-m_{j+1}+m_{j+1}^1-r_{N+1-j}$, where the subscripts are taken modulo $N$. Finally we take $m^0_1:=n-m_2+m^1_2-r_N$.
        \item Let $Z_1$ be the closed subscheme of
	\begin{equation}
		Z\times_X Y^N_{m^1_N}\times_X \cdots\times_X Y^2_{m^1_2}\times_X Y^1_{m^0_1}
	\end{equation}
	defined by the condition
	\begin{equation}
		\mathcal{F}^j_{m^1_j}\subseteq V^{-1}(\mathcal{F}^{j+1,(p)}_{m^1_{j+1}}), \textnormal{ }\forall 1\le j\le N,\quad \textnormal{and } \mathcal{F}^1_{m^0_1}\subseteq V^{-1}(\mathcal{F}^{2,(p)}_{m^1_2}).
	\end{equation}
    (As we have already explained in the beginning of this section, we omit the symbol $\pi^\ast$ from all our notations. The condition is in fact imposed on the pullback of these vector bundles to the fiber product scheme.)
    \end{enumerate}        
	\end{construction} 
	
    Note that if $\mathcal{F}^{j+1}_{m^1_{j+1}}$ is chosen, then $\mathcal{F}^j_{m^1_j}$ will be uniquely determined on a dense open subset, as follows from the definition of slopes. The fiber over the complementary closed subset is isomorphic to a Grassmannian $P'/P''$ for some appropriate parabolic subgroups $P'$ and $P''$ of $G$. Hence, the above construction can be realized as a consecutive sequence of blow-ups.
	
	\vspace{0.25em}
	Step 2. Blow-ups in the direction of $V$.

    \begin{enumerate}[label=\Alph*.,
    leftmargin=4em,   
    labelsep=0.6em,   
    itemsep=0pt,      
    topsep=2pt,       
    parsep=0pt        
    ]
        \item For $l\ge 2$, suppose that for all $1\le s\le l-1$ and $1\le t\le N$, the integers $m^s_t$ are defined and the scheme $Z_{l-1}$ has been constructed, whose moduli problem is given by recording subbundles $\mathcal{F}^t_{m^{s}_t}\subseteq \omega_t$ or rank $m^s_t$, together with certain inclusion relations relevant to the essential Verschiebung map. Consider the slope of $\mathcal{F}^1_{m^{l-1}_1}$ on $Z_{l-1}$. By construction, we have $\vec{r}(Z_{l-1},\mathcal{F}^1_{m^{l-1}_1})\le \vec{r}_0$. 
        \begin{itemize}
            \item If $\vec{r}(Z_{l-1},\mathcal{F}^1_{m^{l-1}_1})<\vec{r}_0$, we terminate the inductive construction and set $\tilde{Z}=Z_{l-1}$;
            \item Otherwise, we first set $m^l_2:=m^{l-1}_1+r_N+m_2-n$, and then set inductively $m^l_{j+1}:=\allowbreak m^l_{j}+r_{N+1-j}+m_{j+1}-n$ for $2\le j\le N$.
        \end{itemize}

        \item Let $Z_l$ be the closed subscheme of 
	\begin{equation}
		Z_{l-1}\times_X Y^1_{m^l_1}\times_X\cdots \times_X Y^N_{m^l_{N}}
	\end{equation}
	defined by the condition
	\begin{equation}
		\mathcal{F}^{j}_{m^l_j}\subseteq V^{-1}(\mathcal{F}^{j+1,(p)}_{m^{l}_{j+1}}),\ \forall 2\le j\le N,\quad \textnormal{and }\mathcal{F}^{1}_{m^{l-1}_1}\subseteq V^{-1}(\mathcal{F}^{2,(p)}_{m^{l}_{2}}).
	\end{equation}
        
    \end{enumerate}
	
	\vspace{0.25em} 

\begin{lem} We have
\begin{equation}
    m^l_1-m^{l-1}_1=r_{\textnormal{tot}}-(n_1+n_2+\cdots+n_N).
\end{equation}
\end{lem}
\emph{Proof.} This follows from our construction that
\begin{equation}
\begin{aligned}
    m^l_1-m^{l-1}_1&=(m^l_1-m^l_{N})+(m^l_N-m^l_{N-1})+\cdots+(m^l_{3}-m^l_2)+(m^l_2-m^l_1)\\
    &=(r_1-n_1)+(r_{2}-n_N)+\cdots+(r_{N-1}-n_3)+(r_N-n_2)\\
    &=r_{\textnormal{tot}}-(n_1+\cdots+n_N). 
\end{aligned}
\end{equation} $\hfill\square$

By the above lemma, one sees that step 2 must terminate at some $l=l_0$, since otherwise we would have
\begin{equation}
    \vec{r}(Z_{l-1},\mathcal{F}^1_{m^{l-1}_1})=\vec{r}_0,
\end{equation}
which forces $m^l_1>m^{l-1}_1$. That is, at the embedding $\tilde\tau_1$, we have subbundles $\mathcal{F}^1_{m^l_1}\subseteq \omega_1$ of arbitrarily large rank, which is impossible.

	When the inductive construction stops, we denote the resulting variety by $\tilde{Z}$. This completes our construction of the auxiliary variety $\tilde{Z}$.
	
	\begin{lem}\label{The auxiliary variety is proper surjective}
		The forgetful map $\tilde{Z}\to Z$ is proper surjective. 
	\end{lem}	
	\noindent\emph{Proof:} It suffices to show that each morphism $Z_l\to Z_{l-1}$ is surjective, which reduces to checking the induced map on closed points is surjective. For $2\le j\le N$, assume we have chosen $\mathcal{F}^j_{m^l_j}$, then by the moduli problem, we need to choose $\mathcal{F}^{j+1}_{m_{j+1}^l}\subseteq\omega_j$ such that
	\begin{equation}
		V\mathcal{F}^{j}_{m^l_j}+\mathcal{F}^{j+1,(p)}_{m^{l-1}_{j+1}}\subseteq \mathcal{F}^{j+1,(p)}_{m^l_{j+1}} \subseteq\omega_j^{(p)}.
	\end{equation}
	The choice of $\mathcal{F}^{j+1,(p)}_{m_{j+1}^l}$ is generically unique, and the fiber over the complementary closed subscheme are Grassmannians $P'^{(p)}/P''^{(p)}$ for some parabolic subgroups. Thus, $Z_l$ is obtained from $Z_{l-1}$ by sequence of consecutive blow-ups and Frobenius twists. $\hfill\square$
	
	\medskip
	The basic idea of Construction \ref{Construction when flagged at essential set} is to "draw parallelograms in the correct direction". We illustrate this construction with the following example.
	
	\begin{example}
		$G=G(U(5,2)\times U(6,1)\times U(4,3))$. Let $Z/Y^1_2$ be a projective variety where $\vec{r}(Z,\mathcal{F}^1_2)=(3,0,4)$. Then $\tilde{Z}$ is obtained by further recording the following vector bundles and their inclusion relations, as indicated by the colored lines.
		
		
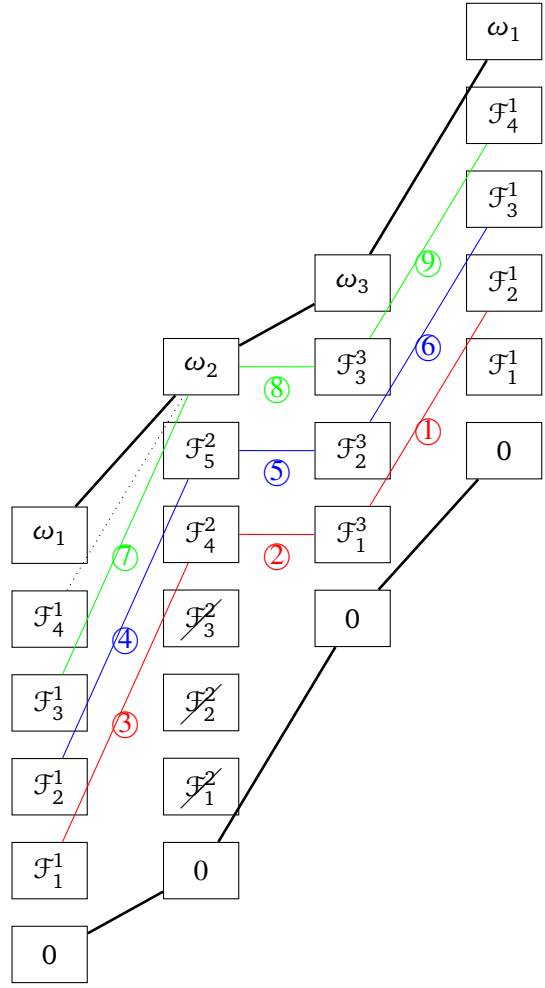
\begin{center}
			\begin{tikzpicture}[node distance=10pt and 1cm, auto]
				
				\begin{scope}
					\node (1) [draw, rectangle, minimum width=1cm, minimum height=0.75cm, align=center] {$\omega_1$};
					\node (2) [draw, rectangle, minimum width=1cm, minimum height=0.75cm, align=center, below=of 1] {${\mathcal{F}^1_4}$};
					\node (3) [draw, rectangle, minimum width=1cm, minimum height=0.75cm, align=center, below=of 2] {$\mathcal{F}^1_3$};
					\node (4) [draw, rectangle, minimum width=1cm, minimum height=0.75cm, align=center, below=of 3] {${\mathcal{F}^1_2}$};
					\node (5) [draw, rectangle, minimum width=1cm, minimum height=0.75cm, align=center, below=of 4] {$\mathcal{F}^1_1$};
					\node (6) [draw, rectangle, minimum width=1cm, minimum height=0.75cm, align=center, below=of 5] {$0$};
					
					\node (7) [draw, rectangle, minimum width=1cm, minimum height=0.75cm, align=center, right=of 5] {0};
					\node (8) [draw, rectangle, minimum width=1cm, minimum height=0.75cm, align=center, above=of 7] {$\cancel{\mathcal{F}^2_1}$};
					\node (9) [draw, rectangle, minimum width=1cm, minimum height=0.75cm, align=center, above=of 8] {$\cancel{\mathcal{F}^2_2}$};
					\node (10) [draw, rectangle, minimum width=1cm, minimum height=0.75cm, align=center, above=of 9] {$\cancel{\mathcal{F}^2_3}$};
					\node (11) [draw, rectangle, minimum width=1cm, minimum height=0.75cm, align=center, above=of 10] {$\mathcal{F}^2_4$};	
					\node (12) [draw, rectangle, minimum width=1cm, minimum height=0.75cm, align=center, above=of 11] {$\mathcal{F}^2_5$};
					\node (13) [draw, rectangle, minimum width=1cm, minimum height=0.75cm, align=center, above=of 12] {$\omega_2$};
					
					\node (14) [draw, rectangle, minimum width=1cm, minimum height=0.75cm, align=center, right=of 10] {$0$};
					\node (15) [draw, rectangle, minimum width=1cm, minimum height=0.75cm, align=center, above=of 14] {$\mathcal{F}^3_1$};
					\node (16) [draw, rectangle, minimum width=1cm, minimum height=0.75cm, align=center, above=of 15] {$\mathcal{F}^3_2$};
					\node (17) [draw, rectangle, minimum width=1cm, minimum height=0.75cm, align=center, above=of 16] {$\mathcal{F}^3_3$};
					\node (18) [draw, rectangle, minimum width=1cm, minimum height=0.75cm, align=center, above=of 17] {$\omega_3$};
					
					\node (19) [draw, rectangle, minimum width=1cm, minimum height=0.75cm, align=center, right=of 16] {$0$};
					\node (20) [draw, rectangle, minimum width=1cm, minimum height=0.75cm, align=center, above=of 19] {$\mathcal{F}^1_1$};
					\node (21) [draw, rectangle, minimum width=1cm, minimum height=0.75cm, align=center, above=of 20] {$\mathcal{F}^1_2$};
					\node (22) [draw, rectangle, minimum width=1cm, minimum height=0.75cm, align=center, above=of 21] {$\mathcal{F}^1_3$};
					\node (23) [draw, rectangle, minimum width=1cm, minimum height=0.75cm, align=center, above=of 22] {$\mathcal{F}^1_4$};
					\node (24) [draw, rectangle, minimum width=1cm, minimum height=0.75cm, align=center, above=of 23] {$\omega_1$};

					\draw [line width=1pt] (1)--(13);
					\draw [line width=1pt] (13)--(18);
					\draw [line width=1pt] (18)--(24);
					\draw [line width=1pt] (6)--(7);
					\draw [line width=1pt] (7)--(14);
					\draw [line width=1pt] (14)--(19);

					\draw [red] (21)--(15) node[midway, below] {$\textcircled{1}$};
					\draw [red] (15)--(11) node[midway, below] {$\textcircled{2}$};
					\draw [red] (11)--(5) node[midway, below] {$\textcircled{3}$};
					\draw [blue] (4)--(12) node[midway, below] {$\textcircled{4}$};
					\draw [blue] (12)--(16) node[midway, below] {$\textcircled{5}$};
					\draw [blue] (16)--(22) node[midway, below] {$\textcircled{6}$};
					\draw [green] (3)--(13) node[midway, below] {$\textcircled{7}$};
					\draw [green] (13)--(17) node[midway, below] {$\textcircled{8}$};
					\draw [green] (17)--(23) node[midway, below] {$\textcircled{9}$};
					
					\draw [dotted] (2)--(13);
					
				\end{scope}
			\end{tikzpicture}
		\end{center}
		
		$\bullet$ The numbers on the arrows represent the order of construction. The red lines correspond to the first step, in which we start from $\mathcal{F}^1_2$ on the right and draw lines "from right to left"; The blue and green lines correspond to the second step, where we take blow-ups up to Frobenius twists "from left to right".
		
		$\bullet$ We have $V^{-1}(\omega_2^{(p)})\cap \mathcal{F}_4^1=\mathcal{F}^1_4$, and the intersection is a vector bundle of rank $4$, indicated as the dotted arrow in the picture. Thus, $\vec{r}(\tilde{Z},\mathcal{F}^1_4)=(3,0,3)<\vec{r}$, and we end our construction after we record this vector bundle in our moduli problem.
		
		$\bullet$ An alternative way to understand the terminal of our construction is that "we cannot draw a line segment starting at $\mathcal{F}^1_4$ which is parallel to \textcircled{7}". 
	\end{example}
	
	Next, we associate a "strata Hasse invariant" with each "parallelogram" in the picture of $\tilde{Z}$. Assume that $\tilde{Z}=Z_{l_0}$ in our construction and let $\mathcal{F}^1_{m^{l_0}_1}$ be the vector bundle of the highest rank at the first place that appears in the construction (which could be $\omega_1$, i.e., $m^{l_0}_1=m_1$). For $2\le l\le l_0-1$ we have the rank condition by construction:
	\begin{equation}
		\begin{aligned}
			\textnormal{rank $\mathcal{F}^N_{m_N^{l+1}}/\mathcal{F}^{N}_{m_N^{l}}$}&=(m^{l+1}_{N-1}+r_{2}+m_N-n)-(m^{l}_{N-1}+r_{2}+m_N-n)=m^{l+1}_{N-1}-m^{l}_{N-1}\\
			&=\textnormal{rank $\mathcal{F}^{N-1}_{m_{N-1}^{l+1}}/\mathcal{F}^{N-1}_{m_{N-1}^{l}}$}=\cdots=\textnormal{rank $\mathcal{F}^2_{m_2^{l+1}}/\mathcal{F}^{2}_{m_2^{l}}$}=\textnormal{rank $\mathcal{F}^1_{m^{l}_1}/\mathcal{F}^1_{m^{l-1}_1}$}.
		\end{aligned}
	\end{equation}
	The equality for $l=1$ can be checked in a similar way. Thus, we obtain a section on this variety by taking $h^l_j$ as the determinant of the following morphisms (The first morphism is induced by inclusion, and the second isomorphism is induced by the essential Verschiebung map)
	\begin{equation}
		\begin{aligned}
			&\tilde{h}^l_j:\mathcal{F}^{j}_{m^{l+1}_{j}}/\mathcal{F}^{j}_{m^{l}_{j}}\stackrel{}{\longrightarrow} V^{-1}(\mathcal{F}^{j+1,(p)}_{m^{l+1}_{j+1}})/V^{-1}(\mathcal{F}^{j+1,(p)}_{m^{l}_{j+1}}) \stackrel{\simeq}{\longrightarrow} \big(\mathcal{F}^{j+1}_{m^{l+1}_{j+1}}/\mathcal{F}^{j+1}_{m^{l}_{j+1}}\big)^{(p)},\quad \textnormal{if $2\le j\le N$},\\
			&\tilde{h}^{l}_1:\mathcal{F}^{1}_{m^l_{1}}/\mathcal{F}^{1}_{m^{l-1}_{1}}\stackrel{}{\longrightarrow} V^{-1}(\mathcal{F}^{2,(p)}_{m^{l+1}_{2}})/V^{-1}(\mathcal{F}^{2,(p)}_{m^{l}_{2}}) \stackrel{\simeq}{\longrightarrow} \big(\mathcal{F}^{2}_{m^{l+1}_{2}}/\mathcal{F}^{2}_{m^{l}_{2}}\big)^{(p)},\quad \textnormal{if $j=1$}.
		\end{aligned}
	\end{equation}
	Therefore, in $(\textnormal{Pic $\tilde{Z}$})_\mathbb{Q}$, we have
	\begin{equation}
		\begin{aligned}
			&c_1(h^l_j)=pc_1\big(\mathcal{F}^{j+1}_{m^{l+1}_{j+1}}/\mathcal{F}^{j+1}_{m^l_{j+1}})\big)-c_1\big(\mathcal{F}^j_{m^{l+1}_j}/\mathcal{F}^j_{m^l_j}\big),\quad \textnormal{if $2\le j\le N$},\\
			&c_1(h^l_1)=pc_1\big(\mathcal{F}^{2}_{m^{l+1}_{2}}/\mathcal{F}^{2}_{m^l_{2}}\big)-c_1\big(\mathcal{F}^1_{m^l_1}/\mathcal{F}^1_{m^{l-1}_1}\big),\quad \textnormal{if $j=1$}.
		\end{aligned}
	\end{equation}
	"The $l$-th" partial Hasse invariant is defined as
	\begin{equation}
		\begin{aligned}
			[h^l]:&=p^{N-1}[h^l_N]+ p^{N-2}[h^l_{N-1}]+\cdots+ p[h^l_2]+ [h^l_1]\\
			&=p^N[\mathcal{F}^1_{m^{l+1}_1}/\mathcal{F}^1_{m^l_1}]-[\mathcal{F}^1_{m^l_1}/\mathcal{F}^1_{m^{l-1}_1}]\in (\textnormal{Pic $\tilde{Z}$})_\mathbb{Q}.
		\end{aligned}
	\end{equation}
	
	\begin{rmk}
		If $r_j=0$ for some $j\ge2$, then rank $V_{es}^{-1}({\mathcal{F}^{j+1,(p)}_{m^l_{j+1}}})=n_{j+1}+m^l_{j+1}=m^l_{j}=$ rank $\mathcal{F}^j_{m^l_j}$, and the inclusion condition is actually an equality $\mathcal{F}^j_{m^l_j}=V_{es}^{-1}(\mathcal{F}^{j+1,(p)}_{m^l_{j+1}})$. In this case, the corresponding section $\tilde{h}^l_j$ is an isomorphism. The case for $j=1$ is similar. Therefore, either $\tilde{h}^l_j$ is a section of a line bundle whose zero section has codimension 1, or it is a nonzero section of the trivial bundle. We refer to both of them as "strata partial Hasse invariants".
	\end{rmk}

    \begin{lem}
        For $0\le l\le l_0-1$, we have $\vec{r}(\tilde{Z},\mathcal{F}^1_{m^l_1})\le \vec{r}_0$. Besides, $\vec{r}(\tilde{Z},\mathcal{F}^1_{m^{l_0}_1})<\vec{r}_0$.
    \end{lem}
    \emph{Proof.} This follows directly from the construction. $\hfill\square$
    
	\begin{lem} For any $1\le l\le l_0-1$ and $1\le j\le N$, let $Z_j^l$ be the zero locus of $h^l_j$ (if it is non-empty). Then $\vec{r}(Z^l_j,\mathcal{F}^1_{m^l_1})<\vec{r}_0$.
	\end{lem}
	\noindent\emph{Proof: }Assume $2\le j\le N$ for simplicity, the argument for $j=1$ is similar. By definition, we have
	\begin{equation}
		\textnormal{rank }\big(\mathcal{F}^{j}_{m^{l+1}_j}/\mathcal{F}^j_{m^l_j}\cap V^{-1}(\mathcal{F}^{j+1,(p)}_{m^{l}_{j+1}})/\mathcal{F}^j_{m^l_j}\big)\ge 1
	\end{equation}
	when viewed as subbundles of $V_{es}^{-1}(\mathcal{F}^{j+1,(p)}_{m^{l+1}_{j+1}})$. Consequently,
	\begin{equation}
		\textnormal{dim }\big(\omega_j\cap V^{-1}(\mathcal{F}^{j+1,(p)}_{m^{l}_{j+1}})\big)\ge \textnormal{dim }\mathcal{F}^{j}_{m^l_j}+1=m^l_j+1>n+m^l_{j+1}-m_{j+1}-r_j.
	\end{equation}
	Therefore, the slope of $\mathcal{F}^1_{m^l_1}$ on $Z^l_j$ is at most $\vec{r}'=(r_N,\dots, r_j-1,\ast,\dots,\ast)<\vec{r}$.$\hfill\square$

    \medskip
	We are now ready to prove that the pullback of $\mathcal{L}_{Y^1_m}(\lambda)$ to $Z$ is nef, under the induction hypothesis stated at the beginning of this subsection. Note that it suffices to prove
    \begin{equation}
        \big(\mathcal{L}_{Y^1_m}(\{k_i\};\alpha)\cdot C\big)\ge0
    \end{equation}
    for any complete curve $C\subseteq \tilde{Z}$.
	
	Let $C$ be a complete curve on $\tilde{Z}$. We notice the following:
    \begin{itemize}[itemsep=0pt,topsep=2pt, parsep=0pt]
        \item If $C$ is contained in the zero locus of $h^s_j$ for some $1\le j\le N$, let $1\le s\le l_0-1$ be the smallest superscript of such $h^s_j$'s. By the above lemma, $\vec{r}(Z(h^s_j),\mathcal{F}^1_{m^s_1})<\vec{r}_0$. It follows from the induction hypothesis that the pullback of $[\mathcal{L}_{Y^1_{m^s_1}}(k_i;\alpha)]$ to $Z(h^s_j)$ is nef, as long as the coefficients satisfy Theorem \ref{Nefness criterion for minimal partial flag space}; 
        \item If $C$ is not contained in the zero locus of any $h^s_j$, we set $s=l_0$. By construction, $\vec{r}(\tilde{Z},\mathcal{F}^{1}_{m^{l_0}_1})<\vec{r}_0$. By induction, $\mathcal{L}_{Y^1_{m^s_1}}(k_i;\alpha)$ is nef as long as the coefficients satisfy Theorem \ref{Nefness criterion for minimal partial flag space}. 
        \item Since we always have $\vec{r}(\tilde{Z},\mathcal{F}^1_{m^0_1})\le\vec{r}$, the pullback of $[\mathcal{L}_{Y^1_{m^0_1}}(k_i;\alpha)]$ is also nef as long as the coefficients satisfy Theorem \ref{Nefness criterion for minimal partial flag space}.
    \end{itemize}
    
	If $s=1$, then $h^1_j$ vanishes for some $1\le j\le N$, and $C$ is contained in a stratum where the slope of $\mathcal{F}^1_m$ is strictly smaller than $\vec{r}$. By the inductive hypothesis, we obtain the non-negativity of the intersection numbers.
	
	If $s\ge2$, consider the positive linear combination
	\begin{equation}
		\begin{aligned}
			[h]&=[h^{s-1}]+\frac{p^N+1}{p^N}[h^{s-2}]+\frac{p^{2N}+p^N+1}{p^{2N}}[h^{s-3}]+\cdots+\frac{p^{(s-2)N}+\cdots+1}{p^{(s-2)N}}[h^1]\\
			&=p^N[\mathcal{F}^1_{m^s_1}]-\frac{p^{sN}-1}{p^{(s-2)N}(p^N-1)}[\mathcal{F}^1_{m^1_1}]+\frac{p^{(s-1)N}-1}{p^{(s-2)N}(p^N-1)}[\mathcal{F}^1_{m^0_1}].
		\end{aligned}
	\end{equation}
	We use $[h]$ to kill the coefficients of $[\mathcal{F}^1_m]=[\mathcal{F}^1_{m^1_1}]$ and get
	\begin{equation}
		\begin{aligned}
			&\quad\ [\mathcal{L}_{Y^1_m}(k_i;\alpha)]\\
            &=k_1[\omega_1/\mathcal{F}^1_m]+\alpha[\mathcal{F}^1_m]+\sum_{l=2}^N k_l[\omega_l]\\
			&=\frac{p^{(s-2)N}(p^N-1)}{p^{sN}-1}(k_1-\alpha)[h]-\frac{p^{(s-1)N}(p^N-1)}{p^{sN}-1}(k_1-\alpha)[\mathcal{F}^1_{m^s_1}] -\frac{p^{(s-1)N}-1}{p^{sN}-1}(k_1-\alpha)[\mathcal{F}^1_{m^0_1}]+\sum_{l=1}^N k_l[\omega_l]\\
			&=\beta[h]+ \bigg(tk_1[\omega_1/\mathcal{F}^1_{m^s_1}]+ \big(tk_1-\frac{p^{(s-1)N}(p^N-1)}{p^{sN}-1}(k_1-\alpha)\big)[\mathcal{F}^1_{m^s_1}] +tk_2[\omega_2]+\cdots+tk_N[\omega_N]\bigg)\\
			&\ +\bigg((1-t)k_1[\omega_1/\mathcal{F}^1_{m^0_1}]+ \big((1-t)k_1-\frac{p^{(s-1)N}-1}{p^{sN}-1}(k_1-\alpha)\big)[\mathcal{F}^1_{m^0_1}] +(1-t)k_2[\omega_2]+\cdots+(1-t)k_N[\omega_N]\bigg)\\
			&=\beta[h]+[\mathcal{L}_{Y^1_{m^s_1}}(\{tk_i\};\alpha^s_1)]+[\mathcal{L}_{Y^1_{m^0_1}}(\{(1-t)k_i\};\alpha^0_1)].
		\end{aligned}
	\end{equation}
	In the last two equations, we wish to express the original class $[\mathcal{L}_{Y^1_{m}}]$ as a linear combination of Hasse invariants and pullbacks of the classes $[\mathcal{L}_{Y^1_{m^s_1}}]$ and $[\mathcal{L}_{Y^1_{m^0_1}}]$. To simplify the notation, we use $\beta,\alpha^s_1,\alpha^0_1$ to represent their corresponding coefficients. 

    \begin{prop}
        There exists a suitable choice of real numbers $\beta\ge0$ and $t\in [0,1]$ fitting into the above equation, such that the weights for the line bundles $[\mathcal{L}_{Y^1_{m^s_1}}(\{tk_i\};\alpha^s_1)]$ and $[\mathcal{L}_{Y^1_{m^0_1}}(\{(1-t)k_i\};\alpha^0_1)]$ satisfy their respective nefness conditions given by Theorem \ref{Nefness criterion for minimal partial flag space}.
    \end{prop}
	\noindent\emph{Proof.} By construction, 
    \begin{equation}
        \beta=\frac{p^{(s-2)N}(p^N-1)}{p^{sN}-1}(k_1-\alpha)\ge0.
    \end{equation}
    It remains to prove the statement for $t$. The proof is based on a case-by-case analysis.
    
	\noindent\emph{Case 1:} If $\#T\ge 2$, then $\alpha$ satisfies
	\begin{equation}
		k_1\ge\alpha,\qquad  p^{a(1)}\alpha\ge k_2,
	\end{equation}
	and we must choose $t$ such that (the remaining nefness conditions can be verified easily)
	\begin{equation}
		\begin{aligned}
			p^{a(1)}\big(tk_1-\frac{p^{(s-1)N}(p^N-1)}{p^{sN}-1}(k_1-\alpha)\big)&\ge tk_2,\\
			p^{a(1)}\big((1-t)k_1-\frac{p^{(s-1)N}-1}{p^{sN}-1}(k_1-\alpha)\big)&\ge (1-t)k_2.
		\end{aligned}
	\end{equation}
	We could take the sum of the inequalities and obtain
	\begin{equation}
		p^{a(1)}(k_1-(k_1-\alpha))=p^{a(1)}\alpha\ge k_2.
	\end{equation}
	This is valid because of the nefness condition for $\mathcal{L}_{Y^1_m}(\{k_i\};\alpha)$. Thus, there exists some $t$ that satisfies both inequalities. We claim that such $t$ must lie in $[0,1]$. We could rearrange the inequalities as follows.
    \begin{equation}
    \begin{aligned}
        &(p^{a(1)}k_1-k_2)t\ge p^{a(1)}\frac{p^{(s-1)N}(p^N-1)}{p^{sN}-1}(k_1-\alpha),\\
        &(p^{a(1)}k_1-k_2)(1-t)\ge p^{a(1)}\frac{p^{(s-1)N}-1}{p^{sN}-1}(k_1-\alpha).
    \end{aligned}
    \end{equation}
	The conditions force $t\in [0,1]$.

	\medskip
	\noindent\emph{Case 2:} If $\# T=1$ but $m_1\notin\{1,n-1\}$, then $a(1)$ is $N$ and $\alpha$ satisfies
	\begin{equation}
		k_1\ge \alpha,\qquad p^N\alpha\ge k_1,
	\end{equation}
	and we must choose $t$ such that
	\begin{equation}
		\begin{aligned}
			p^N\big(tk_1-\frac{p^{(s-1)N}(p^N-1)}{p^{sN}-1}(k_1-\alpha)\big)&\ge tk_1,\\
			p^N\big((1-t)k_1-\frac{p^{(s-1)N}-1}{p^{sN}-1}(k_1-\alpha)\big)&\ge (1-t)k_1.
		\end{aligned}
	\end{equation}
	Similarly, taking the sum of the inequalities, we obtain
	\begin{equation}
		p^N\alpha\ge k_1,
	\end{equation}
	which follows from our assumption. By the same reason as the previous case, there exists some $t\in[0,1]$ that satisfies both inequalities.

	\medskip
	\noindent\emph{Case 3:} If $\# T=1$ and $m_1=n-1$, then $a(1)=N$ and $\alpha$ satisfies
	\begin{equation}
		k_1\ge \alpha,\qquad (p^{N(n-m)}-1)\alpha\ge (p^{N(n-m-1)}-1)k_1.
	\end{equation}
	Let $\delta=r_{\textnormal{tot}}-n_1>0$. By construction, we have $m^0_1=m-\delta$ and $m^s_1=m+(s-1)\delta$. We look for some $t\in [0,1]$ such that
	\begin{equation}
		\begin{aligned}
			(p^{N(n-(m+(s-1)\delta))}-1)\big(tk_1-\frac{p^{(s-1)N}(p^N-1)}{p^{sN}-1}(k_1-\alpha)\big)\ge (p^{N(n-(m+(s-1)\delta)-1)}-1)tk_1,\\
			(p^{N(n-(m-\delta))}-1)\big((1-t)k_1-\frac{p^{(s-1)N}-1}{p^{sN}-1}(k_1-\alpha)\big)\ge (p^{N(n-(m-\delta)-1)}-1)(1-t)k_1.
		\end{aligned}
	\end{equation}
	This follows from a direct but complicated computation, which we present below for the integrality of the proof. The above inequalities are equivalent to 
	\begin{equation}
		\begin{aligned}
			&\mathbb{I}:(p^{N(n-(m+(s-1)\delta))}-p^{N(n-(m+(s-1)\delta)-1)})tk_1\ge  (p^{N(n-(m+(s-1)\delta))}-1) \frac{p^{(s-1)N}(p^N-1)}{p^{sN}-1}(k_1-\alpha),\\
			&\mathbb{II}:(p^{N(n-(m-\delta))}-p^{N(n-(m-\delta)-1)})(1-t)k_1\ge (p^{N(n-(m-\delta))}-1)\frac{p^{(s-1)N}-1}{p^{sN}-1}(k_1-\alpha).
		\end{aligned}
	\end{equation}
	Taking $p^{Ns\delta}\mathbb{I}+\mathbb{II}$, we want to show
	\begin{equation}
		(p^{N(n-(m-\delta))}-p^{N(n-(m-\delta)-1)})k_1\ge \frac{(p^{N(n-(m-\delta))}-p^{Ns\delta})(p^N-1)p^{N(s-1)}+(p^{N(n-(m-\delta))}-1)(p^{N(s-1)}-1)}{p^{sN}-1}(k_1-\alpha).
	\end{equation}
	We could rearrange the terms and denote $q=p^N$ for simplicity. We could simplify the inequalities to be
	\begin{equation}
		\mu\alpha\ge \nu k_1,
	\end{equation}
	where the coefficients for $\alpha,k_1$ respectively are
	\begin{equation}
		\begin{aligned}
			&\mu=(q^{n-m+\delta+s}-q^{s\delta+s}+q^{s\delta+s-1}-q^{n-m+\delta}-q^{s-1}+1),\\
			&\nu=(q^{n-m+s+\delta-1}-q^{n-m+\delta-1}-q^{s\delta+s}+q^{s\delta+s-1}-q^{s-1}+1).
		\end{aligned}
	\end{equation}
	We hope that this is automatically satisfied, i.e., we need
	\begin{equation}
		\frac{\mu}{\nu}\ge\frac{q^{n-m}-1}{q^{n-m-1}-1}.
	\end{equation}
	One could expand both sides, and a direct but boring computation verifies the desired inequality (We need to use the fact $\delta\ge 1$ and $s\ge 2$). 

	\medskip
	\noindent\emph{Case 3':} If $\# T=1$ and $m_1=1$, the proof for this case is essentially the same as above.  $\hfill\square$
	
	\begin{rmk}
		In Construction \ref{Construction when flagged at essential set}, it is possible that $m^{l_0}_1=m_1$ (resp. $m^0_1=0$). In this case the bundle $\mathcal{F}^1_{m^{l_0}_1}$ (resp. $\mathcal{F}^1_{m^0_1}$) is simply $\omega_1$ (resp. $0$). The corresponding line bundle $[\mathcal{L}_{Y^1_{m^0_1}}(\{k_i\};\alpha)]$ (resp. $\mathcal{L}_{Y^1_{m^0_1}}(\{k_i\};\alpha)$) degenerates as the pullback of some automorphic line bundle $[\mathcal{L}_X(\{k_i\})]$ from the underlying Shimura variety. It can be verified (by a case-by-case discussion) that the above computations is still valid in these boundary cases.
        
	\end{rmk}	
	
	Return to our proof. For this complete curve $C$, we have
    \begin{equation}
        \begin{aligned}
            &\quad\big([\mathcal{L}_{Y^1_m}(\{k_i\};\alpha)]\cdot C\big)\\
            &=\big(\beta[h]+[\mathcal{L}_{Y^1_{m^s_1}}(\{tk_i\};\alpha^s_1)]+[\mathcal{L}_{Y^1_{m^0_1}}(\{(1-t)k_i\};\alpha^0_1)]\cdot C\big)\\
            &=\beta\big([h]\cdot C\big)+\big([\mathcal{L}_{Y^1_{m^s_1}}(\{tk_i\};\alpha^s_1)]\cdot C\big)+\big([\mathcal{L}_{Y^1_{m^0_1}}(\{(1-t)k_i\};\alpha^0_1)]\cdot C\big)\ge0.
        \end{aligned}
    \end{equation}
	Here, nonnegativity follows from the choice $\beta\ge0$ and induction. This shows that in the first case (i.e., in the case $r_{\textnormal{tot}}(Z,\mathcal{F}^1_{m})>n_1+n_2+\cdots+n_N$), the bundle $\mathcal{L}_{Y^1_m}(\{k_i\};\alpha)$ is nef over $Z$.

	\medskip
	It remains to deal with the latter case, that is, if
    $$
    r_{\textnormal{tot}}(Z,\mathcal{F}^1_m)=\sum_{j=1}^N r_j=n_1+\dots+n_N.
    $$
    This is the induction base of our argument in this section. 
	
    If the essential degree $t=1$, this is the induction base of our proof in this subsection, and we have shown that the restriction of $\mathcal{L}_{Y^1_m}(k_i;\alpha)$ to $Z$ is nef in Step 2.

    Otherwise, we may assume the essential degree $t\ge 2$. Let $\tilde{Z}$ 
	be the closed subscheme of
	\begin{equation}
		Z\times_X Y^N_{m^1_N}\times_X\cdots\times_X Y^2_{m^1_2}
	\end{equation}
	defined by the condition
	\begin{equation}
		\mathcal{F}^j_{m^1_j}\subseteq V^{-1}_{es}(\mathcal{F}^{j+1,(p)}_{m^1_{j+1}}),\ \forall 1\le j\le N.
	\end{equation}
	Here the integers $m^1_j$ are given in Construction \ref{Construction when flagged at essential set}. Since the total slope takes the minimal possible value, we must have $m^0_1=m^1_1$ and $\mathcal{F}^1_{m^0_1}=\mathcal{F}^1_m\subseteq V^{-1}_{es}(\mathcal{F}^{2,(p)}_{m^1_{2}})$. The above condition implies that the sequence
	\begin{equation}
		\{\mathcal{F}^1_{m}=\mathcal{F}^1_{m^0_1},\mathcal{F}^2_{m^1_2},\mathcal{F}^3_{m^1_3},\dots,\mathcal{F}^{N-1}_{m^{1}_{N-1}},\mathcal{F}^N_{m^1_N},\mathcal{F}^1_{m^0_1}\}
	\end{equation}
	forms an $F_{es},V_{es}$-chain over $\tilde{Z}$. By Theorem \ref{Construction via F_{es},V_{es}-chain}, $\tilde{Z}$ admits a proper surjective morphism to an auxiliary unitary Shimura variety $X'$.
	\begin{lem}
		$X'$ has essential degree $t'=\# T_{X'}\le t=\# T_X$.
	\end{lem}
	\noindent\emph{Proof}: It suffices to show that the essential degree does not increase. By Theorem \ref{Construction via F,V-chain}, $X'$ has signature $\{(m_i',n_i')\}$, where
	\begin{equation}
		m'_i=n-r_{i-1}.
	\end{equation}
	By the definition of slopes, if $m_{i}\in \{0,n\}$, then $r_{i}=0$. Hence
	\begin{equation}
		T_{X'}\subseteq\big\{i|\ m_{i-1}\notin\{0,n\}\big\}.
	\end{equation}
	Thus, $t'\le t$.   $\hfill\square$
	
	Recall the inductive hypothesis at the beginning of this subsection: To prove the nefness criterion for $\mathcal{L}_{Y^1_m}(\{k_i\};\alpha)$ of essential degree $t$, we need to assume that the nefness criterion for $\mathcal{L}_{X'}(\lambda)$ is known for any unitary Shimura variety $X'$ where $t'\le t$. Thus, the line bundle
	\begin{equation}
		\mathcal{L}_{X'}(\{k'_i\})
	\end{equation}
	on $X'$ is nef if and only if the weight $\{k_i\}$ satisfy the condition in Theorem \ref{Nefness criterion for Shimura variety}. In particular,
	\begin{equation}
		\begin{aligned}
			[\mathcal{L}_{X'}(p^{N-1},1,p,\dots,p^{N-2})]&=\sum_{j=2}^{N+1}p^{j-2}[\omega_j']=\sum_{j=1}^Np^{j-2}\Big(\frac{1}{p}[\frac{V\mathcal{F}^{j-1}_{m^1_{j-1}}}{p\mathcal{F}^{j,(p)}_{m^1_j}}]\Big)\\
			&=\sum_{j=2}^{N+1} p^{j-2}\frac{1}{p}\Big([\frac{V\mathcal{F}^{j-1}_{m^1_{j-1}}}{Vp\H(A/X)_{j-1}}]+[\frac{Vp\H(A/X)_{j-1}}{p\mathcal{F}^{j,(p)}_{m^1_j}}]\Big)\\
			&=\sum_{j=2}^{N+1}p^{j-2}\Big([\omega_j/\mathcal{F}^j_{m^1_j}]+\frac{1}{p}[\mathcal{F}^{j-1}_{m^1_{j-1}}]\Big)\\
			&=p^{N-1}[\omega_1/\mathcal{F}^1_m]+\frac{1}{p}[\mathcal{F}^1_m]+[\omega_2]+p[\omega_3]+\cdots+p^{N-2}[\omega_N]
		\end{aligned}
	\end{equation}
	is nef. Since in the previous parts of this paper, we have presented many times how to compute the classes $[\omega_j']$ from the moduli problem, we omit the details here. We remark that some of the $\omega'_j$'s might have rank $0$ or $n$, so their classes $[\omega_j]\in (\textnormal{Pic X'})_{\mathbb Q}$ are 0. It is harmless to include them to maintain the integrality of the formula.     
	We have
	\begin{equation}
		\begin{aligned}
			[\mathcal{L}_{Y^1_m}(\{k_i\};\alpha)]=\frac{p}{p^N-1}(k_1-\alpha)[\mathcal{L}_{X'}(p^{N-1},1,\dots,p^{N-2})]+[\mathcal{L}_X(\lambda')]
		\end{aligned}
	\end{equation}
	on $\tilde{Z}$. Here $\lambda'=\{k_i'\}$, where
	\begin{equation}
		\begin{aligned}
			&k_1'=k_1-\frac{p\cdot p^{N-1}}{p^N-1}(k_1-\alpha)=\frac{1}{p^N-1}(p^N\alpha-k_1),\\
			&k_2'=k_2-\frac{p}{p^N-1}(k_1-\alpha),\\
			&k_3'=k_3-\frac{p\cdot p}{p^N-1}(k_1-\alpha),\\
			&\cdots\cdots\\
			&k_N'=k_N-\frac{p\cdot p^{N-2}}{p^N-1}(k_1-\alpha).
		\end{aligned}
	\end{equation}
	It is straightforward to verify all the inequalities $p^{a(i_s)}k_{i_s}'\ge k_{i_{s+1}}'$ for all $2\le s\le t$. For $s=1$, we have (recall that $i_1=1$)
	\begin{equation}
		p^{a(1)}k_1'-k_{i_2}'=p^{a(1)}k_1-k_{i_2}-\frac{p^{a(1)}}{p^{N}-1}(p^N-1)(k_1-\alpha)=p^{a(1)}\alpha-k_{i_2}\ge0.
	\end{equation}
	Therefore, we have expressed the restriction of $[\mathcal{L}_{Y^1_m}(\{k_i\};\alpha)]$ to $\tilde{Z}$ as a summation of two nef line bundles, so it is also nef.

    \medskip
	Consequently, we have proved the sufficiency part of the nefness criterion for the line bundle $\mathcal{L}_{Y^1_m}(\{k_i\};\alpha)$ over the minimal partial flag space.

	\medskip
	\subsection{ Step 5: The sufficiency part of $\mathcal{L}_{Y^i_m}(\lambda)$ with $i\notin T$}

    In this section, we consider the line bundle 
    $$
    [\mathcal{L}_{Y^i_m}(\lambda)]=[\mathcal{L}_{Y^i_m}(k_{i_1},\dots,k_{i_t};\alpha)]=\sum_{j=1}^tk_j[\omega_{i_j}]-\alpha[\mathcal{F}^i_m],
    $$
    where $Y^i_m$ is a minimal partial flag space with $i\notin T$. We aim to prove the nefness of the above line bundle when the weight $\lambda$ satisfies the condition in Theorem \ref{Nefness criterion for minimal partial flag space}. Our strategy is basically the same as the previous section ---- We construct an auxiliary variety by consecutive blow-ups up to Frobenius twists and make use of the inductive hypothesis at the beginning of this section. Since most ideas are the same, we will omit some details, but make important remarks for the adjustments to the previous subsection.
    
	For simplicity, we assume that $i_1=1\in T$, $i_2=i_1+a(i_1)=1+a(1)$, and $i_1<i<i_2$. 
    Our goal is to show
	\begin{equation}
		\big([\mathcal{L}_{Y^i_m}(\lambda)\cdot C]\big)\ge0
	\end{equation}
	for any complete curve $C$ (or equivalently, its pullback along any proper morphism). 
	

    Let $\pi:Z\to Y^i_m$ be a proper morphism such that $\vec{r}(Z,\pi^\ast\mathcal{F}^i_m)=\vec{r}_0$. The induction hypothesis in subsection 7.1 in particular implies the following assumption.
    
	\noindent\textbf{Induction Hypothesis: }
	
	$\bullet$ The sufficiency part of the nefness criterion for $\mathcal{L}_{Y^i_l}(\lambda)$, where $1\le l< m$, holds for any proper variety over $Y^i_l$ where $\mathcal{F}^i_l$ has slope $\le\vec{r}_0$. More precisely, if $\pi:W\to Y^i_l$ is a proper morphism such that $\vec{r}(W,\pi^\ast\mathcal{F}^i_l)\le \vec{r}_0$, then $\pi^\ast\mathcal{L}_{Y^i_l}(\lambda)$ is nef provided that $\lambda$ satisfies the condition given by Theorem \ref{Nefness criterion for minimal partial flag space}.
	
	$\bullet$ The sufficiency part of the nefness criterion for $\mathcal{L}_{Y^i_l}(\lambda)$, where $m\le l\le n$, holds for any proper variety over $Y^i_l$ where $\mathcal{F}^i_l$ has slope $<\vec{r}_0$. More precisely, if $\pi:W\to Y^i_l$ is a proper morphism such that $\vec{r}(W,\pi^\ast\mathcal{F}^i_l)< \vec{r}_0$, then $\pi^\ast\mathcal{L}_{Y^i_l}(\lambda)$ is nef provided that $\lambda$ satisfies the condition given by Theorem \ref{Nefness criterion for minimal partial flag space}.

    As before, we again distinguish two cases:
    \begin{itemize}[itemsep=0pt,topsep=2pt, parsep=0pt]
        \item Either $r_{\textnormal{tot}}(W,\pi^\ast \mathcal{F}^i_m)>n_1+\cdots+n_N$,
        \item or $r_{\textnormal{tot}}(W,\pi^\ast\mathcal{F}^i_m)=n_1+\cdots+n_N$.
    \end{itemize}
    
    In the first case, we can construct the variety $\tilde{Z}$ in the same manner as Construction \ref{Construction when flagged at essential set}, but with the starting point of each step at $i$. At the $\tilde{\tau}_i$-component, we obtain a sequence of integers 
    $$
    m^0_i<m^1_i<\cdots<m^{l_0}_i
    $$
    together with a sequence
    \begin{equation}
        \mathcal{F}^i_{m^0_i}\subseteq \mathcal{F}^i_{m^1_i}\subseteq \cdots\subseteq \mathcal{F}^i_{m^{l_0}_i}
    \end{equation}
     of subbundles of $\H(\mathcal{A}/X)_i$. As in Lemma \ref{The auxiliary variety is proper surjective}, the forgetful map $\tilde{Z}\to Z$ remains proper and surjective. Moreover, our construction of "strata Hasse invariants" ${h}^l_j$ and $h^l$ for $1\le l\le l_0-1,1\le j\le N$ on $\tilde{Z}$ remains valid (except that we always start from the $\tilde\tau_i$-th component).

    Let $C$ be a complete curve on $\tilde{Z}$. We continue to let $s$ denote the smallest positive integer such that $C$ is contained in the zero locus of $h^s$. (As before, if $C$ is not contained in any of the $Z(h^k)$'s, we simply set $s=l_0$.)
     
    We take
	\begin{equation}
		\begin{aligned}
			[h]&=[h^{s-1}]+\frac{p^N+1}{p^N}[h^{s-2}]+\cdots+\frac{p^{(s-2)N}+\cdots+1}{p^{(s-2)N}}[h^1]\\
			&=p^N[\mathcal{F}^i_{m^s_i}]-\frac{p^{sN}-1}{p^{(s-2)N}(p^N-1)}[\mathcal{F}^i_{m^1_i}]+\frac{p^{(s-1)N}-1}{p^{(s-2)N}(p^N-1)}[\mathcal{F}^i_{m^0_i}].
		\end{aligned}
	\end{equation}
	Recall that $m^1_i=m$. We have the following linear combination
	\begin{equation}
		\begin{aligned}
			&\quad\ [\mathcal{L}_{Y^i_m}(k_{i_1},\dots,k_{i_t};\alpha)]\\
			&=\frac{p^{(s-2)N}(p^N-1)}{p^{sN}-1}\alpha[h] -\frac{p^{(s-1)N}(p^N-1)}{p^{sN}-1}\alpha[\mathcal{F}^i_{m^s_i}] -\frac{p^{(s-1)N}-1}{p^{sN}-1}\alpha[\mathcal{F}^i_{m^0_i}]+\sum_{j=1}^tk_{i_j}[\omega_{i_j}]\\
			&=\beta[h]+ \bigg( -\frac{p^{(s-1)N}(p^N-1)}{p^{sN}-1}\alpha[\mathcal{F}^i_{m^s_1}] +\sum_{j=1}^tuk_{i_j}[\omega_{i_j}]\bigg) +\bigg(-\frac{p^{(s-1)N}-1}{p^{sN}-1}\alpha[\mathcal{F}^i_{m^0_1}] +\sum_{j=1}^t(1-u)k_{i_j}[\omega_{i_j}]\bigg)\\
			&=\beta[h]+[\mathcal{L}_{Y^i_{m^s_i}}(uk_{i_1},\dots,uk_{i_t};\alpha^s_i)]+[\mathcal{L}_{Y^i_{m^0_i}}((1-u)k_{i_1},\dots,(1-u)k_{i_t};\alpha^0_i)].
		\end{aligned}
	\end{equation}
    Here, since $i\notin T$, it would be better to keep track of the set of essential places. So we choose to only record these $\omega_{i_j}$'s in our formula. The key is again to prove the following proposition
    \begin{prop}
        If we are not in the case $t=1,m_1=m=n-1$ or $t=1,m_1=m=1$, then there exists a suitable choice of real numbers $\beta\ge0$ and $u\in [0,1]$ fitting into the above equation, such that the weights for the line bundles $[\mathcal{L}_{Y^i_{m^s_i}}(\{uk_j\};\alpha^s_1)]$ and $[\mathcal{L}_{Y^i_{m^0_1}}(\{(1-u)k_j\};\alpha^0_1)]$ satisfy their respective nefness conditions given by Theorem \ref{Nefness criterion for minimal partial flag space}.
    \end{prop}
    \noindent\emph{Proof.} The number $\beta$ is by definition nonnegative. The choice of $t$ is again based on a case-by-case calculation. Since the ideas of our computations are similar as before, we explain the ideas of our proof in the following.
    
	\noindent\emph{Case 1:} If $t\ge2$, then $\alpha$ satisfies
	\begin{equation}
		\alpha\ge0,\qquad p^{a(1)}\big(k_1-\frac{1}{p^{i-1}}\alpha\big)\ge k_2.
	\end{equation}
	We hope to choose $0\le u\le 1$ such that
	\begin{equation}
		\begin{aligned}
			&p^{a(1)}\big(uk_1-\frac{1}{p^{i-1}}\cdot\frac{p^{(s-1)N}(p^N-1)}{p^{sN}-1}\alpha\big)\ge uk_2,\\
			&p^{a(1)}\big((1-u)k_1-\frac{1}{p^{i-1}}\cdot\frac{p^{(s-1)N}-1}{p^{sN}-1}\alpha\big)\ge (1-u)k_2.
		\end{aligned}
	\end{equation}
	A direct computation of the summation of two equations, similar to the previous step, implies the existence of $u$.

	\medskip
	\noindent\emph{Case 2:} If $t=1$ but $m_1\notin\{1,n-1\}$, then $a(1)=N$, and the arguments are the same as above, with all $k_2$'s replaced by $k_1$.

	\medskip
	\noindent\emph{Case 3:} If $t=1$, $m_1=n-1$, and $m\ne n-1$, then $a(1)=N$ and $\alpha$ satisfies
	\begin{equation}
		\alpha\ge0,\qquad \big(p^{N(n-m)}-1\big)(k_1-\frac{1}{p^{i-1}}\alpha)\ge \big(p^{N(n-m-1)}-1\big)k_1.
	\end{equation}
	Recall that $\delta=r_{\textnormal{tot}}-n_1=r_{\textnormal{tot}}-1$, and that $m^0_i=m-\delta$, $m^s_i=m+(s-1)\delta$ by construction. We hope to choose $0\le u\le 1$ such that
	\begin{equation}
		\begin{aligned}
			&(p^{N(n-(m+(s-1)\delta))}-1)\big(uk_1-\frac{1}{p^{i-1}}\cdot\frac{p^{(s-1)N}(p^N-1)}{p^{sN}-1}\alpha\big)\ge	(p^{N(n-(m+(s-1)\delta)-1)}-1)uk_1,\\
			&(p^{N(n-(m-\delta))}-1)\big((1-u)k_1-\frac{1}{p^{i-1}}\cdot\frac{p^{(s-1)N}-1}{p^{sN}-1}\alpha\big)\ge (p^{N(n-(m-\delta)-1)}-1)(1-u)k_1.
		\end{aligned}
	\end{equation}
	The linear algebra argument then becomes essentially the same as in the previous subsection. We omit the detailed computations here.

	\medskip
	\noindent\emph{Case 4:} If $\# T=1$, $m_1=1$, and $m\ne 1$, this follows dually from the above computation.

    $\hfill\square$

    \begin{prop}
        If $t=1$ and $m_1=m=n-1$ or $1$, then the sufficiency part of the nefness criterion for $\mathcal{L}_{Y^i_m}(k_1;\alpha)$ holds.
    \end{prop}
	
	\medskip
    \noindent\emph{Proof.} The nefness criterion in this case can be proved directly.
    
	\noindent\emph{Case 5:} If $\# T=1$, $m_1=n-1$ and $m=n-1$. Then $a(1)=N$ and $\alpha$ satisfies
	\begin{equation}
		\alpha\ge 0,\qquad k_1-\frac{1}{p^{i-1}}\alpha\ge0.
	\end{equation}
	The condition differs from those in the previous cases. We prove that $\mathcal{L}_{Y^i_{n-1}}(k_1,\dots,k_t;\alpha)$ is nef under the above condition.
    
	Consider the stratum $\tilde{Z}$ defined by
	\begin{equation}
		\tilde{Z}:=[\mathcal{F}^1_{n-2}\subseteq V^{-(i-1)}(\mathcal{F}_{n-1}^{i,(p^{i-1})})]\subseteq Z\times_{Y^i_{n-1}}Y^i_{n-1}\times_X Y^1_{n-2}.
	\end{equation}
	The forgetful map $\tilde{Z}\to Y^i_{n-1}$ is proper and surjective, and in fact, it is a blow-up, since generically $V^{-(i-1)}(\mathcal{F}_{n-1}^{i,(p^{i-1})})$ intersects $\omega_{1}$ in a vector bundle of rank $n-2$.
	
	There exists a strata Hasse invariant on $\tilde{Z}$ given by
	\begin{equation}
		h:\omega_{1}/\mathcal{F}^1_{n-2}\longrightarrow \H(\mathcal{A}/X)_1/V^{-1}_{es}(\mathcal{F}_{n-1}^{i,(p^{i-1})})\stackrel{\simeq}{\longrightarrow} (\H(\mathcal{A}/X)_i/\mathcal{F}^i_{n-1})^{\otimes p^{i-1}}.
	\end{equation}
	The second map is induced by the Verschiebung map. The zero locus of $h$ is the closed stratum defined by $\omega_{1}=V^{-(i-1)}(\mathcal{F}_{n-1}^{i,(p^{i-1})})$. We have
	\begin{equation}
		[h]=-p^{i-1}[\mathcal{F}^i_{n-1}]-[\omega_1/\mathcal{F}^1_{n-2}]
	\end{equation}
	in $\textnormal{(Pic }\tilde{Z}\textnormal{)}_{\mathbb{Q}}$. We consider the following linear combination
	\begin{equation}
		\begin{aligned}
			[\mathcal{L}_{Y^i_{n-1}}(k_1;\alpha)]&=k_1[\omega_1]-\alpha[\mathcal{F}^i_{n-1}]\\
			&=\frac{\alpha}{p^{i-1}}[h]+k_1[\omega_1]+\frac{1}{p^{i-1}}\alpha[\omega_{1}]-\frac{1}{p^{i-1}}\alpha[\mathcal{F}^1_{n-2}]\\
			&=\frac{\alpha}{p^{i-1}}[h]+\big(k_1+\frac{1}{p^{i-1}}\alpha\big)[\omega_{1}/\mathcal{F}^1_{n-2}]+k_1[\mathcal{F}^1_{n-2}]\\
			&=\frac{\alpha}{p^{i-1}}[h]+[\mathcal{L}_{Y^1_{n-2}}(k_1+\frac{1}{p^{i-1}}\alpha;k_1)].
		\end{aligned}
	\end{equation}
	We can view the second factor as the pullback of a line bundle from the minimal partial flag space $Y^1_{n-2}$, and since $Y^1_m$ is flagged at the essential set $T=\{1\}$, we have already proved the sufficiency part of the corresponding nefness criterion in the previous step. It is straightforward to check that
	\begin{equation}
		\begin{aligned}
			k_1+\frac{1}{p^{i-1}}\alpha\ge k_1,\qquad (p^{2N}-1)k_1\ge (p^N-1)(k_1+\frac{1}{p^{i-1}}\alpha)
		\end{aligned}
	\end{equation}
	holds. So $[[\mathcal{L}_{Y^1_{n-2}}(k_1+\frac{1}{p^{i-1}}\alpha,k_1)]$ is a nef class. For a test curve $C\in \tilde{Z}$, we have
	\begin{equation}
		\big([\mathcal{L}_{Y^i_{n-1}}(k_1;\alpha)]\cdot C\big)\ge0
	\end{equation}
	if $C\nsubseteq Z(h)$. If $C\subseteq Z(h)$, then the moduli interpretation of $h$ implies that
	\begin{equation}
		[\omega_1]=[V^{-(i-1)}(\mathcal{F}_{n-1}^{i,(p^{i-1})})]=p^{i-1}[\mathcal{F}_{n-1}^i]
	\end{equation}
	on $\textnormal{Pic (}Z(h)\textnormal{)}_{\mathbb Q}$. Therefore,
	\begin{equation}
		\big([\mathcal{L}_{Y^i_{n-1}}(k_1;\alpha)]\cdot C\big)= \big([\mathcal{L}_{Y^i_{n-1}}(k_1;\alpha)]|_{Z(h)}\cdot C\big)=\big((k_1-\frac{1}{p^{i-1}}\alpha)[\omega_1]\cdot C\big)\ge0.
	\end{equation}
	This completes the proof of the nefness criterion for $\mathcal{L}_{Y_{n-1}^i}(k_1;\alpha)$. 

	\medskip
	\noindent\emph{Case 6:} If $\# T=1$, $m_1=1$, and $m=1$, this follows dually from the above proof.

    $\hfill\square$


    It remains to deal with the case when the pullback of $\mathcal{F}^i_m$ to $Z$ has slope $r_{\textnormal{tot}}(Z,\mathcal{F}^i_m)=n_1+\cdots+n_N$. 

    \begin{prop}\label{Nefness when not flagged at essential set but with base slope}
        If $r_{\textnormal{tot}}(Z,\mathcal{F}^i_m)=n_1+\cdots+n_N$, then $\mathcal{L}_{Y^i_m}(\{k_j\};\alpha)$ is nef over $Z$ as long as the coefficients satisfy the condition of Theorem \ref{Nefness criterion for minimal partial flag space}.
    \end{prop}
	\medskip
    \noindent\emph{Proof.} Similar to the proof in the previous subsection, we apply the induction hypothesis and use the geometric description Theorem \ref{Construction via F_{es},V_{es}-chain}.
    
	\noindent\emph{Case 1:} Either $t=\# T\ge 2$, or $t=1$ but $m_1\notin\{1,n-1\}$.
	
	We need to show that the pullback of $$[\mathcal{L}_{Y^i_m}(k_{i_1},\dots,k_{i_t};\alpha)]=\sum_{j=1}^tk_{i_j}[\omega_{i_j}]-\alpha[\mathcal{F}^i_m]$$ to $Z$ is nef under the condition: If $t\ge2$,
	\begin{equation}
		\begin{aligned}
			&\alpha\ge0,\\
			&p^{a(i_j)}k_{i_j}\ge k_{i_{j+1}},\qquad\textnormal{for $2\le j\le t$},\\
			&p^{a(i_1)}\big(k_{i_1}-\frac{1}{p^{i-1}}\alpha\big)\ge k_{i_2}.
		\end{aligned}
	\end{equation}
	If $t=1$:
	\begin{equation}
		\alpha\ge0,\qquad p^{a(i_1)}\big(k_1-\frac{1}{p^{i-1}}\alpha\big)\ge k_{1}.
	\end{equation}
	
	Let ${u}_i=m$, and for $1\le j\le N-1$, we define inductively
	\begin{equation}
		{u}_{i-j}={u}_{i-j+1}+n_{i-j+1}-r_j.
	\end{equation}
	Here, we view the subscripts of $u$ as modulo $N$. Our choice here is consistent with the previous subsection. Consider the closed subscheme $\tilde{Z}$ defined by the condition
	\begin{equation}
		\tilde{Z}:=[\mathcal{F}^j_{u_j}\subseteq V^{-1}(\mathcal{F}^{j+1,(p)}_{u_{j+1}}),\textnormal{ for all $1\le j\le N$})]
	\end{equation}
	of the variety
	\begin{equation}
        Z\times_{Y^i_m} Y^1_{u_1}\times_X\cdots\times_X Y^N_{u_N}.
	\end{equation}
	The projection map $\tilde{Z}\to Z$ is proper and surjective, for the same reason as before. Thus, it suffices to prove that the pullback of $\mathcal{L}_{Y^i_m}(k_1,\dots,k_t;\alpha)$ to $\tilde{Z}$ is nef. By definition, the sequence
	\begin{equation}
		\mathcal{F}^1_{u_1},\dots, \mathcal{F}^N_{u_N}
	\end{equation}
	forms an $F,V$-chain. 

	By Proposition \ref{Construction via F_{es},V_{es}-chain}, we obtain a proper surjective morphism (up to Frobenius twists) $\tilde{Z}\to X'$, where $X'$ is an auxiliary unitary Shimura variety of signature $(m'_i,n'_i)$. Let $T'$ denote the essential set of $X'$.
	
	\begin{lem}
		Let $t'=\# T'$. Then $t'\le t+1$.
	\end{lem}
	\noindent\emph{Proof:} For $1<j<N$, if $m_{i-j}\in\{0,n\}$, then $r_j=0$ and $\mathcal{F}_{{u}_{i-j}}^{i-j}=V^{-1}(\mathcal{F}^{i-j+1,(p)}_{{u}_{i-j+1}})$. This follows directly from the definition of slopes and ${u}_{i-j}$. So for $s\ne i$, if $m_s=0$ or $n$, then $m'_{s+1}=n$. Since
	\begin{equation}
	    T=\{i_1,i_2,\dots,i_t\}.
	\end{equation}
	This implies that
    \begin{equation}
        T'\subseteq \{i_1+1,i+1,i_2+1,\dots, i_t+1\}.
    \end{equation}
	Hence $t'\le t+1$. $\hfill\square$

    Back to our proof of the proposition. The above lemma implies that we could apply the induction hypothesis at the beginning of this subsection. We know that $\mathcal{L}_{X'}(k'_1,\dots,k'_{t'})$ is nef as long as the coefficients satisfy the corresponding condition in Theorem \ref{Nefness criterion for Shimura variety}. Consider the class (recall that $u_i=m$)
	\begin{equation}
		\begin{aligned}
			[\mathcal{L}_{X'}(\lambda')]:&=\sum_{j=0}^{N-1}p^j[\omega_{i+j+1}']=\sum_{j=0}^{N-1}p^j\cdot\frac{1}{p} \big[\frac{V\mathcal{F}^{i+j}_{u_{i+j}}}{p\mathcal{F}^{i+j+1,(p)}_{u_{i+j+1}}}\big]\\
			&=\sum_{j=0}^{N-1}p^j\cdot\big(\frac{1}{p}[\mathcal{F}^{i+j}_{u_{i+j}}]+[\omega_{i+j+1}/\mathcal{F}^{i+j+1}_{u_{i+j+1}}]\big)\\
			&=\frac{1}{p}[\mathcal{F}^{i}_{m}]+\sum_{j=0}^{N-2}p^{j}[\omega_{i+j+1}]+p^{N-1}[\omega_i/\mathcal{F}^i_{m}].
		\end{aligned}
	\end{equation}
	As before, some $\omega_s$'s might have rank $0$ or $n$. It is non-problematic to keep them in the formula and put some coefficients $k_s$ before them. We omit the computation for the class $[\omega_s']$ via Dieudonn\'e theory.
	
	Now (recall that $i_1=1$)
	\begin{equation}
		\begin{aligned}
			[\mathcal{L}_{Y^i_m}(k_{i_1},\dots,k_{i_t};\alpha)] &=\sum_{j=1}^{t} k_{i_j}[\omega_{i_j}]-\alpha[\mathcal{F}^i_m]\\
			&=\frac{\alpha}{p^{N-1}-\frac{1}{p}}[\mathcal{L}_{X'}(\lambda')] +\sum_{j=2}^{t}\big(k_{i_j}-\frac{\alpha}{p^{N-1}-\frac{1}{p}}\cdot p^{i_j-i-1}\big)[\omega_{i_j}] +\big(k_1-\frac{\alpha}{p^{N-1}-\frac{1}{p}}\cdot p^{N-i}\big)[\omega_1]\\
			&=\frac{\alpha}{p^{N-1}-\frac{1}{p}}[\mathcal{L}_{X'}(\lambda')] +[\mathcal{L}_X(\tilde{\lambda})].
		\end{aligned}
	\end{equation}
	Here $\tilde{\lambda}$ corresponds to the tuple
	\begin{equation}
		(\tilde{k}_{i_1},\dots,\tilde{k}_{i_t})=(k_1-\frac{\alpha p^{N-i}}{p^{N-1}-\frac{1}{p}} ,k_{i_2}-\frac{\alpha p^{i_2-i-1}}{p^{N-1}-\frac{1}{p}},\dots,k_{i_t}-\frac{\alpha p^{i_t-i-1}}{p^{N-1}-\frac{1}{p}}).
	\end{equation}
	It suffices to check that the latter is a nef class. For $2\le j\le t$, we see that
	\begin{equation}
		p^{a(i_j)}\tilde{k}_{i_j}\ge \tilde{k}_{i_{j+1}}.
	\end{equation}
	Moreover (recall that $i_2-1=a(1)$),
	\begin{equation}
		\begin{aligned}
			&p^{a(1)}\big(k_1-\frac{\alpha}{p^{N-1}-\frac{1}{p}}\cdot p^{N-i}\big)\ge k_{i_2}-\frac{\alpha}{p^{N-1}-\frac{1}{p}}\cdot p^{i_2-i-1}\\
			\Leftrightarrow\qquad &p^{a(1)}k_1-\frac{p^{N+a(1)-i}-p^{i_2-i-1}}{p^{N-1}-\frac{1}{p}}\cdot \alpha\ge k_{i_2}\\
			\Leftrightarrow\qquad &p^{a(1)}k_1-p^{a(1)+1-i}\alpha\ge k_{i_2}\\
			\Leftrightarrow\qquad &p^{a(1)}\big(k_1-\frac{1}{p^{i-1}}\alpha\big)\ge k_{i_2},
		\end{aligned}
	\end{equation}
	which holds by our assumption.

	\medskip
	\noindent\emph{Case 2:} We have $t=\# T=1$, $m_1=n-1$, but $m\ne n-1$.
	
	In this case, we need to show that the pullback of $\mathcal{L}_{Y^i_m}(k_1;\alpha)=k_1[\omega_1]-\alpha[\mathcal{F}^i_m]$ to $Z$ is nef under the condition
	\begin{equation}
		\alpha\ge0,\qquad (p^{N(n-m)}-1)(k_1-\frac{1}{p^{i-1}}\alpha)\ge (p^{N(n-m-1)}-1)k_1.
	\end{equation}
	We could again construct the auxiliary variety $\tilde{Z}$ and the $F,V$-chain
	\begin{equation}
		\mathcal{F}^1_{u_1},\dots,\mathcal{F}^N_{u_N}.
	\end{equation}
	Note that $r_{\textnormal{tot}}=\sum_{j=1}^N n_j=1$ in this case. Thus, all the entries in the tuple $\vec{r}(\tilde{Z},\mathcal{F}^i_m)=(r_1,\dots,r_N)$ are 0, except a single 1. By Lemma \ref{Properties of slope}, $r_j\le n_{i-j}$ for $1\le j<N$. The only possible nonzero terms are $r_{i-1}$ and $r_N$.
	
	If $r_{i-1}=0$ and $r_N=1$. Then $u_1=u_2=\cdots=u_i=m$, and the inclusion conditions 
    \begin{equation}
        \mathcal{F}^{i-v}_{u_{i-v}}\subseteq V^{-1}(\mathcal{F}^{i-v+1,(p)}_{u_{i-v+1}})
    \end{equation}
    for $1\le v\le N-1$ are, in fact, equalities. In particular, $u_1=m$ and
	\begin{equation}
		\mathcal{F}^1_{u_1}=V^{-(i-1)}(\mathcal{F}_{m}^{i,(p^{i-1})}).
	\end{equation}
	In $(\textnormal{Pic $\tilde{Z}$})_{\mathbb{Q}}$, we deduce that
	\begin{equation}
		[\mathcal{F}^1_m]={p^{i-1}}[\mathcal{F}^i_m].
	\end{equation}
	So we have
	\begin{equation}
		\begin{aligned}
			[\mathcal{L}_{Y^i_m}(k_1;\alpha)]&=k_1[\omega_1]-\alpha[\mathcal{F}^i_m]=k_1[\omega_1]-\frac{1}{p^{i-1}}\alpha[\mathcal{F}^1_m]\\
			&=k_1[\omega_1/\mathcal{F}^1_m]+\big(k_1-\frac{1}{p^{i-1}}\alpha\big)[\mathcal{F}^1_m]\\
			&=[\mathcal{L}_{Y^1_m}(k_1;k_1-\frac{1}{p^{i-1}}\alpha)].
		\end{aligned}
	\end{equation}
	Note that $\tilde{Z}$ admits a proper morphism to $Y^1_m$, so the last equality means that we can view the above class as the pullback of an automorphic line bundle from $Y^1_m$. Note that $Y^1_m$ is flagged at the essential set $T=\{1\}$, and
	\begin{equation}
		(p^{N(n-m)}-1)(k_1-\frac{1}{p^{i-1}}\alpha)\ge (p^{N(n-m-1)}-1)k_1.
	\end{equation}
	Since we have already proved the sufficiency part of the nefness criterion for $\mathcal{L}_{Y^i_m}(\mu)$, we deduce that $\mathcal{L}_{Y^i_m}(k_1;\alpha)|_Z$ is nef.
	
	If $r_{i-1}=1$ and $r_N=0$. Then $u_1=u_2=\cdots=u_{i-1}=m-1$ and $u_i=\cdots=u_N=m$. Hence,
	\begin{equation}
		\mathcal{F}^i_m=V^{-(N+1-i)}(\mathcal{F}_{m-1}^{1,(p^{N+1-i})}).
	\end{equation}
	By Lemma \ref{Computing Picard Classes}, we have
	\begin{equation}
		[\mathcal{F}^i_m]=-p^{N+1-i}[\omega_1/\mathcal{F}^1_{m-1}]
	\end{equation}
	in $(\textnormal{Pic $\tilde{Z}$})_{\mathbb{Q}}$. So
	\begin{equation}
		\begin{aligned}
			[\mathcal{L}_{Y^i_m}(k_1;\alpha)]&=k_1[\omega_1]+p^{N+1-i}\alpha[\omega_1/\mathcal{F}^1_{m-1}]\\
			&=\big(k_1+p^{N+1-i}\alpha\big)[\omega_1/\mathcal{F}^1_{m-1}]+k_1[\mathcal{F}^1_{m-1}]\\
			&=[\mathcal{L}_{Y^1_{m-1}}(k_1+p^{N+1-i}\alpha;k_1)].
		\end{aligned}
	\end{equation}
	It suffices to check
	\begin{equation}
		(p^{N(n-m+1)}-1)k_1\ge (p^{N(n-m)}-1)(k_1+p^{N+1-i}\alpha).
	\end{equation}
	This is straightforward. (The case $m=1$ is a little bit different but much easier. In this case, the condition is $[\mathcal{F}^i_m]=-p^{N+1-i}[\omega_1]$. We have
    \begin{equation}
        [\mathcal{L}_{Y^i_m}(k_1;\alpha)]=(k_1+p^{N+1-i}\alpha)[\omega_1].
    \end{equation}
    The nefness follows directly from the nefness of the Hodge line bundle $[\omega_1]$.) This completes the proof in Case 2.
	
	\medskip
	\noindent\emph{Case 2':} If $r=\# T=1$ and $m_1=1$. This follows from a dual argument.
	
	\medskip
	By the above case-by-case discussion, we have completed the proof of Proposition \ref{Nefness when not flagged at essential set but with base slope}. $\hfill\square$

	\medskip
	Combining all the discussions above, we have shown the nefness of the automorphic line bundles $\mathcal{L}_{Y^i_m}(\{k_j\};\alpha)$ over minimal partial flag spaces, under the numerical assumption of Theorem \ref{Nefness criterion for minimal partial flag space}. 

	\medskip
    Combining the arguments in Section 7, we have completed the proof of the ampleness criterion for automorphic line bundle over unitary Shimura varieties, i.e., Theorem \ref{The ampleness criterion for the flag space, case 1}, \ref{The ampleness criterion for the flag space, case 2}, \ref{The ampleness criterion for the flag space, case 2'}, \ref{The ampleness criterion for the flag space, case 3}. $\hfill\square$

	\medskip
	\noindent \textbf{Acknowledgements:} I want to thank my advisor, Liang Xiao, for suggesting working on this problem and for his constant encouragement. His insightful suggestions have been valuable throughout this research. I also want to thank Matthew Emerton, Yichao Tian, Zhiyu Tian, Kai Wen Lan, Fred Diamond, Payman Kassaei, George Boxer, Rong Zhou for various helpful conversations and their thoughtful perspectives on various aspects of this ampleness criterion. Additionally, I want to thank Wushi Goldring and Jean-Stefan Koskivirta for sharing their understanding of zip cones, strata Hasse invariants, and their work on general ampleness criterion.

	\medskip

	\bibliographystyle{plain}
	\bibliography{Reference}

@article{AlexandreVanishing,
  title={Vanishing results for the coherent cohomology of automorphic vector bundles over the {S}iegel variety in positive characteristic},
  author={Alexandre, Thibault},
  journal={Algebra \& Number Theory},
  volume={19},
  number={1},
  pages={143--193},
  year={2024},
  publisher={Mathematical Sciences Publishers}
}

@article{Andreatta-Goren,
  title={Hilbert modular varieties of low dimension},
  author={Andreatta, Fabrizio and Goren, Eyal Z},
  journal={Geometric aspects of Dwork theory},
  volume={1},
  pages={113--175},
  year={2004},
  publisher={Walter de Gruyter}
}

@article{Buzzard-Diamond-Jarvis,
  title={On Serre's conjecture for mod $\ell$ Galois representations over totally real fields},
  author={Buzzard, Kevin and Diamond, Fred and Jarvis, Frazer},
  year={2010}
}

@incollection{BrionFlagVarieties,
  title={Lectures on the geometry of flag varieties},
  author={Brion, Michel},
  booktitle={Topics in Cohomological Studies of Algebraic Varieties: Impanga Lecture Notes},
  pages={33--85},
  year={2005},
  publisher={Springer}
}

@article{BGKS,
  title={Ampleness of automorphic bundles on zip-schemes},
  author={Y. Brunebarbe and J.-S. Koskivirta and W. Goldring and B. Stroh},
  journal={in preparation}
}

@article{CalegariGeraghty,
  title={Modularity lifting beyond the {T}aylor--{W}iles method},
  author={Calegari, Frank and Geraghty, David},
  journal={Inventiones mathematicae},
  volume={211},
  number={1},
  pages={297--433},
  year={2018},
  publisher={Springer}
}

@article{Caraiani-Tamiozzo,
  title={On the {\'e}tale cohomology of {H}ilbert modular varieties with torsion coefficients},
  author={Caraiani, Ana and Tamiozzo, Matteo},
  journal={Compositio Mathematica},
  volume={159},
  number={11},
  pages={2279--2325},
  year={2023},
  publisher={London Mathematical Society}
}

@article{Caraiani-ScholzeGenericCohomology,
  title={On the generic part of the cohomology of non-compact unitary Shimura varieties},
  author={Caraiani, Ana and Scholze, Peter},
  journal={Annals of Mathematics},
  volume={199},
  number={2},
  pages={483--590},
  year={2024},
  publisher={Department of Mathematics, Princeton University Princeton, New Jersey, USA}
}

@article{Caraiani-ScholzeCompactGenericPart,
  title={On the generic part of the cohomology of compact unitary Shimura varieties},
  author={Caraiani, Ana and Scholze, Peter},
  journal={Annals of Mathematics},
  volume={186},
  number={3},
  pages={649--766},
  year={2017},
  publisher={Department of Mathematics, Princeton University Princeton, New Jersey, USA}
}

@article{CaraianiShin,
  title={Recent progress on {L}anglands reciprocity for {GLn}: {S}himura varieties and beyond},
  author={Caraiani, Ana and Shin, Sug Woo},
  journal={preprint},
  year={2023}
}

@article{Diamond-Sasaki,
  title={A {Serre} weight conjecture for geometric {Hilbert} modular forms in characteristic $p$},
  author={Diamond, Fred and Sasaki, Shu},
  journal={arXiv preprint arXiv:1712.03775},
  year={2017}
}

@article{Emerton-Reduzzi-Xiao,
  title={Galois representations and torsion in the coherent cohomology of Hilbert modular varieties},
  author={Emerton, Matthew and Reduzzi, Davide A and Xiao, Liang},
  journal={Journal f{\"u}r die reine und angewandte Mathematik (Crelles Journal)},
  volume={2017},
  number={726},
  pages={93--127},
  year={2017},
  publisher={De Gruyter}
}

@book{Fulton-HarrisBook,
  title={Representation theory: a first course},
  author={Fulton, William and Harris, Joe},
  volume={129},
  year={2013},
  publisher={Springer Science \& Business Media}
}

@article{GeorgeBoxerThesis,
  title={Torsion in the coherent cohomology of $S$himura varieties and $G$alois representations},
  author={Boxer, George},
  journal={Ph.D. thesis},
  publisher={Harvard University, Cambridge, MA},
  year={2015}
}

@article{GK-strata,
  title={Strata {Hasse} invariants, {Hecke} algebras and {Galois} representations},
  author={Goldring, Wushi and Koskivirta, Jean-Stefan},
  journal={Inventiones mathematicae},
  volume={217},
  pages={887--984},
  year={2019},
  publisher={Springer}
}

@article{GK-stratification,
  title={Stratifications of flag spaces and functoriality},
  author={Goldring, Wushi and Koskivirta, Jean-Stefan},
  journal={International Mathematics Research Notices},
  volume={2019},
  number={12},
  pages={3646--3682},
  year={2019},
  publisher={Oxford University Press}
}

@article{Goldring-Imai-Koskivirta,
  title={Weights of mod $p$ automorphic forms and partial {H}asse invariants},
  author={Goldring, Wushi and Imai, Naoki and Koskivirta, Jean-Stefan},
  journal={arXiv preprint arXiv:2211.16207},
  year={2022}
}

@book{GrothendieckCrystallineDieudonne,
  title={Groupes de {B}arsotti-{T}ate et cristaux de {D}ieudonn{\'e}},
  author={Grothendieck, Alexandre and Grothendieck, Alexandre and Mathematician, Germany and Grothendieck, Alexandre and Grothendieck, Alexandre},
  volume={45},
  year={1974},
  publisher={Presses de l'Universit{\'e} de Montr{\'e}al}
}

@article{Hamann-LeeVanishing,
  title={Torsion vanishing for some Shimura varieties},
  author={Hamann, Linus and Lee, Si Ying},
  journal={arXiv preprint arXiv:2309.08705},
  year={2023}
}

@inproceedings{HarrisSurvey,
  title={Automorphic {G}alois representations and the cohomology of {S}himura varieties},
  author={Harris, Michael},
  booktitle={Proceedings of the International Congress of Mathematicians, Seoul},
  volume={2},
  pages={367--389},
  year={2014}
}

@article{HarrisTWMethod1,
  title={The {Taylor-Wiles} method for coherent cohomology.},
  author={Harris, Michael},
  journal={Journal f{\"u}r die Reine und Angewandte Mathematik},
  volume={2013},
  number={679},
  year={2013}
}

@article{AtanasovHarrisTWMethod2,
  title={{The Taylor-Wiles method for coherent cohomology, II}},
  author={Atanasov, Stanislav and Harris, Michael},
  journal={American Journal of Mathematics},
  volume={147},
  number={5},
  pages={1383--1432},
  year={2025},
  publisher={Johns Hopkins University Press}
}

@article{Helm,
  title={A geometric {Jacquet-Langlands} correspondence for {$U(2)$} Shimura varieties},
  author={Helm, David},
  journal={Israel Journal of Mathematics},
  volume={187},
  pages={37--80},
  year={2012},
  publisher={Springer}
}

@article{Helm-Tian-Xiao,
  title={Tate cycles on some unitary {Shimura} varieties mod $p$},
  author={Helm, David and Tian, Yichao and Xiao, Liang},
  journal={Algebra \& Number Theory},
  volume={11},
  number={10},
  pages={2213--2288},
  year={2017},
  publisher={Mathematical Sciences Publishers}
}

@article{IK-partial,
  title={Partial {Hasse} invariants for {Shimura} varieties of {Hodge}-type},
  author={Imai, Naoki and Koskivirta, Jean-Stefan},
  journal={arXiv preprint arXiv:2109.11117},
  year={2021}
}

@book{JantzenRepresentationBook,
  title={Representations of algebraic groups},
  author={Jantzen, Jens Carsten},
  volume={107},
  year={2003},
  publisher={American Mathematical Soc.}
}

@article{KisinIntegralModel,
  title={Integral models for Shimura varieties of abelian type},
  author={Kisin, Mark},
  journal={Journal of the American Mathematical Society},
  volume={23},
  number={4},
  pages={967--1012},
  year={2010}
}

@article{Koskivirta-Wedhorn,
  title={Generalized $\mu$-ordinary {Hasse} invariants},
  author={Koskivirta, Jean-Stefan and Wedhorn, Torsten},
  journal={Journal of Algebra},
  volume={502},
  pages={98--119},
  year={2018},
  publisher={Elsevier}
}

@article{Kottwitz,
  title={Points on some {Shimura} varieties over finite fields},
  author={Kottwitz, Robert E},
  journal={Journal of the American Mathematical Society},
  volume={5},
  number={2},
  pages={373--444},
  year={1992},
  publisher={JSTOR}
}

@article{KoskivirtaNormalization,
  title={Normalization of closed ekedahl—oort strata},
  author={Koskivirta, Jean-Stefan},
  journal={Canadian Mathematical Bulletin},
  volume={61},
  number={3},
  pages={572--587},
  year={2018},
  publisher={Cambridge University Press}
}

@book{LanCompactification,
  title={Arithmetic compactifications of {PEL-type Shimura} varieties},
  author={Lan, Kai-Wen},
  year={2008},
  publisher={Harvard University}
}

@article{Lan-SuhCompactVanishing,
  title={Vanishing theorems for torsion automorphic sheaves on compact {PEL}-type {S}himura varieties},
  author={Lan, Kai-Wen and Suh, Junecue},
  journal={Duke Mathematical Journal},
  volumn={161. No. 6},
  year={2012}
}

@article{Lan-SuhGeneralVanishing,
  title={Vanishing theorems for torsion automorphic sheaves on general PEL-type Shimura varieties},
  author={Lan, Kai-Wen and Suh, Junecue},
  journal={Advances in Mathematics},
  volume={242},
  pages={228--286},
  year={2013},
  publisher={Elsevier}
}

@book{Lazarsfeld,
  title={Positivity in algebraic geometry I: {Classical} setting: line bundles and linear series},
  author={Lazarsfeld, Robert K},
  volume={48},
  year={2017},
  publisher={Springer}
}

@book{Mazur-Messing,
  title={Universal extensions and one dimensional crystalline cohomology},
  author={Mazur, Barry and Messing, William},
  volume={370},
  year={2006},
  publisher={Springer}
}

@article{MessingCrystals,
  title={The crystals associated to {B}arsotti-{T}ate groups},
  author={Messing, William},
  journal={The crystals associated to Barsotti-Tate groups: with applications to abelian schemes},
  pages={112--149},
  year={2006},
  publisher={Springer}
}

@incollection{MoonenGroupsSchemes,
  title={Group schemes with additional structures and {W}eyl group cosets},
  author={Moonen, Ben},
  booktitle={Moduli of abelian varieties},
  pages={255--298},
  year={2001},
  publisher={Springer}
}

@article{Moonen-Wedhorn,
  title={Discrete invariants of varieties in positive characteristic},
  author={Moonen, Ben and Wedhorn, Torsten},
  journal={International Mathematics Research Notices},
  volume={2004},
  number={72},
  pages={3855--3903},
  year={2004},
  publisher={Hindawi Publishing Corporation}
}

@book{MumfordAV,
  title={Abelian varieties},
  author={Mumford, David},
  volume={13},
  publisher={Tata Institute of Fundamental Research Publications},
  year={2012}
}

@inproceedings{OdaDeRham,
  title={The first de Rham cohomology group and Dieudonn{\'e} modules},
  author={Oda, Tadao},
  booktitle={Annales scientifiques de l'{\'E}cole Normale Sup{\'e}rieure},
  volume={2},
  pages={63--135},
  year={1969}
}

@inproceedings{MadapusiCompactification,
  title={Toroidal compactifications of integral models of Shimura varieties of Hodge type},
  author={Madapusi Pera, Keerthi},
  booktitle={Annales scientifiques de l'{\'E}cole Normale Sup{\'e}rieure},
  volume={52},
  number={2},
  pages={393--514},
  year={2019}
}

@article{Pink-Wedhorn-ZieglerAlgebraicZipData,
  title={Algebraic zip data},
  author={Ziegler, P and Wedhorn, Torsten and Pink, R},
  journal={Documenta Mathematica},
  volume={16},
  year={2011},
  publisher={Documenta Mathematica}
}

@article{Pink-Wedhorn-ZieglerF-zips,
  title={F-zips with additional structure},
  author={Pink, Richard and Wedhorn, Torsten and Ziegler, Paul},
  journal={Pacific Journal of Mathematics},
  volume={274},
  number={1},
  pages={183--236},
  year={2015},
  publisher={Mathematical Sciences Publishers}
}

@inproceedings{Reduzzi-XiaoPartialHasse,
  title={Partial Hasse invariants on splitting models of Hilbert modular varieties},
  author={Reduzzi, Davide A and Xiao, Liang},
  booktitle={Annales scientifiques de l'{\'E}cole Normale Sup{\'e}rieure},
  volume={50},
  pages={579--607},
  year={2017}
}

@article{Tian-Xiao,
  title={On {Goren--Oort} stratification for quaternionic Shimura varieties},
  author={Tian, Yichao and Xiao, Liang},
  journal={Compositio Mathematica},
  volume={152},
  number={10},
  pages={2134--2220},
  year={2016},
  publisher={London Mathematical Society}
}

@article{VasiuIntegralModel,
  title={Integral canonical models of Shimura varieties of preabelian type},
  author={Vasiu, Adrian},
  journal={Asian Journal of Mathematics},
  volume={3},
  pages={401--517},
  year={1999},
  publisher={INTERNATIONAL PRESS}
}

@article{Viehmann-WedhornEONewton,
  title={Ekedahl--Oort and Newton strata for Shimura varieties of PEL type},
  author={Viehmann, Eva and Wedhorn, Torsten},
  journal={Mathematische Annalen},
  volume={356},
  number={4},
  pages={1493--1550},
  year={2013},
  publisher={Springer}
}

@article{Wedhorn-ZieglerTautological,
  title={Tautological rings of {S}himura varieties and cycle classes of {E}kedahl--{O}ort strata},
  author={Wedhorn, Torsten and Ziegler, Paul},
  journal={Algebra \& Number Theory},
  volume={17},
  number={4},
  pages={923--980},
  year={2023},
  publisher={Mathematical Sciences Publishers}
}

@article{Xiao-Zhu,
  title={Cycles on Shimura varieties via geometric Satake},
  author={Xiao, Liang and Zhu, Xinwen},
  journal={arXiv preprint arXiv:1707.05700},
  year={2017}
}

@article{YangAmpleness,
  title={Ampleness of Automorphic Line Bundles on $U(2)$ Shimura Varieties},
  author={Yang, Deding},
  journal={arXiv preprint arXiv:2309.00286},
  year={2023}
}

@unpublished{YangYang,
  author = {Yang, Deding and Yang, Siqi},
  title  = {Local-global compatibility for quaternionic automorphic forms},
  note   = {In progress}
}

@article{ZhangEOStrata,
  title={Ekedahl-Oort strata for good reductions of Shimura varieties of Hodge type},
  author={Zhang, Chao},
  journal={Canadian Journal of Mathematics},
  volume={70},
  number={2},
  pages={451--480},
  year={2018},
  publisher={Cambridge University Press}
}
	\addcontentsline{toc}{section}{References}

\end{document}